\documentclass[10pt]{elsarticle}
\usepackage{lineno,hyperref}
\usepackage{booktabs}
\usepackage{multirow}
\usepackage{amsthm}
\usepackage{graphicx}
\usepackage{subfigure}
\usepackage{float}
\usepackage{bm}
\usepackage{lipsum}
\usepackage{listings} 
\usepackage{color}
\usepackage{amsfonts,enumerate,amsmath,amssymb}
\def\qed{\strut\hfill $\Box$}
\newtheorem{thm}{Theorem}[section]

\newtheorem{prop}[thm]{Proposition}
\newtheorem{lem}[thm]{Lemma}

\newtheorem{rem}[thm]{Remark}

\newcommand{\thmref}[1]{Theorem~{\rm \ref{#1}}}
\newcommand{\lemref}[1]{Lemma~{\rm \ref{#1}}}
\newcommand{\propref}[1]{Proposition~{\rm \ref{#1}}}

\newcommand{\vertiii}[1]{{\left\vert\kern-0.25ex\left\vert\kern-0.25ex\left\vert #1 
		\right\vert\kern-0.25ex\right\vert\kern-0.25ex\right\vert}}
\def\para#1{\vskip .4\baselineskip\noindent{\bf #1}}
\numberwithin{equation}{section}
\begin{document}
	\begin{frontmatter}

		\title{Precise Laplace approximation for mixed rough differential equation}
		
		\author[mymainaddress]{Xiaoyu Yang}
		\ead{yangxiaoyu@yahoo.com}
		
		\author[mymainaddress,mythirdaddress]{Yong Xu\corref{mycorrespondingauthor}}\cortext[mycorrespondingauthor]{Corresponding author}
		\ead{hsux3@nwpu.edu.cn}
		
		\author[mymainaddress]{Bin Pei}
		\ead{binpei@nwpu.edu.cn}

		\address[mymainaddress]{School of Mathematics and Statistics, Northwestern Polytechnical University, Xi'an, 710072, China}
		\address[mythirdaddress]{MIIT Key Laboratory of Dynamics and Control of Complex Systems, Northwestern Polytechnical University, Xi'an, 710072, China}

		\begin{abstract}
			This work focuses on the Laplace approximation for  the   rough differential equation (RDE) driven by mixed rough path $(B^H, W)$ with  $H \in (1/3,1/2)$ as $\varepsilon \to 0$. Firstly, based on geometric rough path lifted from mixed fractional Brownian motion (fBm),  the Schilder-type large deviation principle (LDP) for the law of the first level path of the solution to the RDE is given. {Due to the particularity of mixed rough path, the main difficulty in carrying out the Laplace approximation is  to prove the Hilbert-Schmidt property for the Hessian matrix  of the It\^o map restricted on the Cameron-Martin space of the mixed fBm. To this end, we imbed the Cameron-Martin space into a larger Hilbert space, then the Hessian is computable.  Subsequently,  the probability representation for the Hessian is shown.} Finally,  the  Laplace approximation is constructed, which asserts   the more precise asymptotics in the exponential scale.
			
			\vskip 0.08in
			\noindent{\bf Keywords.}
			 Large deviation principle, Laplace approximation, Mixed rough path, fractional Brownian motion
			\vskip 0.08in	
			\noindent {\bf Mathematics subject classification.} 60F10, 60G22, 60H10.			
		\end{abstract}		
	\end{frontmatter}

	\section{Introduction}\label{sec-1}
	In this paper, we consider the   rough differential equation (RDE) driven by mixed rough paths as follows,
	\begin{eqnarray}\label{1-1}
	dY_{t}^{\varepsilon}=\left[\sigma \left( Y_{t}^{\varepsilon} \right)|\hat\sigma \left( Y_{t}^{\varepsilon} \right)\right] \varepsilon d(B^H,W)_t+\beta \left(\varepsilon,Y_{t}^{\varepsilon} \right) dt, \quad Y_{0}^{\varepsilon}=\text{0,}
	\end{eqnarray}
	where  $\varepsilon>0$ is a small parameter. And $\left[\sigma |\hat\sigma \right]$ denotes the block matrix with $\sigma \in C_{b}^{\infty}\left(\mathbb{R}^{n}, \operatorname{Mat}(n, d_1)\right)$ and $\hat\sigma \in C_{b}^{\infty}\left(\mathbb{R}^{n}, \operatorname{Mat}(n, d_2)\right)$, and $\beta \in C_{b}^{\infty}\left([0,1]\times\mathbb{R}^{n}, \mathbb{R}^{n}\right)$ with $C_{b}^{\infty}$ being
	the set of bounded smooth functions with bounded derivatives. $(B^H,W)_t\in G \Omega_{p}(\mathbb{R}^{d_1+d_2})$, where $G \Omega_{p}(\mathbb{R}^{d_1+d_2})$ is the geometric rough path space, represents the mixed geometric rough path,  which will be introduced in Section 2.  Let $\big(\Omega, \mathcal{F},\left\{\mathcal{F}_{t}\right\}_{t \geq 0}, \mathbb{P}\big)$ be a completely probability space,
	  $(b^H_t)_{t\ge 0}$ is $\mathbb{R}^{d_1}$-valued fractional Brownian motion (fBm)  with Hurst index $H\in(1/3,1/2)$, $(w_t)_{t\ge 0}$  is $\mathbb{R}^{d_2}$-valued Brownian motion (Bm). For each time $t$, we denote by $\mathcal{F}_{t}$ the $\sigma$-field generated by the random variables $\{(b_s^H,w_s), 0\le s \le t\}$ and all $\mathbb{P}$-null sets. The expectation with respect to $\mathbb{P}$ is denoted by $\mathbb{E}$.
	  More details about fBm could see \cite{2008Biagini,2008Mishura}.   
	  
	  Let $Y^{\varepsilon}=\hat{\Phi}_{\varepsilon}(\varepsilon (B^H, W), \lambda): G \Omega_{p}\left(\mathbb{R}^{d_1+d_2+1}\right) \mapsto G \Omega_{p}\left(\mathbb{R}^{n}\right)$ denote the It\^o map corresponding to (\ref{1-1}) with $\lambda_{t}=t$.  The purpose of this paper is to prove the Laplace approximation for $Y^{\varepsilon,1}$ under natural assumptions, the   first level path of the solution map,
	\begin{eqnarray}\label{1-3}
	J(\varepsilon):=\mathbb{E}\left[\exp \left(-\frac{F\big(Y_{t}^{\varepsilon,1}\big)}{\varepsilon^{2}}\right)\right],
	\end{eqnarray} 
	where $F$ is a suitable real-valued bounded continuous function. See  Section 2 for precise assumptions.

	Laplace approximation is devised by Pierre-Simon Laplace in 1977. For SDE theory, the origins can be traced back to the Azencott's work on stochastic Taylor expansion  \cite{1982Azencott}. Based on this, Ben Arous gave some Laplace short time expansion for Wiener functional \cite{1988Arous}.
	Since then, Laplace approximation has been a central topic in the probability field. Watanabe studied the precise asymptotics of the Schilder-type for some classes of generalized Wiener functionals by Malliavin calculus \cite{1987Watanabe,1993Takanobu}.  Besides, Kusuoka and Stroock presented the asymptotic expansion of certain Wiener functionals as the variance of the Wiener goes to 0 \cite{1991Kusuoka}.  Osajima focused on the asymptotic expansion of the density function of Wiener functionals \cite{2008Osajima}. For stochastic partial differential equations (SPDEs), it could refer to  \cite{2000Rovira,2021Friz,2022Friz}.
	
     However, the abovementioned references focused on the Bm. Different from Bm, fBm is self-similar and possesses long-range dependence, which has been applied in some complex systems \cite{2008Biagini}. Since the fBm is neither a semi-martingale nor a Markov process, it can not be solved by conventional stochastic analysis \cite{2008Mishura}. Rough path,  proposed by Terry Lyons in 1998,  considers the path itself and the iterated integral together \cite{1998Lyons}. It does not need martingale integration theory, Markov property, and filtration theory, but focuses on the real analysis. Hence, it is effective to analyze fBm. More importantly, the central result in rough path theory, Lyons' universal limit theorem,  states that the solution map is continuous (and locally Lipschitz) with respect to the topology of geometric rough path space \cite{2002Lyons}. Up to now,  there exist three main formulations to rough path theory, Lyons' original formulation \cite{2002Lyons}, Gubinelli’s controlled path theory \cite{2020Friz}, and Davie’s formulation \cite{2010Friz,2015Bailleul}.
    With the aid of rough path theory, Gaussian processes, including Bm and fBm, can be lifted to geometric rough paths. Therefore, it has been applied to analyze the	LDP for Gaussian rough paths, including the Brownian rough path and fBm \cite{2002Ledoux,2006Millet,2007Friz,2006InahamaKawabi}. Further, starting with Aida \cite{2006Aida}, Inahama and Kawabi studied the Taylor expansion and Laplace approximation for Gaussian rough paths, where the fractional rough paths can be covered \cite{2006Inahama,2007Inahama,2013Inahama}.

	Therefore, we wish to investigate the precise Laplace asymptotics for  the first level path of $Y_t^\varepsilon$. {Firstly, we prove that the mixed fBm can be lifted to the mixed geometric rough path. Then, we give the Schilder-type LDP for the laws of the first level path of the solution to the above RDE (\ref{1-1}), which provides  the asymptotics of $J(\varepsilon)$ with $\varepsilon \to 0$ on logarithmic scale.  Furthermore, it proceeds to show the Laplace approximation, providing more precise asymptotics in the exponential scale. Due to the particularty of mixed rough path, the main difficulty in carrying out the Laplace asymptotics is  to analyze Hilbert-Schmidt property for the Hessian matrix  of the It\^o map restricted on the Cameron-Martin space. To this end, we imbed the Cameron-Martin space into a larger Hilbert space.  We prove the Hilbert-Schmidt property with  an orthonormal basis in this Hilbert space. Then, the  probability representation for the Hessian is analyzed. Finally,  the precise Laplace asymptotics for RDE driven by mixed rough path is constructed.}


	The plan of this paper is as follows. In the next section, we establish notation, give some precise conditions and main result. In Section \ref{sec-3}, we state and prove the  Schilder-type LDP for the laws of the first level path of the solution to the RDE (\ref{1-1}). We prove the Hilbert-Schmidt property and probability representation of the Hessian matrix in Section \ref{sec-4}. Finally,  Section \ref{sec-5} states the main proof.
	Throughout this paper, $c$, $C$, $C_\star$ denote certain positive constants that may vary from line to line. Throughout this paper, we take $t\in [0,1]$. Analogous results hold for any  finite time interval.
	
	\section{Preliminaries and main results}\label{sec-2}
	Before introducing the fBm and rough path, we illustrate some information for spaces.
	Let $B$ be a Banach space with $dimB<\infty$, such as $B=\mathbb{R}^{d_1 +d_2}$ or $B=Mat(n,d_1 +d_2)$. Denote that 
	$$C=C([0,1], B)=\{k:[0,1] \rightarrow B \mid \text { continuous }\},$$
	the space of $B$-valued continuous functions with the usual sup-norm. For $p\ge 1$, 
	$C^{p\mathrm{-var}}$ is the set of $k\in C$ such that 
	$$\|k\|_{p-\mathrm{var}} :=\left|k_{0}\right|+\left(\sup _{\mathcal{P}} \sum_{i=1}^{n}\left|k_{t_{i}}-k_{t_{i-1}}\right|^{p}\right)^{1 / p}<\infty,$$
	where $\mathcal{P}$ runs over all the finite partition of $[0,1]$.

	Denote $W^{\delta,p}$ with $p>1, 0<\delta<1$  the Besov space, for a measurable function $k:[0,1]\to B$,
	\begin{eqnarray*}
		\|k\|_{W^{\delta, p}}=\|k\|_{L^{p}}+\left(\iint_{[0,1]^{2}} \frac{\left|k_{t}-k_{s}\right|^{p}}{|t-s|^{1+\delta p}} d s d t\right)^{1 / p}.
	\end{eqnarray*}
	Refer to \cite{2003Adams}, we can see that $W^{\delta,p}$ is given by the interpolation of $W^{1,p}$ and $W^{0,p}=L^p$. When $p=2$, $W_{0}^{\delta, 2} \cong L_{0}^{\delta, 2}=[W^{1,2},L^2]_{1-\delta}$, which can be defined as follows, 
	$$L^{\delta, 2}=\big\{f=c_{0}+\sum_{n=1}^{\infty} c_{n} \sqrt{2} \cos (n \pi x)|c_{n} \in \mathbf{C}, \sum_{n=0}^{\infty}\left(1+n^{2}\right)^{\delta} |c_{n}|^{2}<\infty\big\},$$
	then $W_{0}^{\delta, 2} $ and $ L_{0, \text { real }}^{\delta, 2}=[W^{1,2},L^2]_{1-\delta}$ are equivalent Hilbert spaces.

	\subsection{FBm and Bm}\label{sec-2-1}
	
	Consider the $\mathbb{R}^{d_1}$-valued fBm $(b^H_t)_{t\ge 0}$ with Hurst parameter $H\in (1/3,1/2)$,
	$$b^H_t=(b_t^{H,1},b_t^{H,2},\cdots,b_t^{H,d_1}),$$
	where  	 $(b_t^{H,i})_{t\ge 0}, i\in \{1,\cdots,d_1\}$ are independent one-dimensional fBms.
	The above $\mathbb{R}^{d_1}$-valued fBm $(b^H_t)_{t\ge 0}$ is a centred Gaussian process, satisfying that
	$$\mathbb{E}\big[b_{t}^{H} b_{s}^{H}\big]=\frac{1}{2}\left[t^{2 H}+s^{2 H}-|t-s|^{2 H}\right]\times I_{d_1}, \quad(s, t \geq 0),$$
	and 
	$$\mathbb{E}\big[\left(b_{t}^{H}-b_{s}^{H}\right)^{2}\big]=|t-s|^{2 H} \times I_{d_1},\quad (s, t \geq 0),$$
	where $I_{d_1}$ is the identity matrix in $\mathbb{R}^{d_1\times d_1}$. 
	When $H=1/2$, it is a standard Bm in $\mathbb{R}^{d_1}$. The reproducing kernel Hilbert space for the fBm $(b^H_t)_{t\ge 0}$ denoted by $\mathcal{H}^{H,d_1}$.

	 Then, 	consider the $\mathbb{R}^{d_2}$-valued  Bm $(w_t)_{t\ge 0}$,
	 	$$w_t=(w_t^1,w_t^2,\cdots,w_t^{d_2}),$$
	 where  	 $(w_t^i)_{t\ge 0}, i\in \{1,\cdots,d_2\}$ are independent one-dimensional  Bms. The reproducing kernel Hilbert space for  $(w_t)_{t\ge 0}$, denoted by $\mathcal{H}^{\frac{1}{2},d_2}$, which is defined as follows, 
	 $$\mathcal{H}^{\frac{1}{2},d_2}:=\big\{\hat k \in P(\mathbb{R}^{d_2}) \mid \hat k _{t}=\int_{0}^{t} \hat k _{s}^{\prime} d s \text { for all } t \text { with }\|\hat k \|_{\mathcal{H}^{\frac{1}{2},d_2}}^{2}:=\int_{0}^{1}|\hat k _{t}^{\prime}|_{\mathbb{R}^{d_2}}^{2} d t<\infty\big\},$$
	 where $P(\mathbb{R}^{d_2}):=\big\{\hat k \in C([0,1], \mathbb{R}^{d_2}) \mid \hat k _{0}=0\big\}$.
	 Due to [\cite{2013Inahama}, Proposition 3.4], it has ${\mathcal{H}^{H,d_1}} \hookrightarrow W_{0}^{\delta, 2} \cong L_{0, \text { real }}^{\delta, 2}$. 
	 
	 Let $\mathcal{H}:={{\mathcal{H}^{H,d_1}}\oplus{\mathcal{H}^{\frac{1}{2},d_2}}}$ be the Cameron-Martin subspace of the mixed fBm $(b_t^H, w_t)_{0\le t\le 1}$. Hence, $(k,\hat k)\in \mathcal{H}$ is of  finite $q$-variation with $(H+1/2)^{-1}<q<2$.

	\subsection{Rough path}\label{sec-2-2}

 Next, we introduce the geometric rough path. Set $2<p<3 $, and $\Delta=\{(s, t) \mid 0 \leq s \leq t \leq 1\}$, for the continuous map $A$ from $\Delta$ to the Banach space $B$, define that 
\begin{eqnarray*}
\|A\|_{p\mathrm{-var}}=\big(\sup _{\mathcal{P}} \sum_{i=1}^{n}\left|A_{t_{i-1}, t_{i}}\right|_{B}^{p}\big)^{1 / p},
\end{eqnarray*}
where $\mathcal{P}$ stands for the finite partition of $[s,t]$. 

Denote $B$ be the Banach space, a continuous map
\begin{eqnarray*}
X=\big(1, X^{1}, X^{2}\big): \Delta \rightarrow T^{2}(B)=\mathbb{R} \oplus B \oplus B^{\otimes 2} ,
\end{eqnarray*}
 is said to be a $B$-valued rough path of roughness $2$ if it satisfies the following  conditions,
 
 (Condition A): For any $s \leq u \leq t$, $X_{s, t}=X_{s, u} \otimes X_{u, t}$ where $\otimes$ stands for the tensor product.
 
 (Condition B) For all $1 \leq j \leq[p]$, $\left\|X^{j}\right\|_{p / j \mathrm{-var}}<\infty$.
 
 The $0$-th component $1$ is omitted. Therefore, we denote the rough path by $X=\left(X^{1},  X^{2}\right)$. The set of all the $B-$valued rough paths of roughness $2<p<3$ is denoted by $\Omega_{p}(B)$. With the distance $d_{p}(X, Y)=\sum_{i=1}^{\lfloor p\rfloor}\left\|X^{j}-Y^{j}\right\|_{p / j\mathrm{-var}}$, it is a complete space.
 
 For rough path $X$, denote that
 \begin{eqnarray*}
 	{\vertiii{X}_{p\mathrm{-var}}}
 	=\left\|X^{1}\right\|_{p-\mathrm{var}}+\left\|X^{2}\right\|_{p / 2-\mathrm{var}}^{1 / 2}.
 \end{eqnarray*}

 Next, for rough paths $X$, $Y$,  for $\kappa>p-1$,	define that
 \begin{eqnarray}\label{2-2}
 D_{j,p}\left( X,Y \right) =D_{j,p}\left( X^j,Y^j \right) =\left( \sum_{n=1}^{\infty}n^\kappa\sum_{l =1}^{2^n}{\left| X_{t_{l -1}^{n},t_{l}^{n}}^{j}-Y_{t_{l -1}^{n},t_{l}^{n}}^{j} \right|}^{p/j} \right) ^{j/p},
 \end{eqnarray}
 and
 \begin{eqnarray}\label{2-29}
 d_p\left( X,Y \right) \leq C\max \left( D_{\text{1,}p}\left( X,Y \right) ,D_{\text{1,}p}\left( X,Y \right) \left( D_{\text{1,}p}\left( X \right) +D_{\text{1,}p}\left( Y \right) \right) ,D_{\text{2,}p}\left( X,Y \right) \right). 
 \end{eqnarray}
 Then, from [\cite{2002Lyons}, Section 4.1], it has
 \begin{eqnarray*}
 \left\|X^{1}-Y^{1}\right\|_{p\mathrm{-var}}^{p} \leq c_{1} D_{1, p}(X, Y)^{p},
 \end{eqnarray*}
 and
 \begin{eqnarray*}
 &\left\|X^{2}-Y^{2}\right\|_{p / 2  {\mathrm{-var} }}^{p / 2} \cr
 &\quad \leq c_{1}\left[D_{2, p}(X, Y)^{p / 2}+D_{1, p}(X, Y)^{p / 2}\left(D_{1, p}(X)^{p}+D_{1, p}(Y)^{p}\right)^{1 / 2}\right],
 \end{eqnarray*}
 where $c_1$ is a constant.

 And $B$-valued continuous path $x$ with finite variation  can be lifted to the rough path  $X$, where the $j$-th level path is defined in the following way,
 	$$X_{s,t}^{j}=\int_{s\le t_1\le \cdots \le t_j\le t}{d}x_{t_1}\otimes dx_{t_2}\otimes \cdots \otimes dx_{t_j}.$$

The rough path obtained
 in above way is called smooth rough path. A
 rough path constructed as the $d_p$-limit of a sequence of smooth rough path is called a geometric rough path. The set of all the geometric rough paths is denoted by $G \Omega_{p}(B)$, and it is a complete separable metric space \cite{2002Lyons}.

For any $m \in \mathbb{N}$, consider the $m$-dyadic grid that  $t_{l}^{m} =\frac{l}{2^m}$ and set $\Delta_{l}^{m} (b^H,w)^\mathrm{T}=(b^H,w)^\mathrm{T}_{t_{l}^{m}}-(b^H,w)^\mathrm{T}_{t_{l-1}^{m}} $ for $0\le l\le  2^m$. Then denote by   $(b^H(m),w(m))_{t}$ the process obtained by linear interpolation of $(b^H,w)_t^\mathrm{T}$ on the $m$th dyadic grid. So $(b^H(m),w(m))_{0}=0$ and for $t \in[ t_{l-1}^{m}, t_{l}^{m}]$,
$$(b^H(m),w(m))^\mathrm{T}_{t}=(b^H(m),w(m))^\mathrm{T}_{t_{l-1}^{m}}+2^{m}\left(t-t_{l-1}^{m}\right) \Delta_{l}^{m} (b^H,w)^\mathrm{T}.$$ 
The corresponding smooth rough path $(B^H(m),W(m))_{s, t}=\left(1, (B^H(m),W(m))_{s, t}^{1}, (B^H(m),W(m))_{s, t}^{2}\right)$ is built by taking its iterated path integrals, that is 
 	$$(B^H(m),W(m))_{s,t}^{j}=\int_{s\le t_1\le \cdots \le t_j\le t}{d}(b^H(m),w(m))^\mathrm{T}_{t_1}\otimes d(b^H(m),w(m))^\mathrm{T}_{t_2}\otimes \cdots \otimes d(b^H(m),w(m))^\mathrm{T}_{t_j}.$$
For $n \le m$, $l=1, \ldots, 2^{n}$
$$(B^H(m),W(m))_{t_{l-1}^{n}, t_{l}^{n}}^{1}=\Delta_{l}^{n} (b^H,w)^\mathrm{T} \quad \text { for } n \le m,$$
and when $n \ge m$
$$(B^H(m),W(m))_{t_{l-1,}^{n}, t_{l}^{n}}^{1}=2^{m-n} \Delta_{\tilde l}^{m} (b^H,w)^\mathrm{T} \quad \text { for } n \ge m,$$
where $\tilde l$ is the unique integer $\tilde l$ among $1,\cdots,2^m$.

The second level path $(B^H(m),W(m))_{s, t}^{2}$ could be defined as follows,
when $n \le m$,
\begin{eqnarray*}
	(B^H(m),W(m))_{t_{l-1}^{n}, t_{l}^{n}}^{2}&=&\frac{1}{2} \Delta_{l}^{n} (b^H,w)^\mathrm{T} \otimes \Delta_{l}^{n} (b^H,w)^\mathrm{T}\cr
	&&+\frac{1}{2} \sum_{r, s=2^{m-n}(l-1)+1}^{2^{m-n} l}\left(\Delta_{r}^{m} (b^H,w)^\mathrm{T} \otimes \Delta_{s}^{m} (b^H,w)^\mathrm{T}-\Delta_{s}^{m} (b^H,w)^\mathrm{T} \otimes \Delta_{r}^{m} (b^H,w)^\mathrm{T}\right).
\end{eqnarray*}
While for the case of $n\ge m$,
$$(B^H(m),W(m))_{t_{l-1}^{n}, t_{l}^{n}}^{2}=2^{2(m-n)-1}\left(\Delta_{\tilde l}^{m} (b^H,w)^\mathrm{T} \otimes \Delta_{\tilde l}^{m} (b^H,w)^\mathrm{T}\right).$$
The  smooth rough path $(B^H(m),W(m))_{s, t}=\left(1, (B^H(m),W(m))_{s, t}^{1}, (B^H(m),W(m))_{s, t}^{2}\right)$ with order $2$ is  constructed.

In this paper, the mixed fBm $(b^H,w)^\mathrm{T}\in \mathbb{R}^{d_1+d_2}$ can be lifted to a geometric rough path $(B^H,W)\in G \Omega_{p}(\mathbb{R}^{d_1+d_2})$ with roughness $2<p<3$. 
\begin{prop}\label{prop2-1}
Let $(b_t^H,w_t)_{t\ge 0}^\mathrm{T}\in \mathbb{R}^{d_1+d_2}$ be the mixed fBm with Hurst parameter $1/3<H<1/2$. Then, for any $2<p<3$ such that $hp>1$, the sequence of its dyadic polygonal approximations $(B^H(m),W(m))=\left(1, (B^H(m),W(m))^{1}, (B^H(m),W(m))^{2}\right)$  converges to a unique geometric rough path $(B^H,W)=\left(1, (B^H,W)^{1}, (B^H,W)^{2}\right)$ almost surely according to the $p$-variation. 
\end{prop}	
\para{Proof}.  For $i=1,2$ and $0 \leq s \leq t \leq 1$, we proceed to prove that
\begin{eqnarray}\label{2-3}
\mathbb{E}\big[\left\|(B^H(m),W(m))^{j}-(B^H,W)^{j}\right\|^{p / j}_{p / j\mathrm{-var}}\big]\leq C\left(\frac{1}{2^{m}}\right)^{hp /2-1 / 2},
\end{eqnarray}
for some constant $C$. Since $p>2$, this implies that 
\begin{eqnarray}\label{2-4}
\sum_{m=1}^{\infty} \left\|(B^H(m),W(m))^{1}-(B^H,W)^{1}\right\|_{p\mathrm{-var}}<\infty, \quad a . s ..
\end{eqnarray}
In particular,  the sequence of its dyadic polygonal approximations $(B^H(m),W(m))$  converges to a unique geometric rough path $(B^H,W)$ almost surely according to the $p$-variation. 

The subsequent proof consists of two steps.	

\textbf{Step 1}.  Prove that (\ref{2-3}) and (\ref{2-4}) hold for $i=1$.
In particular, $(B^H(m),W(m))^1$ converges to  $(B^H,W)^1$  in $p$-variation almost surely.

Note that for  $n \le m$, $l=1, \ldots, 2^{n}$, it has 
$(B^H(m),W(m))_{t_{l-1}^{n}, t_{l}^{n}}^{1}-(B^H,W)_{t_{l-1}^{n}, t_{l}^{n}}^{1}=0$. On the other hand, for  $n > m$ and $\gamma>p / i-1$, we can conclude that
\begin{eqnarray}\label{2-5}
&&\mathbb{E}\big[\left\|(B^H(m),W(m))^{1}-(B^H,W)^{1}\right\|^{p }_{p \mathrm{-var}}\big]\cr&\leq&  C \sum_{n=m+1}^{\infty}n^\kappa\sum_{l =1}^{2^n}\mathbb{E}\big[\big| (B^H(m),W(m))_{t_{l-1}^{n},t_{l}^{n}}^{1}-(B^H,W)_{t_{l-1}^{n},t_{l}^{n}}^{1} \big|^{p}\big]\cr
&\leq&  C 2^{p-1}\sum_{n=m+1}^{\infty}n^\kappa\sum_{l =1}^{2^n}\mathbb{E}\big[\big| (B^H(m),W(m))_{t_{l-1}^{n},t_{l}^{n}}^{1}\big|^{p}\big]+\mathbb{E}\big[\big|(B^H,W)_{t_{l-1}^{n},t_{l}^{n}}^{1} \big|^{p}\big]   \cr
&\leq&C \sum_{n=m+1}^{\infty} n^{\kappa}\left(\frac{1}{2^{n}}\right)^{hp -1}\cr
&\leq&C\left(\frac{1}{2^{m}}\right)^{hp / 2-1 / 2} \sum_{n=m+1}^{\infty} n^{\kappa}\left(\frac{1}{2^{n}}\right)^{hp / 2-1 / 2},
\end{eqnarray}
the last series converges to a finite constant $C$. Moreover, we can assert that (\ref{2-3}) holds for $i=1$.

Then, with H\"older inequality and (\ref{2-5}), we have
\begin{eqnarray}\label{2-6}
&&\sum_{m=1}^{\infty}\mathbb{E}\big[\left\|(B^H(m),W(m))^{1}-(B^H,W)^{1}\right\|_{p\mathrm{-var}}\big] \cr&\le&\sum_{m=1}^{\infty}\mathbb{E}\bigg[ \sup _{\pi \in \Pi([0,1])}\big(\sum_{l}\big|(B^H(m),W(m))_{t_{L-1}, t_{L}}^{1}-(B^H,W)_{t_{L-1}, t_{L}}^{1}\big|^{p}\big)^{1 / p}\bigg]\cr
&\le&\sum_{m=1}^{\infty}\bigg[\mathbb{E}\big[ \sup _{\pi \in \Pi([0,1])}\sum_{l}\big|(B^H(m),W(m))_{t_{L-1}, t_{L}}^{1}-(B^H,W)_{t_{L-1}, t_{L}}^{1}\big|^{p}\big]\bigg]^{1 / p}\cr
&\leq&C\sum_{m=1}^{\infty}\left(\frac{1}{2^{m}}\right)^{h / 2-1 / 2p}.
\end{eqnarray}
Since $h / 2-1 / 2p>0$, the last series converges. Hence it deduces that (\ref{2-4}) holds for $i=1$. In particular, $(B^H(m),W(m))^1$ converges to  $(B^H,W)^1$  in $p$-variation almost surely.

\textbf{Step 2}. Prove that  (\ref{2-3}) and (\ref{2-4}) hold for $i=2$.
In particular, $(B^H(m),W(m))^2$ converges to  $(B^H,W)^2$  in $p$-variation almost surely.

For  $n\ge m$,
\begin{eqnarray*}\label{2-42}
	(B^H(m+1), W( m+1 )) _{t_{l-1}^{n},t_{l}^{n}}^{2}-(B^H(m), W( m )) _{t_{l-1}^{n},t_{l}^{n}}^{2}=2^{2\left( m+1-n \right) -1}\big( \Delta _{\tilde l}^{m+1}(b^H,w)^\mathrm{T} \big) ^{\otimes 2}-2^{2\left( m-n \right) -1}\left( \Delta _{\tilde l}^{m}(b^H,w)^\mathrm{T} \right) ^{\otimes 2}.
\end{eqnarray*}
Then by the triangle inequality, it has
\begin{eqnarray}\label{2-9}
&&\sum_{l=1}^{2^{n}} \mathbb{E}\big[\big|(B^H(m+1), W( m+1 )) _{t_{l-1}^{n},t_{l}^{n}}^{2}-(B^H(m), W( m )) _{t_{l-1}^{n},t_{l}^{n}}^{2}\big|^{p / 2}\big]\cr
&\le& {C} 2^{(2(m-n)+1) p / 2} \sum_{\tilde l=1}^{2^{m+1}} \mathbb{E}\big[\big|\big( \Delta _{\tilde l}^{m+1}(b^H,w)^\mathrm{T} \big) ^{\otimes 2}\big|^{p/2}\big]\cr
&&+{C} 2^{(2(m-n)-1) p / 2} \sum_{\tilde l=1}^{2^{m+1}}\frac{1}{2^{p}} \mathbb{E}\big[\big|\big( \Delta _{\tilde l}^{m}(b^H,w)^\mathrm{T} \big) ^{\otimes 2}\big|^{p/2}\big]\cr
&\leq& C 2^{(2(m-n)+1) p / 2} \sum_{\tilde l=1}^{2^{m+1}}\left(\frac{1}{2^{m+1}}\right)^{hp}+\frac{1}{2^{p}}\left(\frac{1}{2^{m}}\right)^{hp}\cr
&\leq&C\left(\frac{2^{m}}{2^{n}}\right)^{p-hp}\left(\frac{1}{2^{m}}\right)^{hp-1}.
\end{eqnarray}
 For $n\le m$,
 \begin{eqnarray*}\label{2-12}
 	&(B^H(m+1), W( m +1)) _{t_{l -1}^{n},t_{l}^{n}}^{2}-(B^H(m), W( m )) _{t_{l -1}^{n},t_{l}^{n}}^{2}\cr
 	\quad &=\frac{1}{2}\sum_{r=2^{m-n}\left( l -1 \right) +1}^{2^{m-n}l}{\left( \Delta _{2r-1}^{m+1}(b^H,w)^\mathrm{T}\otimes \Delta _{2r}^{m+1}(b^H,w)^\mathrm{T}-\Delta _{2r}^{m+1}(b^H,w)^\mathrm{T}\otimes \Delta _{2r-1}^{m+1}(b^H,w)^\mathrm{T} \right)}.
 \end{eqnarray*}
Denote $$
\Theta\left( n,m,l \right) =(B^H(m+1), W( m +1)) _{t_{l-1}^{n},t_{l}^{n}}^{2}-(B^H(m), W( m ))  _{t_{l-1}^{n},t_{l}^{n}}^{2},
$$
then, 
$$
\Theta\left( n,m,l \right) =\frac{1}{2}\sum_{r=2^{m-n}\left( l-1 \right) +1}^{2^{m-n}l}{\left( \Delta _{2r-1}^{m+1}(b^H,w)^\mathrm{T}\otimes \Delta _{2r}^{m+1}(b^H,w)^\mathrm{T}-\Delta _{2r}^{m+1}(b^H,w)^\mathrm{T}\otimes \Delta _{2r-1}^{m+1}(b^H,w)^\mathrm{T} \right)},$$
By the hypercontractivity inquality \cite{2007Friz}, it follows that
\begin{eqnarray}\label{2-13}
 \mathbb{E}\big[\big| \Theta\left( n,m,l \right) ^{i,j} \big|^{{p/2}}\big] ^{{2/p}}\leq C{\frac{p}{2}} \mathbb{E}\big[\big| \Theta\left( n,m,l \right) ^{i,j} \big|^2 \big]^{\frac{1}{2}}.
\end{eqnarray}
Moreover, 
\begin{eqnarray}\label{2-14}
\mathbb{E}\big[\big| \Theta\left( n,m,l \right) ^{i,j} \big|^2\big]\leq C\big( \Theta_{1}\left( n,m,l \right) ^{i,j}+\Theta_{2}\left( n,m,l \right) ^{i,j} \big), 
\end{eqnarray}
where
\begin{eqnarray*}\label{2-15}
	\Theta_{1}\left( n,m,l \right) ^{i,j}&=&\sum_{r=2^{m-n}\left( l-1 \right) +1}^{2^{m-n}l}\mathbb{E}\big[\left| \Delta _{2r-1}^{m+1}(b^H,w)^{\mathrm{T},i}\Delta _{2r}^{m+1}(b^H,w)^{\mathrm{T},j}-\Delta _{2r-1}^{m+1}(b^H,w)^{\mathrm{T},j}\Delta _{2r}^{m+1}(b^H,w)^{\mathrm{T},i} \right|^2\big],
\end{eqnarray*}
and
\begin{eqnarray}\label{2-16}
\Theta_{2}\left( n,m,l \right) ^{i,j}&=&\sum_{l\ne r}\mathbb{E}\big[ \left( \Delta _{2r-1}^{m+1}(b^H,w)^{\mathrm{T},i}\Delta _{2r}^{m+1}(b^H,w)^{\mathrm{T},j}-\Delta _{2r-1}^{m+1}(b^H,w)^{\mathrm{T},j}\Delta _{2r}^{m+1}(b^H,w)^{\mathrm{T},i} \right) \big.\cr
&&\quad\big. \left( \Delta _{2l-1}^{m+1}(b^H,w)^{\mathrm{T},i}\Delta _{2l}^{m+1}(b^H,w)^{\mathrm{T},j}-\Delta _{2l-1}^{m+1}(b^H,w)^{\mathrm{T},j}\Delta _{2l}^{m+1}(b^H,w)^{\mathrm{T},i} \right) \big].
\end{eqnarray}
On the one hand,   due to the  property of fBm and standard Bm, it follows that	
\begin{eqnarray}\label{2-51}
	\Theta_{1}\left( n,m,l \right) ^{i,j}
	&\leq& C\sum_{r=2^{m-n}\left( l-1 \right) +1}^{2^{m-n}l}\mathbb{E}\big[\left| \Delta _{2r-1}^{m+1}(b^H,w)^{\mathrm{T},i} \right|^2\big]\mathbb{E}\big[\left| \Delta _{2r}^{m+1}(b^H,w)^{\mathrm{T},j} \right|^2\big]\cr
	&\leq& C2^{-\frac{n}{2}}2^{-2m\left( H-\frac{1}{4} \right)},
\end{eqnarray}
On the other hand,  let
\begin{eqnarray}\label{2-17}
A&\equiv& \mathbb{E}\big[ \bigtriangleup _{2r-1}^{m+1}(b^H,w)^{\mathrm{T},i}\bigtriangleup _{2r}^{m+1}(b^H,w)^{\mathrm{T},j}\bigtriangleup _{2l-1}^{m+1}(b^H,w)^{\mathrm{T},i}\bigtriangleup _{2l}^{m+1}(b^H,w)^{\mathrm{T},j} \big]\cr
&&-\mathbb{E}\big[ \bigtriangleup _{2r}^{m+1}(b^H,w)^{\mathrm{T},i}\bigtriangleup _{2r-1}^{m+1}(b^H,w)^{\mathrm{T},j}\bigtriangleup _{2l-1}^{m+1}(b^H,w)^{\mathrm{T},i}\bigtriangleup _{2l}^{m+1}(b^H,w)^{\mathrm{T},j} \big]\cr
&&-\mathbb{E}\big[ \bigtriangleup _{2r-1}^{m+1}(b^H,w)^{\mathrm{T},i}\bigtriangleup _{2r}^{m+1}(b^H,w)^{\mathrm{T},j}\bigtriangleup _{2l}^{m+1}(b^H,w)^{\mathrm{T},i}\bigtriangleup _{2l-1}^{m+1}(b^H,w)^{\mathrm{T},j} \big]\cr
&&+\mathbb{E}\big[ \bigtriangleup _{2r}^{m+1}(b^H,w)^{\mathrm{T},i}\bigtriangleup _{2r-1}^{m+1}(b^H,w)^{\mathrm{T},j}\bigtriangleup _{2l}^{m+1}(b^H,w)^{\mathrm{T},i}\bigtriangleup _{2l-1}^{m+1}(b^H,w)^{\mathrm{T},j} \big]\cr
&=&2 \mathbb{E}\big[ \bigtriangleup _{2r-1}^{m+1}(b^H,w)^{\mathrm{T},i}\bigtriangleup _{2l-1}^{m+1}(b^H,w)^{\mathrm{T},i} \big]  \mathbb{E} \big[\bigtriangleup _{2r}^{m+1}(b^H,w)^{\mathrm{T},j}\bigtriangleup _{2l}^{m+1}(b^H,w)^{\mathrm{T},j} \big]\cr
&&-2\mathbb{E}\big[ \bigtriangleup _{2r}^{m+1}(b^H,w)^{\mathrm{T},i}\bigtriangleup _{2l-1}^{m+1}(b^H,w)^{\mathrm{T},i} \big] \mathbb{E}\big[\bigtriangleup _{2r-1}^{m+1}(b^H,w)^{\mathrm{T},j}\bigtriangleup _{2l}^{m+1}(b^H,w)^{\mathrm{T},j} \big].
\end{eqnarray}
When $i\in \{\text{1,}\cdots ,d_1\}$, $(b^H,w)^{\mathrm{T},i}$ are independent fBms for different $i$. 
With the result in [\cite{2002Coutin}, Section 3.3], it has
\begin{eqnarray}\label{2-18}
\mathbb{E}\big[\left( \bigtriangleup _{l}^{n}(b^H,w)^{\mathrm{T},i}\bigtriangleup _{{ l'}}^{n}(b^H,w)^{\mathrm{T},i} \right)\big] \le \left( 2h-1 \right) C\frac{|{ l'}-l|^{2h-2}}{\left( 2^n \right) ^{2h}},i\in \{\text{1,}\cdots ,d_1\}.
\end{eqnarray}
Besides, when $i\in \{d_1+\text{1,}\cdots ,d_1+d_2\}$, $(b^H,w)^{\mathrm{T},i}$ are independent  Bms for different $i$, we have
\begin{eqnarray}\label{2-19}
\mathbb{E}\big[ \bigtriangleup _{l}^{n}(b^H,w)^{\mathrm{T},i}\bigtriangleup _{{ l'}}^{n}(b^H,w)^{\mathrm{T},i} \big] =0,
\end{eqnarray}
for ${ l'}\ne l$. Meanwhile, it has
\begin{eqnarray}\label{2-20}
\mathbb{E}\big[\bigtriangleup _{l}^{n}(b^H,w)^{\mathrm{T},i}\bigtriangleup _{{ l'}}^{n}(b^H,w)^{\mathrm{T},j}\big]   =\text{0, }i\ne j.
\end{eqnarray}
From the above it follows that,  
\begin{eqnarray}\label{2-52}
	\Theta_{2}\left( n,m,l \right) ^{i,j}
	&\leq& C2^{-\frac{n}{2}}2^{-2m\left( H-\frac{1}{4} \right)},
\end{eqnarray}
for 	 $i,j\in \{1,\cdots ,d_1+d_2\}$.  

Therefore,  by the  (\ref{2-51}) and (\ref{2-52}), it has
\begin{eqnarray}\label{2-21}
 \mathbb{E}\big[\big| \Theta\left( n,m,l \right) ^{i,j} \big|^{2} \big] ^{{1/2}}\leq C2^{-\frac{n}{2}}2^{-2m\left( H-\frac{1}{4} \right)},
\end{eqnarray}
for 	 $i,j\in \{1,\cdots ,d_1+d_2\}$.  

When $n=m$, we get
\begin{eqnarray}\label{2-22}
 \mathbb{E}\big[\big| \Theta\left( n,n,l \right) ^{i,j} \big|^{2} \big] ^{{1/2}}\leq C2^{-2nH}.
\end{eqnarray}
Combined with (\ref{2-29}), (\ref{2-9}), and (\ref{2-21}), it deduces that
\begin{eqnarray}\label{2-23}
&&\mathbb{E}\big[\left\|(B^H(m),W(m))^{2}-(B^H,W)^{2}\right\|^{p / 2}_{p / 2\mathrm{-var}}\big]\cr
&\leq&C \sum_{n=1}^{m} n^{\kappa}\left(\frac{1}{2^{n+m}}\right)^{hp / 2-1 / 2}+C\sum_{n=m+1}^{\infty} n^{\kappa}\left(\frac{2^{m}}{2^{n}}\right)^{p-hp}\left(\frac{1}{2^{m}}\right)^{hp-1}\cr
&&+C\bigg[\sum_{n=m+1}^{\infty} n^{\kappa}\left(\frac{1}{2^{n}}\right)^{hp-1}\bigg]^{1/2}\cr
&\leq&C\left(\frac{1}{2^{m}}\right)^{(hp-1)/2 },
\end{eqnarray}
hence, (\ref{2-3}) holds for $i=2$. By the  straightforward computation, it follows that (\ref{2-4}) holds for $i=2$.
Therefore, $(B^H(m),W(m))^{2}$ is a Cauchy sequence in $p$-variation, and the conclusion follows.

In above, 
 the sequence of its dyadic polygonal approximations $(B^H(m),W(m))$  converges to a unique function  $(B^H,W)$ almost surely in $p$-variation distance. According to the definition, $(B^H,W)$ is a geometric rough path.
The proof is completed. \qed	

\begin{prop}\label{prop2-2}
	Let $(b_t^H,w_t)_{t\ge 0}^\mathrm{T}\in \mathbb{R}^{d_1+d_2}$ be the mixed fBm with $\mathbb{R}^{d_1}$-valued fBm $(b_t^H)_{t\ge 0}$ ($1/3<H<1/2$) and $\mathbb{R}^{d_2}$-valued standard Bm $(w_t)_{t\ge 0}$.  Then,
	\begin{eqnarray}\label{2-28}
	(B^H,W)^1_{s t}=\left(b^H_{s t}, w_{s t}\right)^{\mathrm{T}}, \quad (B^H,W)^2_{s t}=\left(\begin{array}{ll}
	B_{s t}^{H,2} & I[b_H, w]_{s t} \\
	I[w, b_H]_{s t} & W_{s t}^{2}
	\end{array}\right),
	\end{eqnarray}
	 where  the $(B^{H,1}, B^{H,2})$ is a canonical geometric rough path, and the $(W^{1}, W^{2})$ is a geometric rough path in Stratonovich sense. Then
	 \begin{eqnarray}\label{2-30}
     I[b^H, w]_{s t} \triangleq \int_{s}^{t} b^H_{s u} \otimes \mathrm{d}^{\mathrm{I}} w_{u},
	 \end{eqnarray}
	 \begin{eqnarray}\label{2-31}
     I[w,b^H]_{s t} \triangleq w_{s t} \otimes b^H_{s t}-\int_{s}^{t} \mathrm{~d}^{\mathrm{I}} w_{u} \otimes b^H_{s u},
	 \end{eqnarray}
	 	 where $\int \cdots \mathrm{d}^{\mathrm{I}} w$ stands for the It\^o intrgral.
\end{prop}	
\para{Proof}.  $(B^{H,1}, B^{H,2})$ is a canonical geometric rough path (see [\cite{2020Friz},Section 10.3]). And $(W^{1}, W^{2})$ is a geometric rough path in Stratonovich sense (see [\cite{2020Friz},Section 3]). It remains to show (\ref{2-30}) and (\ref{2-31}) hold. Take $s=0,t=1$ for simplicity. Then denote by   $(b_t^H(m))_{0\le t\le1}$ the process obtained by linear interpolation of $(b_t^H)_{0\le t\le1}$ on the $m$th dyadic grid. Similar for $(w_t(m))_{0\le t\le1}$. To show (\ref{2-30}) hold, it turns to prove that 
\begin{eqnarray}\label{2-24}
\underset{m\rightarrow \infty}{\lim}\int_0^1{b_{t}^{H,i}\left( m \right) dw_t^j\left( m \right)}=\int_0^1{b_{t}^{H,i}dw^j_t},
\end{eqnarray}
for $i\in\{1,2,\cdots,d_1\}$ and $j\in\{1,2,\cdots,d_2\}$. 
The right hand side of  (\ref{2-24}) is in the It\^o sense.

To this end, define a step function $\hat b^{H,i}_t=b^{H,i}_{t^m_{k-1}}$ for $t\in [t^m_{k-1},t^m_{k}]$ and $1\le k\le 2^m$. Then it decuces that 

\begin{eqnarray}\label{2-25}
\mathbb{E} \big[|\int_0^1{(b_{t}^{H,i}-\hat b_{t}^{H,i}) dw^j_t}|^2\big]\to 0,
\end{eqnarray}
as $m\to \infty$. Hence, the right hand side of  (\ref{2-24}) is in the It\^o sense.
According to the definition of dyadic approximation, the left hand side of  (\ref{2-24}) can be rewritten as follows,
\begin{eqnarray}\label{2-26}
\int_0^1{b_{t}^{H,i}\left( m \right) dw^j_t\left( m \right)}
&=& \sum_{k=1}^{2^m} \int_{t^m_{k-1}}^{t^m_k} b_u^{H,i}(m) \big(\frac{w^j_{t^m_k}-w^j_{t^m_{k-1}}}{1/2^m}\big)du  .
\end{eqnarray}
Next, we get that 
\begin{eqnarray}\label{2-27}
&&\mathbb{E}\left[\bigg|\sum_{k=1}^{2^m} \bigg(2^m \int_{t^m_{k-1}}^{t^m_k} b_u^{H,i}(m) du-b_{t^m_{k-1}}^{H,i}(m)\bigg)\big(w^j_{t^m_k}-w^j_{t^m_{k-1}}\big)\bigg|^2\right] \cr
&=& \sum_{k=1}^{2^m} \mathbb{E}\bigg[\bigg|2^m \int_{t^m_{k-1}}^{t^m_k} b_u^{H,i}(m) du-b_{t^m_{k-1}}^{H,i}(m)\bigg|^2\bigg]\mathbb{E}\big[\big(w^j_{t^m_k}-w^j_{t^m_{k-1}}\big)^2\big] \cr
&=& \sum_{k=1}^{2^m}2^m \mathbb{E}\bigg[\bigg| \int_{t^m_{k-1}}^{t^m_k} \big(b_u^{H,i}(m) -b_{t^m_{k-1}}^{H,i}(m)\big)du\bigg|^2\bigg]\cr
&\le& \sum_{k=1}^{2^m} \int_{t^m_{k-1}}^{t^m_k}\mathbb{E}\bigg[\bigg|  b_u^{H,i}(m) -b_{t^m_{k-1}}^{H,i}(m)du\bigg|^2\bigg]\cr
&\le& C \sum_{k=1}^{2^m} \int_{t^m_{k-1}}^{t^m_k} (u-t^m_{k-1})^{2H}du\cr
&\le& C 2^{-2mH},
\end{eqnarray}
then, it follows that (\ref{2-24}) holds. Likewise, it asserts that (\ref{2-31}) holds. 	 The above definition  coincides with  [\cite{2021Pei}, Proposition 1.1].

The proof is completed. \qed

	\subsection{Assumptions and main results}\label{sec-2-3}

	Denote $C_0^{p\mathrm{-var}}(\mathbb{R}^n)$ the space of continuous paths in $\mathbb{R}^n$ with $p$-variation ($2<p<3$) starting at 0. Let $Y^{\varepsilon}=\hat{\Phi}_{\varepsilon}(\varepsilon (B^H, W), \lambda): G \Omega_{p}\left(\mathbb{R}^{d_1+d_2+1}\right) \mapsto G \Omega_{p}\left(\mathbb{R}^{n}\right)$ denote the It\^o map corresponding to (\ref{1-1}) with $\lambda_{t}=t$. Denote the $Y^{\varepsilon,1}$ is the first level path of the solution map.
	 Consider the solution map $\hat{\Phi}_{\varepsilon}(\varepsilon (B^H,W)+(\gamma, \eta), \lambda)$,
	\begin{eqnarray}\label{3-2-2}
	d\tilde Y_{t}^{\varepsilon}=\big[\sigma \big( \tilde Y_{t}^{\varepsilon} \big)|\hat\sigma \big( \tilde Y_{t}^{\varepsilon} \big)\big] (\varepsilon d(B^H,W)_t+(\gamma_t, \eta_t)^\mathrm{T})+\beta \big(\varepsilon,\tilde Y_{t}^{\varepsilon} \big) dt, \quad \tilde Y_{0}^{\varepsilon}=\text{0}.
	\end{eqnarray}
	Denote $\phi^{(\varepsilon)}=\tilde{Y}^{\varepsilon,1}$ the first level path of the solution map.
	Now $\hat{\Phi}_{0}((\gamma, \eta), \lambda)$ is the solution map lying above $\phi^{0}=\Psi(\gamma, \eta) \in C_{0}^{q-\operatorname{var}}\left(\mathbb{R}^{n}\right)$ with $1< q<2$, satisfying that
	\begin{eqnarray}\label{3-2-3}
	d \phi_{t}^{0}=\big[\sigma \left(\phi_{t}^{0}\right)|\hat\sigma \left(\phi_{t}^{0}\right)\big] d (\gamma_t, \eta_t)^\mathrm{T}+\beta\left(0,\phi_{t}^{0}\right) d t, \quad \phi_{0}^{0}=0.
	\end{eqnarray}
	
	We assume
	\begin{itemize}
		\item[\textbf{A1}.]  The function $F$ is real-valued bounded continuous  on $C_0^{p'\mathrm{-var}}(\mathbb{R}^n)$ with $p'>1/H$.
		\item[\textbf{A2}.] The function $F_{\Lambda}:=F\circ\Phi+||(\cdot,\cdot)||^2_{\mathcal{H}} /2$ attains its minimum at a unique point $(\gamma,\eta)^\mathrm{T}
		\in{\mathcal{H}}$, with $\Phi(\gamma,\eta)=\phi^0$.
		\item[\textbf{A3}.]The function  $F$ is $m+3$ times Fr\'echet differentiable on a neighborhood $U(\phi^0)$ with $\phi^0\in C_0^{p'\mathrm{-var}}(\mathbb{R}^n)$. Furthermore, there exist positive constants $M_1,M_2,\cdots$ such that
		$$\left|\nabla^{j} F(\eta)\langle z, \ldots, z\rangle\right| \leq M_{j}\|z\|_{p'\mathrm{-var}}^{j} \quad(j=1, \ldots, m+3),$$
		hold for any $\eta\in U(\phi^0)$ and $z\in C_0^{p'\mathrm{-var}}(\mathbb{R}^n)$.
		\item[\textbf{A4}.]  The bounded self-adjoint operator $A$ on ${\mathcal{H}}$, which is related to the Hessian matrix  $\left.\nabla^{2}(F \circ \Psi)(\gamma,\eta)\right|_{{\mathcal{H}} \times {\mathcal{H}}}$, is strictly larger than $-Id_{\mathcal{H}}$.
	\end{itemize}

Then, here follows our main result.
	
	\begin{thm}\label{thm2-1}
		Under  Assumptions {\rm (A1)}-{\rm (A4)}, we have the following asymptotic expansion with $\varepsilon \rightarrow 0$. There exist constants $c, \alpha _0, \alpha_1,\cdots$ s.t.
		\begin{eqnarray}\label{2-1}
		\mathbb{E}\left[ \exp \left( -F\left( Y^{\varepsilon ,1} \right) /\varepsilon ^2 \right) \right] 
		=\exp \left( -F_{\Lambda}(\gamma,\eta) /\varepsilon ^2 \right)\left( -c/\varepsilon  \right) \cdot \left( \alpha _0+\alpha _1\varepsilon +\cdots +\alpha _m\varepsilon ^m+O\left( \varepsilon ^{m+1} \right) \right),
		\end{eqnarray}
		for any $m\ge0$.
	\end{thm}	
   \begin{rem}\label{rem2-2}
    		Denote  $\tau_{p,p'}$ the injection from $C^{p \mathrm{-} var}(\mathbb{R}^n)$ to $C^{p' \mathrm{-} var}(\mathbb{R}^n)$. Then, we focus on the Laplace approximation for the first level path to the solution map, that is $F\big(\tau_{p,p'}(Y_{t}^{\varepsilon,1})\big)$. But for simplicity, we write in the sense that $F\big(Y_{t}^{\varepsilon,1}\big)$. The another way of saying that one views $F$ as a function on $C^{p \mathrm{-} var}(\mathbb{R}^n)$, even if it can be done and assumptions {\rm (A1)}-{\rm (A4)} remains equivalent.
    \end{rem}

	\section{Schilder-type large deviation principle}\label{sec-3}

	In this section, we prove the Schilder-type LDP for the law of the first leve path of the solution map to the RDE (\ref{1-1}) with $\varepsilon \to 0$.

	\begin{prop}\label{prop3-2-1}
		With $\varepsilon \to 0$, the law of $\varepsilon(B^H, W)$, $\left\{\mathbb{P}_{\varepsilon}^{H}\right\}_{\varepsilon>0}$ satisfies the LDP with the good rate function $I$, which is defined as follows,
		\begin{eqnarray}\label{eq-3-2}
		I(B^H, W)= \begin{cases}\frac{1}{2}\|(k,\hat k)\|_{\mathcal{H}}^{2} & \big(\text { if } (B^H, W) \text { is lying above } (k,\hat k)^\mathrm{T} \in {\mathcal{H}}\big), \\ \infty & (\text { otherwise }).\end{cases}
		\end{eqnarray}
	\end{prop}	
		
	\para{Proof}. 
	The subsequent proof consists of several steps.	
	
	\textbf{Step 1}.  Prove that the smooth rough paths $(B^H(m),W(m))$ are exponentially good approximations of $(B^H,W)$  in the sense that, for every $\delta>0$,
	\begin{eqnarray}\label{eq-3-3}
	\lim _{m \rightarrow \infty} \lim \sup _{\varepsilon \rightarrow 0} \varepsilon^{2} \log \mathbb{P}\left(d_{p}(F(\varepsilon (B^H(m), W( m ))), F(\varepsilon (B^H, W )))>\delta\right)=-\infty.
	\end{eqnarray}
	Let $pH>1$, it will be enough to show that there is a sequence $c(m)\to 0$ such that for ${\varrho}> 1$, 
		\begin{eqnarray}\label{eq-3-4}
		 \mathbb{E}\big[D_{\text{1,}p}\left( (B^H(m), W( m )) ,(B^H, W) \right) ^{\varrho} \big]   ^{\text{1/}{\varrho}}\leq c\left( m \right) \sqrt{{\varrho}},	
		\end{eqnarray}
		and
		\begin{eqnarray}\label{3-4}
\mathbb{E}\big[D_{1, p}(B^H(m),W(m))^{{\varrho}}\big]^{1 / {\varrho}} \leqslant C \sqrt{{\varrho}},
		\end{eqnarray}
and
		\begin{eqnarray}\label{eq-3-5}
		 \mathbb{E}\big[ D_{\text{2,}p}\left( (B^H(m), W( m )) ,(B^H, W) \right) ^{\varrho} \big]  ^{\text{1/}{\varrho}}\leq c\left( m \right) {\varrho},	
		\end{eqnarray}
		where $c(m)$ tends to zero as $m \to \infty$.
		
	 Accordingly, by the Chebyshev inequality, for ${\varrho}={\varrho}(\varepsilon)=\varepsilon^{-2}$, and $j=1,2$, we can conclude  
	 $$\mathbb{P}\left( D_{j,p}\left( (B^H(m), W( m )) ,(B^H, W) \right) >\delta \varepsilon ^{-j} \right) \leq \left( \delta ^{-1}\varepsilon ^j \right) ^{\varrho}c\left( m \right) ^{\varrho}{\varrho}^{j{\varrho}/2}\leq \left( \delta ^{-1}c\left( m \right) \right) ^{\varrho}.	$$
	 It deduces that
	 \begin{eqnarray}\label{3-2}
	 \lim_{\varepsilon \rightarrow 0} \text{sup} \varepsilon ^2\log \mathbb{P}\left( D_{j,p}\left((B^H(m), W( m )) ,(B^H, W) \right) >\delta \varepsilon ^{-j} \right) \leq \log \left( \delta ^{-1}c\left( m \right) \right) .	
	 \end{eqnarray}
	 Besides, since that for some constant $C$, we have $\mathbb{E}\left(D_{1, p}(B^H,W)^{{\varrho}}\right)^{1 / {\varrho}} \leqslant C \sqrt{{\varrho}}$. Then by the Cauchy-Schwarz inequality and  triangle inequalities, for every ${\varrho}$, it has
	 \begin{eqnarray*}
     &&\mathbb{E}\big[\left[D_{1, p}((B^H(m), W( m )), (B^H, W) )\left(D_{1, p}(B^H(m), W( m ))+D_{1, p}(B^H, W) \right)\right]^{{\varrho}}\big]^{1 / {\varrho}} \cr
     &&\quad\leqslant 2 c(m)(2 C+c(m)) {\varrho},
	 \end{eqnarray*}
  in consequence,
    \begin{eqnarray}\label{3-3}
	&&\limsup _{\varepsilon \rightarrow 0} \varepsilon^{2} \log \mathbb{P}\left(D_{1, p}((B^H(m), W( m )), (B^H, W))\left(D_{1, p}(B^H(m), W( m ))+D_{1, p}(B^H, W)\right)>\delta \varepsilon^{-2}\right) \cr
	&&\quad \leqslant \log \left(\delta^{-1} 2 c(m)(2 C+c(m))\right),
    \end{eqnarray}
    combined with (\ref{3-2}), we can assert that (\ref{eq-3-3}) holds.
    
	Therefore, we proceed to prove (\ref{eq-3-4}) and (\ref{eq-3-5}). First, for any ${\varrho}>p$, it has
   \begin{eqnarray}\label{eq-3-6}
   \mathbb{E}\big[ D_{\text{1,}p}\left( (B^H(m), W( m )) ,(B^H, W)\right) ^{\varrho} \big] \leq A\left( m,{\varrho} \right) \sum_{n=m+1}^{\infty}{a}_{n}^{{\varrho}/p}\sum_{l =1}^{2^n}{\mathbb{E}}\left[ \big| 2^{m-n}\Delta _{{\tilde l}}^{m}(b^H,w)^\mathrm{T}-\Delta _{l}^{n}(b^H,w)^\mathrm{T} \big|^{\varrho} \right], 
   \end{eqnarray}
	where $A\left( m,{\varrho} \right) =\big( \sum_{n=m+1}^{\infty}{2}^n\big( \frac{n^\kappa}{a_n} \big) ^{{\varrho}/\left( {\varrho}-p \right)} \big) ^{\left( {\varrho}-p \right) /p}.$
	
	Then due to the fact that 	when $i\in \{\text{1,}\cdots ,d_1\}$, $(b^H,w)^{\mathrm{T},i}$ are independent fBms, when $i\in \{d_1+\text{1,}\cdots ,d_1+d_2\}$, $(b^H,w)^{\mathrm{T},i}$ are independent  Bms. Choose $
	A\left( m,{\varrho} \right) =\left( \sum_{n=m+1}^{\infty}{2}^n\left( n^\kappa/a_n \right) ^{\frac{{\varrho}}{{\varrho}-p}} \right) ^{\frac{{\varrho}-p}{p}},a_n=2^{np\left( H-\frac{1}{{\varrho}}-\beta_1 \right)}
	$ and $n^\kappa \le C2^{np\beta_2}$ with $\beta_2>0$  and $\beta_1+\beta_2\in(0,(H-1/p)/2)$, we get
	\begin{eqnarray}\label{eq-3-7}
	\mathbb{E}\big[ D_{\text{1,}p}\left( (B^H(m), W( m )) ,(B^H, W) \right) ^{\varrho} \big] &\leq& A\left( m,{\varrho} \right) \sum_{n=m+1}^{\infty}{a_{n}^{p/{\varrho}}}\sum_{l=1}^{2^n}\mathbb{E}\big[ \left| 2^{m-n}\Delta _{k}^{m}(b^H,w)^\mathrm{T}-\Delta _{l}^{n}(b^H,w)^\mathrm{T} \right|^{\varrho} \big]\cr
	&\leq& A\left( m,{\varrho} \right) \left( 2d_1 \right) ^{\varrho}{\varrho}^{\frac{{\varrho}}{2}}\sum_{n=m+1}^{\infty}{a_{n}^{p/{\varrho}}}2^n\left( 2^{\left( m-n \right) {\varrho}}2^{-m{\varrho}H}+2^{-n{\varrho}H} \right)\cr
	&&+A\left( m,{\varrho} \right) \left( 2d_2 \right) ^{\varrho}{\varrho}^{\frac{{\varrho}}{2}}\sum_{n=m+1}^{\infty}{a_{n}^{p/{\varrho}}}2^n\left( 2^{\left( m-n \right) {\varrho}}2^{-m{\varrho}/2}+2^{-n{\varrho}/2} \right)\cr
	&\leq& A\left( m,{\varrho} \right) \left( 2d_1+2d_2 \right) ^{\varrho}{\varrho}^{\frac{{\varrho}}{2}}\sum_{n=m+1}^{\infty}{a_{n}^{p/{\varrho}}}2^n\left( 2^{\left( m-n \right) {\varrho}}2^{-m{\varrho}H}+2^{-n{\varrho}H} \right)\cr
	&\leq& C\left( 2d_1+2d_2 \right)^{\varrho}{\varrho}^{\frac{{\varrho}}{2}}2^{-m{\varrho}\beta_1},
    \end{eqnarray}
    where the second inequality is due to Lemma 4 in \cite{2002Ledoux} and $A\left( m,{\varrho} \right)$ is a bounded double-sequence in $m$ and ${\varrho}$.
    
     Similarly, it has
     \begin{eqnarray}\label{3-7}
     \mathbb{E}\big[D_{1, p}(B^H(m), W( m )) ^{{\varrho}}\big] \leqslant A(0, {\varrho})\left( 2d_1+2d_2 \right)^{\varrho} {\varrho}^{\frac{{\varrho}}{2}} \sum_{n=1}^{\infty} a_{n}^{\frac{{\varrho}}{p}} 2^{n(1-{\varrho} H)}<\infty.
     \end{eqnarray}
    Combined with (\ref{eq-3-7}) and the triangle inequality, we can prove that
    \begin{eqnarray}\label{3-8}
         \mathbb{E}\big[D_{1, p}(B^H,W)^{{\varrho}}\big]  <\infty.
    \end{eqnarray}
    
	When $j=2$, we proceed to assert that (\ref{eq-3-5}) holds. Consider the case of  $n\ge m$,
	\begin{eqnarray*}\label{eq-3-8}
	(B^H(m), W( m )) _{t_{l-1}^{n},t_{l}^{n}}^{2}=2^{2\left( m-n \right) -1}\left( \Delta _{\tilde l}^{m}(b^H,w)^\mathrm{T} \right) ^{\otimes 2}.
	\end{eqnarray*}
	then due to the  property of fBm and standard Bm, it follows that	
	\begin{eqnarray}\label{eq-3-9}
 \mathbb{E}\big[\big|(B^H(m), W( m )) _{t_{l-1}^{n},t_{l}^{n}}^{2} \big|^{\varrho} \big] ^{1/\varrho}\leq C{\varrho}2^{-2n}2^{-2m\left( H-1 \right)},
	\end{eqnarray}
	and 
	\begin{eqnarray}\label{eq-3-10}
	\mathbb{E}\big[\big| (B^H, W)_{t_{l-1}^{n},t_{l}^{n}}^{2} \big|^{\varrho} \big] ^{1/\varrho}\leq C{\varrho}2^{-2nH}.
    \end{eqnarray}
	Likewise, we can see
	\begin{eqnarray}\label{eq-3-11}
	 \mathbb{E}\big[| (B^H, W)_{t_{l-1}^{n},t_{l}^{n}}^{2}-(B^H(m), W( m ))  _{t_{l-1}^{n},t_{l}^{n}}^{2}|^{\varrho} \big] &\leq& 2^{\left( {\varrho}-1 \right)}\left( \mathbb{E}| (B^H, W)_{t_{l-1}^{n},t_{l}^{n}}^{2}| ^{\varrho}+| (B^H(m), W( m )) _{t_{l-1}^{n},t_{l}^{n}}^{2}| ^{\varrho} \right) 
	\cr
	&\leq& C2^{\left( {\varrho}-1 \right)}{\varrho}2^{-2n{\varrho}H}+C2^{\left( {\varrho}-1 \right)}{\varrho}2^{-2\left( n-m \right) {\varrho}}2^{-2m{\varrho}H}
	\cr
	&\leq& C{\varrho}2^{-2n{\varrho}H}.
	\end{eqnarray}
	When $n\le m$,
	\begin{eqnarray*}\label{eq-3-12}
	&(B^H(m+1), W( m +1)) _{t_{l -1}^{n},t_{l}^{n}}^{2}-(B^H(m), W( m )) _{t_{l -1}^{n},t_{l}^{n}}^{2}\cr
	\quad &=\frac{1}{2}\sum_{r=2^{m-n}\left( l -1 \right) +1}^{2^{m-n}l}{\left( \Delta _{2r-1}^{m+1}(b^H,w)^\mathrm{T}\otimes \Delta _{2r}^{m+1}(b^H,w)^\mathrm{T}-\Delta _{2r}^{m+1}(b^H,w)^\mathrm{T}\otimes \Delta _{2r-1}^{m+1}(b^H,w)^\mathrm{T} \right)}.
	\end{eqnarray*}
   Take same manner from (\ref{2-13}) to (\ref{2-22}), it follows that 
	\begin{eqnarray}\label{eq-3-13}
   \mathbb{E}\big[\big|(B^H(m+1), W( m +1)) _{t_{l -1}^{n},t_{l}^{n}}^{2}-(B^H(m), W( m )) _{t_{l -1}^{n},t_{l}^{n}}^{2}\big|^{\varrho}\big]^{\frac{1}{{\varrho}}}\leq C2^{-\frac{n}{2}}2^{-2m\left( H-\frac{1}{4} \right)}.
	\end{eqnarray}
	Fix $M$ large enough, we can obtain that
	\begin{eqnarray*}\label{eq-3-23}
	\mathbb{E}\big[\big| (B^H(M), W( M )) _{t_{l-1}^{n},t_{l}^{n}}^{2}-(B^H(m), W( m )) _{t_{l-1}^{n},t_{l}^{n}}^{2} \big|^{\varrho} \big] ^{\frac{1}{{\varrho}}}&\leq& C{\varrho}2^{-\frac{n}{2}}\sum_{N=m}^{M-1}{2}^{-2N\left( H-\frac{1}{4} \right)}\cr
	&\leq& C{\varrho}2^{-\frac{n}{2}}2^{-2m\left( H-\frac{1}{4} \right)}.
	\end{eqnarray*}
	By the construction of the rough path lying above,
	\begin{eqnarray*}\label{eq-3-24}
	\lim_{M\rightarrow \infty} (B^H(M), W( M )) _{t_{l-1}^{n},t_{l}^{n}}^{2}=(B^H, W )_{t_{l-1}^{n},t_{l}^{n}}^{2}.
	\end{eqnarray*}
	Hence, we get 
	\begin{eqnarray}\label{eq-3-25}
 \mathbb{E}\big[\big| (B^H, W)_{t_{l-1}^{n},t_{l}^{n}}^{2}-(B^H(m), W( m ))_{t_{l-1}^{n},t_{l}^{n}}^{2} \big|^{\varrho} \big] ^{1/\varrho}\leq C{\varrho}2^{-\frac{n}{2}}2^{-2m\left( H-\frac{1}{4} \right)}.
	\end{eqnarray}
	It means that
	\begin{eqnarray*}\label{eq-3-26}
	\mathbb{E}\big[ D_{\text{2,}p}\left( (B^H(m), W( m )) ,(B^H, W) \right) ^{\varrho} \big] &\leq& A\left( {\varrho} \right) \sum_{n=1}^{\infty}{a_{n}^{2{\varrho}/p}}\sum_{l=1}^{2^n}\mathbb{E}\big[\big| (B^H(m), W( m )) _{t_{l-1}^{n},t_{l}^{n}}^{2}-(B^H, W)_{t_{l-1}^{n},t_{l}^{n}}^{2} \big|^{\varrho}\big]
	\cr
	&\leq& CA\left( {\varrho} \right) {\varrho}^{\varrho}\bigg[ \sum_{n=1}^m{a_{n}^{2{\varrho}/p}}2^{-n\left( \frac{{\varrho}}{2}-1 \right)}2^{-2m{\varrho}\left( H-\frac{1}{4} \right)}+\sum_{n=m+1}^{\infty}{a_{n}^{2{\varrho}/p}}2^{-n\left( 2{\varrho}H-1 \right)} \bigg] ,
	\end{eqnarray*}
	where $A\left( {\varrho} \right) =\big( \sum_{n=1}^{\infty}{2^n}\left( n^\kappa/a_n \right) ^{\frac{2{\varrho}}{2{\varrho}-p}} \big) ^{\frac{2{\varrho}-p}{p}}.$
	
	Then choose approriate $a_{n}=2^{-n p\left(\beta_3-H+\frac{\beta_4}{2}+\frac{1}{2 {\varrho}}\right)}$ and $n^\kappa \le C2^{np\beta_5}$ with $\beta_3>0$, $\beta_4\in(0,2H-1/2)$, and $\beta_3+\frac{\beta_4}{2}+\beta_5<H-\frac{1}{p}$, the series $\sum_{n} a_{n}^{\frac{2 {\varrho}}{p}} 2^{-n[{\varrho}(2 H-\beta_4)-1]}$ converges. $A\left( {\varrho} \right)$ is a bounded double-sequence in  ${\varrho}$. Moreover,
	\begin{eqnarray}\label{eq-3-27}
	\mathbb{E}\big[D_{2, p}((B^H(m), W( m )), (B^H, W ))^{{\varrho}}\big]^{1/\varrho} \leq C \varrho 2^{-m \rho},
	\end{eqnarray}
	where $\rho>0$.
	
	In above, (\ref{eq-3-4}), (\ref{3-4}) and (\ref{eq-3-5}) are proved. Furthermore,  the smooth rough paths $(B^H(m),W(m))$ are exponentially good approximations of the geometric rough path $(B^H,W)$.

	\textbf{Step 2}. Prove that the rough path $F(K,\hat K)$ above any element $(k,\hat k)$ in the
	Cameron-Martin space ${\mathcal{H}}$ is defined as the limit of $F(K(m),\hat K( m ))$,
	\begin{eqnarray}\label{eq-3-29}
	\lim _{m \rightarrow \infty} \sup _{|(k,\hat k)|_{\mathcal{H}}\leq \alpha} d_{p}(F(K(m),\hat K( m )), F(K,\hat K))=0.
	\end{eqnarray}
	 That is to say proving  	\begin{eqnarray}\label{eq-3-41}
	\lim_{m,m'\rightarrow \infty} \mathop{\text{sup}}_{|(k,\hat k)|_{\mathcal{H}}\leq \alpha} D_{j,p}\big( (K(m),\hat K( m )) ,(K(m'),\hat K( m' )) \big) =0,
	\end{eqnarray}
	for $j=1,2$,
	and 
	\begin{eqnarray}\label{3-41}
	\mathop{\text{sup}}_{m\in \mathbb{N}}\mathop{\text{sup}}_{|(k,\hat k)|_{\mathcal{H}}\leq \alpha} D_{j,p}\big( K(m),\hat K( m ) \big) <\infty ,
		\end{eqnarray}
	with $j=1$. 
	
	For $j=1$, when $m\ge n$, it has $D_{\text{1,}p}\big( (K\left( m \right),\hat K\left( m \right)),(K,\hat K) \big) ^p =0$.
	Then turn to the case that $	m\le n$, 
	\begin{eqnarray}\label{eq-3-30}
	&&\mathop{\text{sup}}_{|(k,\hat k)|_{\mathcal{H}}\leq \alpha}D_{\text{1,}p}\big( (K\left( m \right),\hat K\left( m \right)),(K,\hat K) \big) ^{p} \cr
	&\leq& \mathop{\text{sup}}_{|(k,\hat k)|_{\mathcal{H}}\leq \alpha}\sum_{n=m+1}^{\infty}n^\kappa\sum_{l =1}^{2^n}{\left| 2^{m-n}\Delta _{{\tilde l}}^{m}(k,\hat k)^\mathrm{T}-\Delta _{l}^{n}(k,\hat k)^\mathrm{T} \right|}^{p}\cr
	&\leq& C_{p}\sum_{n=m+1}^{\infty}n^\kappa\sum_{l =1}^{2^n}{\left( 2^{\left( m-n \right) {p}}2^{-m{p}H}+2^{-{p}nH} \right)}\cr
	&\leq& C \alpha^{p} 2^{-m(p H-1-\beta_6)},
	\end{eqnarray}
	in the final line, take suitable $\beta_6>0$ such that $p H-1-\beta_6>0$. Hence, (\ref{eq-3-41}) holds for $j=1$. 
	
	In the same manner we can see that,
	\begin{eqnarray}\label{3-42}
	\mathop{\text{sup}}_{m\in \mathbb{N}}\mathop{\text{sup}}_{|(k,\hat k)|_{\mathcal{H}}\leq \alpha} D_{1,p}\big( K,\hat K \big) <\infty ,
	\end{eqnarray}
	Together with (\ref{eq-3-30}), we can prove that (\ref{3-41}) holds for $j=1$.
	
	For $j=2$, when $	m\le n$,
	\begin{eqnarray}\label{eq-3-31}
	&&\mathop{\text{sup}}_{|(k,\hat k)|_{\mathcal{H}}\leq \alpha}\big| (K( m+1 ),\hat K( m+1 )) _{t_{l-1}^{n},t_{l}^{n}}^{2}-(K( m ),\hat K( m )) _{t_{l-1}^{n},t_{l}^{n}}^{2} \big|\cr
	&\leq& \mathop{\text{sup}}_{|(k,\hat k)|_{\mathcal{H}}\leq \alpha}\big( \big| (K( m+1 ),\hat K( m+1 )) _{t_{l-1}^{n},t_{l}^{n}}^{2} \big|+\big| (K( m ),\hat K( m )) _{t_{l-1}^{n},t_{l}^{n}}^{2} \big| \big)\cr
	&\leq& C2^{2\left( m-n \right)}\big( \big| \big( \Delta _{{\tilde l}}^{m+1}(k,\hat k)^\mathrm{T} \big) ^{\otimes 2} \big|+\big| \big( \Delta _{{\tilde l}}^{m}(k,\hat k)^\mathrm{T} \big) ^{\otimes 2} \big| \big)\cr
	&\leq& C2^{-2nH}.
	\end{eqnarray}
	When $m\ge n$, it deduces that
	\begin{eqnarray*}\label{eq-3-32}
	&&(K( m+1 ),\hat K( m+1 )) _{t_{l-1}^{n},t_{l}^{n}}^{2}-(K( m ),\hat K( m )) _{t_{l-1}^{n},t_{l}^{n}}^{2}\cr
	&=&\frac{1}{2}\sum_{{\tilde l}=2^{m-n}\left( l-1 \right) +1}^{2^{m-n}l}{\left( \Delta _{2{\tilde l}-1}^{m+1}(k,\hat k)^\mathrm{T}\otimes \Delta _{2{\tilde l}}^{m+1}(k,\hat k)^\mathrm{T}-\Delta _{2{\tilde l}}^{m+1}(k,\hat k)^\mathrm{T}\otimes \Delta _{2{\tilde l}-1}^{m+1}(k,\hat k)^\mathrm{T} \right)}.
	\end{eqnarray*}
	Clearly, 
	\begin{eqnarray}\label{eq-3-33}
	\big| (K(m+1),\hat K(m+1)) _{t_{l-1}^{n},t_{l}^{n}}^{\text{2,}i,j}-(K(m),\hat K(m)) _{t_{l-1}^{n},t_{l}^{n}}^{\text{2,}i,j} \big|\leq C\big( \Xi_{n,m,l}^{i,j}+\Xi_{n,m,l}^{j,i} \big). 
	\end{eqnarray}
	where $i,j\in \{1,\cdots ,d_1+d_2\}$ and $i\ne j$,
	\begin{eqnarray*}\label{eq-3-34}
		\Xi_{n,m,l}^{i,j}&=& \bigg|\sum_{{\tilde l}=2^{m-n}\left( l-1 \right) +1}^{2^{m-n}l}{\Delta _{2{\tilde l}-1}^{m+1}(k,\hat k)^{\mathrm{T},i}\Delta _{2{\tilde l}}^{m+1}(k,\hat k)^{\mathrm{T},j}}\bigg|,
	\end{eqnarray*}
Therefore, due to the H\"older inequality, one can get 
\begin{eqnarray*}\label{eq-3-35}
	\Xi_{n,m,l}^{i,j}&\le&\big[ \sum_{{\tilde l}=2^{m-n}\left( l-1 \right) +1}^{2^{m-n}l}|\Delta _{2{\tilde l}-1}^{m+1}(k,\hat k)^{\mathrm{T},i}|^2\big]^{\frac{1}{2}}\big[ \sum_{{\tilde l}=2^{m-n}\left( l-1 \right) +1}^{2^{m-n}l}|\Delta _{2{\tilde l}}^{m+1}(k,\hat k)^{\mathrm{T},j}|^2\big]^{\frac{1}{2}}.
\end{eqnarray*}
With Lemma 5.1 in \cite{2007Friz}, it has $\sum_{{\tilde l}=2^{m-n}\left( l-1 \right) +1}^{2^{m-n}l}|\Delta _{2{\tilde l}-1}^{m+1}(k,\hat k)^{\mathrm{T},i}|^2 \le C 2^{-m( 2H-1/2)-n/2}$ for $i\in \{\text{1,}\cdots ,d_1\}$, and $\sum_{{\tilde l}=2^{m-n}\left( l-1 \right) +1}^{2^{m-n}l}|\Delta _{2{\tilde l}-1}^{m+1}(k,\hat k)^{\mathrm{T},i}|^2 \le C 2^{-m/2-n/2}$ for $i\in \{d_1+\text{1,}\cdots ,d_1+d_2\}$.
Then it deduces 
$$\Xi_{n,m,l}^{i,j}\leq C 2^{-m(2 H-1/2)-n/2},$$
and similar results hold for the $\Xi_{n,m,l}^{j,i}$.

    Together with (\ref{eq-3-31}), we  have
    \begin{eqnarray*}\label{eq-3-40}
    \mathop{\text{sup}}_{|(k,\hat k)|_{\mathcal{H}}\leq \alpha}D_{\text{2,}p}\left( (K( m+1 ),\hat K\left( m+1 \right)) ,(K( m ),\hat K\left( m \right))  \right) \leq C2^{-m\beta_7},
    \end{eqnarray*}
    where $\beta_7>0$. This gives (\ref{eq-3-29}) combined with (\ref{eq-3-30}).

	\textbf{Step 3}. Combine Step 1 and Step 2, with $\varepsilon \to 0$,  the law of $\varepsilon(B^H, W)$, $\left\{\mathbb{P}_{\varepsilon}^{H}\right\}_{\varepsilon>0}$ satisfies the Schilder-type LDP with the good rate function $I$ defined in (\ref{eq-3-2}) by means of an extension of the contraction principle [\cite{1998Dembo}, Theorem 4.2.23].

	The proof is completed.\qed
	
	\begin{lem}\label{prop3-2-3}Let $\big(B^{H},W\big)$ be the mixed geometric rough path with $H\in (1/3,1/2)$, there exists some constant $c>0$ such that
		\begin{eqnarray*}
			\mathbb{E}\left[\exp \left(c {\vertiii{\varepsilon(B^H,W)}^2_{p\mathrm{-var}}}\right)\right]=\int_{G \Omega_{p}\big(\mathbb{R}^{d_1+d_2}\big)} \exp \big(c {\vertiii{\varepsilon(B^H,W)}^2_{p\mathrm{-var}}}\big) \mathbb{P}^{H}\big(B^{H},W\big)<\infty.
		\end{eqnarray*}
	\end{lem}		
	\para{Proof}. 
	\propref{prop3-2-1} shows that, there is a sequence $c(m)\to 0$ such that for ${\varrho}> 1$, 
	\begin{eqnarray}\label{3-43}
	\mathbb{E}\big[D_{\text{1,}p}\left( (B^H(m), W( m )) ,(B^H, W) \right) ^{\varrho} \big]   ^{\text{1/}{\varrho}}\leq c\left( m \right) \sqrt{{\varrho}},	
	\end{eqnarray}
	and
	\begin{eqnarray}\label{3-44}
	\mathbb{E}\big[D_{1, p}(B^H(m),W(m))^{{\varrho}}\big]^{1 / {\varrho}} \leqslant C \sqrt{{\varrho}},
	\end{eqnarray}
	and
	\begin{eqnarray}\label{3-45}
	\mathbb{E}\big[ D_{\text{2,}p}\left( (B^H(m), W( m )) ,(B^H, W) \right) ^{\varrho} \big]  ^{\text{1/}{\varrho}}\leq c\left( m \right) {\varrho},	
	\end{eqnarray}
	where $c(m)$ tends to zero as $m \to \infty$.
	
	Then, it deduces that 
	\begin{eqnarray*}
		\mathbb{E}\big[e^{c D_{j, p}((B^{H}(m),W(m)), (B^{H},W))^{2 / j}}\big] & \leq& \sum_{N=0}^{\infty} e^{c(N+1)} \mathbb{P}\big(N<D_{j, p}((B^{H}(m),W(m)), (B^{H},W))^{2 / j} \leq N+1\big) \cr
		& \leq&\left(e^{c}+\cdots+c^{4 c}\right)\cr
		&&+e^{c} \sum_{N=4}^{\infty} e^{c N} \mathbb{P}\left(N<D_{j, p}((B^{H}(m),W(m)), (B^{H},W))^{2 / j}\right) \cr
		& \leq&\left(e^{c}+\cdots+c^{4 c}\right)\cr
		&&+e^{c} \sum_{N=4}^{\infty} e^{c N} N^{-N} \mathbb{E}\left[(D_{j, p}((B^{H}(m),W(m)), (B^{H},W))^{2 / j})^{N}\right] \cr
		& \leq&\left(e^{c}+\cdots+c^{4 c}\right)+e^{c} \sum^{\infty}_4 \exp \left[N\left(c+\log c_{2}+\log a_{m}\right)\right],
	\end{eqnarray*}
	for $j=1,2$.
	We choose $m_0$ large enough such that for any $m>m_0$, it has $\left[N\left(c+\log c_{2}+\log a_{m}\right)\right]<0$, moreover, it has
	\begin{eqnarray*}
		\sup _{m \geq m_{0}} \mathbb{E}\big[e^{c D_{j, p}((B^{H}(m),W(m)), (B^{H},W))^{2 / j}}\big] < \infty.
	\end{eqnarray*}
	Besides, it is clear to see that for any fixed $m_0$, there exists a constant $C(m_0)$ such that
	$
	D_{j, p}(B^H(m_0),W(m_0))^{1 / {j}} \leqslant C(m_0)\|(b^H,w)\|_{\infty}
	$. By the conventional Fernique theorem for Gaussian measures, it follows that $D_{j, p}(B^H(m_0),W(m_0))^{1 / {j}}$ is square exponential integrable. Then with the triangle inequality, the conclusion follows.
	
	The proof is completed.\qed

	Then, for simplicity, the RDE (\ref{sec-1}) is rewritten as follows,
	\begin{eqnarray}\label{3-11}
	dY_{t}^{\varepsilon}=\left[\sigma \left( Y_{t}^{\varepsilon} \right)|\hat\sigma \left( Y_{t}^{\varepsilon} \right)\right] \varepsilon d(B^H,W)_t+\beta \left(Y_{t}^{\varepsilon} \right) \lambda^\varepsilon dt, \quad Y_{0}^{\varepsilon}=\text{0.}
	\end{eqnarray}
	Denote  $\delta_{\lambda^\varepsilon}$ the law of $\lambda^\varepsilon_t$.

	\begin{prop}\label{prop3-2-2}
			With $\varepsilon \to 0$, the law of $Y^{\varepsilon,1}$ satisfies the LDP with the good rate function $I$,
			\begin{eqnarray}\label{3-1}
			I(y)= \begin{cases}\inf \left\{\|(k,\hat k)\|_{{\mathcal{H}}}^{2} / 2 \mid y=\hat{\Phi}_{0}((k,\hat k), \lambda)^{1}\right\} &  \big(\text{if }  y=\hat{\Phi}_{0}((k,\hat k), \lambda)^{1} \text{for some } (k,\hat k)^\mathrm{T} \in {\mathcal{H}}\big), \\ \infty & \text { (otherwise). }\end{cases}
			\end{eqnarray}
	\end{prop}	
	
	\para{Proof}. 
The \lemref{prop3-2-1} and \lemref{prop3-2-3} state that the law of the geometric rough path $\varepsilon(B^H,W)$, $\{\mathbb{P}_{\varepsilon}^{H} \}_{0<\varepsilon\le1}$ is exponential tight and satisfies the LDP with the good rate function (\ref{eq-3-2}). Obviously, the deterministic family $\{\delta_{\lambda^\varepsilon}\}_{0<\varepsilon\le 1}$ is exponential tight and satisfies LDP on $C_0^{1\mathrm{-var}} ([0,1], \mathbb{R})$ with the good rate function $+\infty\cdot\mathrm{I}_{{0}^c}$ with the convention that $+\infty\cdot 0=0$. 
	By the result for LDP of product measure [\cite{1998Dembo}, Exercise 4.2.7], the product measure $\{\mathbb{P}_{\varepsilon}^{H} \otimes\delta_{\lambda^\varepsilon}\}_{0<\varepsilon\le1}$ satisfies the LDP on the $G\Omega_p(\mathbb{R}^{d_1+d_2}) \times 
	C_0^{1\mathrm{-var}} ([0,1], \mathbb{R})$ with the good rate function as follows,
	\begin{eqnarray}\label{3-12}
	\hat{I}_{1}(B^H,W)+\infty \cdot \mathbf{1}_{\{0\}^{c}}(\lambda)= \begin{cases}\frac{1}{2}\|k,\hat k\|_{\mathcal{H}}^{2} & \big(\text { if } (B^H, W) \text { is lying above } (k,\hat k)^\mathrm{T} \in {\mathcal{H}}\text { and } \lambda=0\big), \\ \infty & (\text { otherwise }).\end{cases}
	\end{eqnarray}
	Clearly, the continuity theorem of It\^o map [\cite{2002Lyons}, Section 6.3] and the contraction principle yield that the above proposition holds.

	The proof is completed.\qed
	
\begin{rem}\label{rem3-2-2}
		Consider the following RDE,
		\begin{eqnarray}\label{3-13}
		d\tilde Y_{t}^{\varepsilon}=\big[\varepsilon^\upsilon\sigma \big( \tilde Y_{t}^{\varepsilon} \big)|\varepsilon^{\upsilon'}\hat\sigma \big( \tilde Y_{t}^{\varepsilon} \big)\big]  d(B^H,W)_t+\beta \big(\varepsilon,\tilde Y_{t}^{\varepsilon} \big) dt, \quad \tilde Y_{0}^{\varepsilon}=\text{0,}
		\end{eqnarray}
		where $\upsilon,\upsilon'> 0$.
		If $\upsilon<\upsilon'$, the above RDE (\ref{3-13}) can be rewritten as follows,
		\begin{eqnarray}\label{3-16}
		d\tilde Y_{t}^{\varepsilon}=\big[\varepsilon^{\upsilon-\upsilon'}\sigma \big( \tilde Y_{t}^{\varepsilon} \big)|\hat\sigma \big( \tilde Y_{t}^{\varepsilon} \big)\big] \varepsilon^{\upsilon'} d(B^H,W)_t+\beta \big(\varepsilon,\tilde Y_{t}^{\varepsilon} \big) dt, \quad \tilde Y_{0}^{\varepsilon}=\text{0.}
		\end{eqnarray}
		Let $\tilde Y^{\varepsilon}=\hat{\Psi}_{\varepsilon}(\varepsilon^{\upsilon'} (B^H, W), \lambda): G \Omega_{p}\left(\mathbb{R}^{d_1+d_2+1}\right) \mapsto G \Omega_{p}\left(\mathbb{R}^{n}\right)$ denote the It\^o map corresponding to (\ref{3-13}) with $\lambda_{t}=t$. 
		Consider the solution map $\hat{\Psi}_{\varepsilon}(\varepsilon^{\upsilon'} (B^H,W)+(\gamma, \eta), \lambda)$,
		\begin{eqnarray}\label{3-22}
		d\tilde Y_{t}^{\varepsilon}=\big[\varepsilon^{\upsilon-\upsilon'}\sigma \big( \tilde Y_{t}^{\varepsilon} \big)|\hat\sigma \big( \tilde Y_{t}^{\varepsilon} \big)\big] ( \varepsilon^{\upsilon'} d(B^H,W)_t+(\gamma_t, \eta_t)^\mathrm{T})+\beta \big(\varepsilon,\tilde Y_{t}^{\varepsilon} \big) dt, \quad \tilde Y_{0}^{\varepsilon}=\text{0}.
		\end{eqnarray}
		Now $\hat{\Psi}_{0}(\eta, \lambda)$ is the solution map lying above $\psi^{0}=\hat{\Psi}_{0}(\eta, \lambda) \in C_{0}^{q-\operatorname{var}}\left(\mathbb{R}^{n}\right)$ with $1< q<2$, satisfying that
		\begin{eqnarray}\label{3-23}
		d \psi_{t}^{0}=\hat \sigma \left(\psi_{t}^{0}\right) d \eta_t+\beta\left(0,\psi_{t}^{0}\right) d t, \quad \psi_{0}^{0}=0.
		\end{eqnarray}
		The \lemref{prop3-2-1} and \lemref{prop3-2-3} state that the law of the geometric rough path $\varepsilon^{\upsilon'}(B^H,W)$, $\{\mathbb{P}_{\varepsilon^{\upsilon'}}^{H} \}_{0<\varepsilon\le1}$ is exponential tight and satisfies the LDP with the good rate function (\ref{eq-3-2}). The deterministic family $\{\delta_{\lambda^\varepsilon}\}_{0<\varepsilon\le 1}$ is exponential tight and satisfies LDP on $C_0^{1\mathrm{-var}} ([0,1], \mathbb{R})$ with the good rate function $+\infty\cdot\mathrm{I}_{{0}^c}$ with the convention that $+\infty\cdot 0=0$.
		Similar to \lemref{prop3-2-2}, with aid of the  [\cite{1998Dembo}, Exercise 4.2.7], it deduces that  the product measure $\{\mathbb{P}_{\varepsilon^{\upsilon'}}^{H} \otimes\delta_{\lambda^\varepsilon}\}_{0<\varepsilon\le1}$ satisfies the LDP on the $G\Omega_p(\mathbb{R}^{d_1+d_2}) \times 
		C_0^{1\mathrm{-var}} ([0,1], \mathbb{R})$ with the good rate function as follows,
		\begin{eqnarray}\label{3-17}
		\hat{I}_{1}(B^H,W)+\infty \cdot \mathbf{1}_{\{0\}^{c}}(\lambda)= \begin{cases}\frac{1}{2}\|k,\hat k\|_{\mathcal{H}}^{2} & \big(\text { if } (B^H, W) \text { is lying above } (k,\hat k)^\mathrm{T} \in {\mathcal{H}}\text { and } \lambda=0\big), \\ \infty & (\text { otherwise }).\end{cases}
		\end{eqnarray}
		Clearly, the continuity theorem of It\^o map [\cite{2002Lyons}, Section 6.3], the contraction principle and the fact that $\upsilon-\upsilon'>0$ yield that  with $\varepsilon \to 0$,  the law of $\tilde Y^{\varepsilon,1}$ satisfies the LDP with the good rate function $I$,
		\begin{eqnarray}\label{3-14}
		I(\tilde y)= \begin{cases}\inf \left\{\|\hat k\|_{{\mathcal{H}^{H,d_1}}}^{2} / 2 \mid \tilde y=\hat{\Psi}_{0}(\hat k, \lambda)^{1}\right\} &  \big(\text{if }  \tilde y=\hat{\Psi}_{0}(\hat k, \lambda)^{1} \text{for some } \hat k \in {\mathcal{H}^{H,d_1}}\big), \\ \infty & \text { (otherwise). }\end{cases}
		\end{eqnarray}
		If $\upsilon<\upsilon'$, 
		the above RDE (\ref{3-13}) can be rewritten as follows,
		\begin{eqnarray}\label{3-10}
		d\tilde Y_{t}^{\varepsilon}=\big[\sigma \big( \tilde Y_{t}^{\varepsilon} \big)|\varepsilon^{\upsilon'-\upsilon}\hat\sigma \big( \tilde Y_{t}^{\varepsilon} \big)\big] \varepsilon^{\upsilon} d(B^H,W)_t+\beta \big(\varepsilon,\tilde Y_{t}^{\varepsilon} \big) dt, \quad \tilde Y_{0}^{\varepsilon}=\text{0.}
		\end{eqnarray}
		Let $\tilde Y^{\varepsilon}=\hat{\Psi}_{\varepsilon}(\varepsilon^{\upsilon} (B^H, W), \lambda): G \Omega_{p}\left(\mathbb{R}^{d_1+d_2+1}\right) \mapsto G \Omega_{p}\left(\mathbb{R}^{n}\right)$ denote the It\^o map corresponding to (\ref{3-10}) with $\lambda_{t}=t$.  Consider the solution map $\hat{\Psi}'_{\varepsilon}(\varepsilon^\upsilon (B^H,W)+(\gamma, \eta), \lambda)$,
		\begin{eqnarray}\label{3-24}
		d\tilde Y_{t}^{\varepsilon}=\big[\sigma \big( \tilde Y_{t}^{\varepsilon} \big)|\varepsilon^{\upsilon'-\upsilon}\hat\sigma \big( \tilde Y_{t}^{\varepsilon} \big)\big] ( \varepsilon^{\upsilon} d(B^H,W)_t+(\gamma_t, \eta_t)^\mathrm{T})+\beta \big(\varepsilon,\tilde Y_{t}^{\varepsilon} \big) dt, \quad \tilde Y_{0}^{\varepsilon}=\text{0}.
		\end{eqnarray}
		Now $\hat{\Psi}'_{0}(\gamma, \lambda)$ is the solution map lying above $\hat \psi^{0}=\hat{\Psi}'_{0}(\gamma, \lambda) \in C_{0}^{q-\operatorname{var}}\left(\mathbb{R}^{n}\right)$ with $1< q<2$, satisfying that
		\begin{eqnarray}\label{3-25}
		d \hat\psi_{t}^{0}=\sigma \big(\hat\psi_{t}^{0}\big) d \gamma_t+\beta\big(0,\hat\psi_{t}^{0}\big) d t, \quad \hat\psi_{0}^{0}=0.
		\end{eqnarray}
		with $\varepsilon \to 0$, take the same manner as above, it can deduces that the law of $\tilde Y^{\varepsilon,1}$ satisfies the LDP with the good rate function $I$,
		\begin{eqnarray}\label{3-15}
		I(\tilde y)= \begin{cases}\inf \left\{\|k\|_{{\mathcal{H}^{\frac{1}{2},d_2}}}^{2} / 2 \mid \tilde y=\hat{\Psi}'_{0}(k, \lambda)^{1}\right\} &  \big(\text{if }  \tilde y=\hat{\Psi}'_{0}(k, \lambda)^{1} \text{for some } k \in {\mathcal{H}^{\frac{1}{2},d_2}}\big), \\ \infty & \text { (otherwise). }\end{cases}
		\end{eqnarray}
		If $\upsilon=\upsilon'$, the resuts is shown in \lemref{prop3-2-2}.
\end{rem}

	Here follows the result for the Taylor expansion of $\phi^{(\varepsilon)}-\phi^{0}$ in (\ref{3-2-2}) with $m$th term $\phi^{m}$ around $(\gamma, \eta)^\mathrm{T} \in C_{0}^{q\mathrm{-var}}\left(\mathbb{R}^{d_1+d_2}\right)$ with $1 / p+1 / q>1$. 
	\begin{lem}\label{lem3-2-2}
Let $p\ge 2, 1\le q<2$, for any $m=1,2,\cdots,$ 
		\begin{eqnarray*}\label{rate}
        \phi^{(\varepsilon)}=\phi^{0}+\varepsilon \phi^{1}+\cdots+\varepsilon^{m} \phi^{m}+R_{\varepsilon}^{m+1}.
		\end{eqnarray*}
The maps $((B^H,W),(\gamma,\eta)^\mathrm{T})\in G \Omega_{p}\left(\mathbb{R}^{d_1+d_2}\right) \times C_{0}^{q\mathrm{-var}}\left(\mathbb{R}^{d_1+d_2}\right) \mapsto \phi^{k}$, $R_{\varepsilon}^{m+1} \in C_{0}^{p\mathrm{-var}}\left(\mathbb{R}^{n}\right)$ are continuous for $0\le  m'\le m$. Moreover, we have the following properties:

(i) For any $r_{1}>0$, there exists $C_1>0$ depending on $r_{1}$ such that if $\|(\gamma,\eta)\|_{q\mathrm{-var}} \leq r_{1}$, then $\left\|\phi^{\hat m}\right\|_{p  {\mathrm{-var} }} \leq C_{1}(1+{\vertiii{(B^H,W)}_{p\mathrm{-var}}})^{\hat m}$.

(ii) For any $r_{2},r_{3}>0$, there exists $C_2>0$ depending on $r_{2}$ and $r_{3}$ such that if $\|(\gamma,\eta)\|_{q\mathrm{-var}} \leq r_{2}$ and ${\vertiii{\varepsilon(B^H,W)}_{p\mathrm{-var}}}\leq r_{3}$, then $\left\|R_{\varepsilon}^{m+1}\right\|_{p\mathrm{-var}} \leq C_{2}(\varepsilon+{\vertiii{(B^H,W)}_{p\mathrm{-var}}})^{m+1}$.
	\end{lem}	
	
	\para{Proof}. 
	This lemma could be covered by reference \cite{2010Inahama}.\qed

	Now, we give the Cameron-Martin theorem for mixed geometric rough path.	
	\begin{thm}\label{thm3-1-1}
		{{(Cameron-Martin theorem for mixed rough path)}} For any $(k,\hat k)^\mathrm{T}\in{\mathcal{H}}$, the law of $(B^H, W)$, $\big\{\mathbb{P}_{\varepsilon}^{H}\big\}$ and $\big\{\mathbb{P}_{\varepsilon}^{H}((\cdot,\cdot)+(K,\hat K))\big\}$ are mutually absolutely continuous, moreover, for any bounded Borel function $f$ on $G \Omega_{p}\left(\mathbb{R}^{d_1+d_2}\right)$, 
		\begin{eqnarray}\label{eq-3-1}
		&\int_{G \Omega_{p}\left(\mathbb{R}^{d_1+d_2}\right)} f((B^H, W)+(K,\hat K)) \mathbb{P}_{\varepsilon}^{H}(d (B^H, W))\cr
		&=\int_{G \Omega_{p}\left(\mathbb{R}^{d_1+d_2}\right)} f(B^H, W) \exp \big(\frac{1}{\varepsilon}\big\langle (k,\hat k)^\mathrm{T}, (B^H, W)^{1}\big\rangle-\frac{1}{2 \varepsilon^{2}}\|(k,\hat k)\|_{\mathcal{H}}^{2}\big) \mathbb{P}_{\varepsilon}^{H}(d (B^H, W)).
		\end{eqnarray}
	\end{thm}	
	
	\para{Proof}. $(B^H(m), W( m ))+(K(m), \hat K( m ))$ is the smooth rough path constructed by $(b^H(m), w( m ))+(k(m), \hat k( m ))$. Meanwhile, as $m \to \infty$, $(B^H(m), W( m ))\to (B^H, W)$ in $G\Omega_p(\mathbb{R}^{d_1+d_2})$ and $(k(m),\hat k(m)) \to (k,\hat k)$ with $q\mathrm{-}$norm,  it is easy to see that (\ref{eq-3-1}) holds based on the Cameron-Martin theorem for  fBm and Bm $(b^H_t, w_t)_{t\ge 0}$.
	
	The proof is completed.\qed

	Next, we give the Fernique type theorem for the mixed geometric rough path for later use.

	\section{Computation of  Hessian}\label{sec-4}
	In this section, we set conditions for parameters. Firstly, for the fBm, the Hurst parameter $H\in(1/3,1/2)$. Then, choose $p$ and $q$ such that 
		\begin{eqnarray}\label{4-0}
&&\frac{1}{p'} \vee \frac{1}{[1 / H]+1}<\frac{1}{p}<H, 1<\frac{1}{q}<H+\frac{1}{2}\cr
&&\frac{1}{p}+\frac{1}{q}>1, \quad \frac{1}{q}-\frac{1}{p}>\frac{1}{2} .
		\end{eqnarray}
	For example,  choose that $1 / p=H-2\varepsilon$ and $1 / q=H+1 / 2-\varepsilon$ for small parameter $\varepsilon$, then the above condition (\ref{4-0}) is satisfied.
		\subsection{Hilbert-Schmidt property of  Hessian}\label{sec-41}
		In this subsection, we show the Hilbert-Schmidt property for the Hessian matrix  of the It\^o map restricted on the Cameron-Martin space $
		{\mathcal{H}}$. Throught this section, set $\beta_0(y)=\beta(0,y)$.
		
		For fixed $(\gamma,\eta)^\mathrm{T}\in C_0^{p\mathrm{-var}}(\mathbb{R}^{d_1 +d_2})$, set
		\begin{eqnarray}\label{4-1}
		d \Omega_{t}=[\sigma(\phi_{t}^{0})|\hat \sigma(\phi_{t}^{0})]\left\langle\cdot, d (\gamma_{t},\eta_t)^\mathrm{T}\right\rangle+\nabla \beta_0\left(\phi_{t}^{0}\right)\langle\cdot,1\rangle d t,
		\end{eqnarray}
		where $\phi_{0}$ is the solution to the ODE (\ref{3-2-3}). And $M_t$ satisfies the ODE  as follows,
		\begin{eqnarray*}\label{4-81}
			dM_t=d \Omega_{t}\cdot M_t,\quad M_{0}=\operatorname{Id}_{n},
		\end{eqnarray*}
	where $\operatorname{Id}_{n}$ is an identity matrix in $\mathbb{R}^{n\times n}$.
	
		Similarly, its inverse exists and satisfies the ODE,
		\begin{eqnarray*}\label{4-82}
			d M_{t}^{-1}=-M_{t}^{-1} \cdot d \Omega_{t}, \quad M_{0}^{-1}=\operatorname{Id}_{n}.
		\end{eqnarray*}
		More details about $M_t$, $M_t^{-1}$ and $\Omega_{t}$ refer to \cite{2010Inahama}.
		
		Set $\chi(k,\hat k)=(\nabla \Psi)(\gamma,\eta)\langle (k,\hat k)\rangle \in C^{p\mathrm{-var}}_0(\mathbb{R}^n)$, which satisfies the following ODE,
		\begin{eqnarray}\label{4-2}
		d \chi_{t}-[\nabla \sigma(\phi_{t}^{0})|\nabla\hat \sigma(\phi_{t}^{0})]\left\langle\chi_{t}, d (\gamma_t,\eta_t)^\mathrm{T}\right\rangle-\nabla \beta_0(\phi_{t}^{0})\left\langle\chi_{t},1\right\rangle dt=[\sigma(\phi_{t}^{0})|\hat \sigma(\phi_{t}^{0})] d (k_{t},\hat k_t)^\mathrm{T}, \quad \chi_{0}=0.
		\end{eqnarray}
		With aid of the results in  \cite{2010Inahama}, the solution could be written explicitly as follows,
		\begin{eqnarray}\label{4-3}
		\chi(k,\hat k)_{t}=(\nabla \Psi)(\gamma,\eta)\langle (k,\hat k)\rangle_{t}=M_{t} \int_{0}^{t} M_{s}^{-1} [\sigma(\phi_{s}^{0})|\hat \sigma(\phi_{s}^{0})] d (k_{s},\hat k_s)^\mathrm{T}.
		\end{eqnarray}
		According to the Young integral theory, the above integral can be defined in the Young integral sense. And $(k,\hat k)^\mathrm{T} \in C^{p\mathrm{-var}}_0(\mathbb{R}^{d_1 +d_2})\mapsto \chi(k,\hat k)\in C^{p\mathrm{-var}}_0(\mathbb{R}^n)$ is a continuous map.	
		
		Likewise, $\psi^{2}((k,\hat k), (\gamma,\eta))=		\nabla ^2\Psi (\gamma,\eta) \big< (f,\hat f),(k,\hat k) \big> _t \in C^{p\mathrm{-var}}_0(\mathbb{R}^n)$ satisfies the following ODE,
		\begin{eqnarray}\label{4-4}
		&&d\psi _t- [\nabla \sigma(\phi_{t}^{0})|\nabla\hat \sigma(\phi_{t}^{0})] \langle \psi _t,d(\gamma_t,\eta_t)^\mathrm{T} \rangle -\nabla \beta_0\left( \phi _{t}^{0} \right) \langle \psi _t,1 \rangle dt\cr
		&=&2[\nabla \sigma(\phi_{t}^{0})|\nabla\hat \sigma(\phi_{t}^{0})] \big< \chi(k,\hat k)_t,d(k_t,\hat k_t)^\mathrm{T} \big> \cr
		&&+[\nabla^2 \sigma(\phi_{t}^{0})|\nabla^2\hat \sigma(\phi_{t}^{0})] \big< \chi(k,\hat k)_t,\chi(k,\hat k)_t,d(\gamma_t,\eta_t)^\mathrm{T} \big>\cr
		&&+\nabla ^2\beta_0\left( \phi _{t}^{0} \right) \big< \chi(k,\hat k) _t,\chi(k,\hat k) _t \big> dt,\quad
		\psi _0=0.
		\end{eqnarray}
		Then, its solution can be written in the following sense,
		\begin{eqnarray}\label{4-5}
		\nabla ^2\Psi \left( \gamma \right) \big< (f,\hat f),(k,\hat k) \big> _t
		&=&M_t\int_0^t{M}_{s}^{-1}[\nabla \sigma(\phi_{s}^{0})|\nabla\hat \sigma(\phi_{s}^{0})] \big< \chi (f,\hat f) _s,d(k_s,\hat k_s)^\mathrm{T} \big> \cr
		&&\qquad +[\nabla \sigma(\phi_{s}^{0})|\nabla\hat \sigma(\phi_{s}^{0})] \big< \chi (k,\hat k) _s,d(f_s,\hat f_s)^\mathrm{T} \big>\cr
		&&+M_t\int_0^t{M}_{s}^{-1}[\nabla^2 \sigma(\phi_{s}^{0})|\nabla^2\hat \sigma(\phi_{s}^{0})] \big< \chi (f,\hat f) _s,\chi (k,\hat k)_s,d(\gamma,\eta)^\mathrm{T} _s \big>\cr
		&&\qquad +\nabla ^2\beta\left( \phi _{s}^{0} \right) \big< \chi (f,\hat f)_s,\chi(k,\hat k) _s \big> ds\cr
		&=:&V_1\big( (f,\hat f),(k,\hat k) \big) _t+V_2\big( (f,\hat f),(k,\hat k) \big) _t.
		\end{eqnarray}
		Next, define $\nabla^{2}(F \circ \Psi)(\gamma,\eta)\langle(\cdot,\cdot), (\cdot,\cdot)\rangle$  as follows,
		\begin{eqnarray*}
			\nabla^{2}(F \circ \Psi)(\gamma,\eta)\langle(\cdot,\cdot), (\cdot,\cdot)\rangle &=&\nabla^{2}(F (\Psi(\gamma,\eta)) )\big\langle\nabla^2\Psi(\gamma,\eta)\langle(f,\hat f),(k,\hat k)\rangle\big\rangle\cr
			&&+\nabla(F (\Psi(\gamma,\eta)) )\big\langle\nabla\Psi(\gamma,\eta)\langle(f,\hat f)\rangle,\nabla\Psi(\gamma,\eta)\langle(k,\hat k)\rangle\big\rangle,
		\end{eqnarray*}
		where $\nabla\Psi(\gamma,\eta)\langle(f,\hat f)$ and $\nabla^2\Psi(\gamma,\eta)\langle(f,\hat f),(k,\hat k)\rangle$ are defined in (\ref{4-3}) and (\ref{4-5}) respectively. 
		\begin{thm}\label{thm4-1-1} (Hilbert-Schmidt property of Hessian)
			 $\nabla^{2}(F \circ \Psi)(\gamma,\eta)\langle(\cdot,\cdot), (\cdot,\cdot)\rangle$ is a symmetric Hilbert-Schmidt bilinear form on $
			{\mathcal{H}}$ for all $(\gamma,\eta)^\mathrm{T} 
			\in{\mathcal{H}}$.
			\end{thm}

   We have divided the proof for \thmref{thm4-1-1} into a sequence of lemmas.
		
	\begin{lem}\label{thm4-1-2} 
		Choose appropriate $p$ and $q$. Then,  
			if $p'\ge p$, $\nabla F\left(\phi^{0}\right) \circ V_{2}$ is of trace class for a Fr\'echet differentiable function $F$.  A similar fact holds for $ \nabla^2 F\left(\phi^{0}\right)\big< \chi(k,\hat k) _t,\chi(k,\hat k) _t \big> $.
	\end{lem}
	\para{Proof}. 
According to the Young integral, $V_{2}$ can be defined in the Young integral sense. Then, applying \lemref{lem3-2-2}, it follows that $|| M_{u}^{-1} [\sigma(\phi_{u}^{0})|\hat \sigma(\phi_{u}^{0})]||_{q\mathrm{-var}}<\infty$. With aid of the Young integral, 
\begin{eqnarray*}\label{4-80}
		\big\|\chi(k,\hat k)_{.}\big\|_{q\mathrm{-var}} &\leq& c_{1}\left\|M_{.}^{-1} [\sigma(\phi_{.}^{0})|\hat \sigma(\phi_{\cdot}^{0})]\right\|_{q\mathrm{-var}}\|(k,\hat k)\|_{p\mathrm{-var}} \cr
		&\leq& c_{2}\|(k,\hat k)\|_{p\mathrm{-var}},
\end{eqnarray*}
for some constants $c_1, c_2$. Hence, $(k,\hat k)^\mathrm{T}\in C^{p\mathrm{-var}}_0(\mathbb{R}^{d_1+d_2} )\longmapsto\chi(k,\hat k)\in C^{p\mathrm{-var}}_0(\mathbb{R}^{n} )$ is a bounded linear map.

 Moreover, $((f,\hat f)^\mathrm{T},(k,\hat k)^\mathrm{T})\in C^{p\mathrm{-var}}_0(\mathbb{R}^{d_1+d_2} )\times C^{p\mathrm{-var}}_0(\mathbb{R}^{d_1+d_2} ) \longmapsto V_2\in C^{p\mathrm{-var}}_0(\mathbb{R}^{n} )$ is a bounded bilinear map. Due to the fact that the space $C^{p\mathrm{-var}}_0(\mathbb{R}^{d_1+d_2} )$ is not separable, consider the abstract Wiener space $(\mathcal{X},
{\mathcal{H}},\mu^H)$. The Goodman's theorem (Theorem 4.6, [\cite{1997Janson}]) shows that, its restriction on the Cameron-Martin space ${\mathcal{H}}$ is of trace class.
	The same reasoning applies to the case of   $ \nabla^2 F\left(\phi^{0}\right)\big< \chi(k,\hat k) _t,\chi(k,\hat k) _t \big> $.
	
	The proof is completed. \qed

	Our next concern will be  $V_1\big< (f,\hat f),(k,\hat k) \big>$, and it can be rewritten as follows,
  $$V_1\big< (f,\hat f),(k,\hat k) \big> =R_1\big< (f,\hat f),(k,\hat k) \big> +R_1\big< (k,\hat k),(f,\hat f) \big> -\big( R_2\big< (f,\hat f),(k,\hat k) \big> +R_2\big< (k,\hat k),(f,\hat f) \big> \big), $$
  where
  \begin{eqnarray}\label{4-6}
  R_1\big< (f,\hat f),(k,\hat k) \big> _t=M_t\int_0^t{M}_{s}^{-1}[\nabla \sigma(\phi_{s}^{0})|\nabla\hat \sigma(\phi_{s}^{0})] \big< [\sigma(\phi_{s}^{0})|\hat \sigma(\phi_{s}^{0})] (f_s,\hat f_s)^\mathrm{T},d(k_s,\hat k_s)^\mathrm{T} \big>,
   \end{eqnarray}
   and
  \begin{eqnarray}\label{4-7}
  R_2\big< (f,\hat f),(k,\hat k) \big> _t=M_t\int_0^t{M}_{s}^{-1}[\nabla \sigma(\phi_{s}^{0})|\nabla\hat \sigma(\phi_{s}^{0})] \big< M_s\int_0^s{d}\left[ M_{u}^{-1}[\sigma(\phi_{u}^{0})|\hat \sigma(\phi_{u}^{0})] \right] (f_u,\hat f_u)^\mathrm{T},d(k_s,\hat k_s)^\mathrm{T} \big>.
  \end{eqnarray}

   \begin{lem}\label{thm4-1-3} 
	As a bilinear form on Cameron-Martin space $
	{\mathcal{H}}$, $\nabla F\left(\phi^{0}\right)  \circ R_2$ is of trace class. Moreover, if $\vartheta_l$ is $weak \star$ convergent to $\vartheta$
	as $l \to \infty$ in $C_0^{p\mathrm{-var}}(\mathbb{R}^n) ^\star$
	, then $\vartheta_l \circ R_2$ converges to $\vartheta \circ R_2$ as  $l \to \infty$ in the
	Hilbert-Schmidt norm.
    \end{lem}
	\para{Proof}.
	Applying \lemref{lem3-2-2}, it deduces that $|| M_{.}^{-1} [\sigma(\phi_{.}^{0})|\hat \sigma(\phi_{.}^{0})]||_{q\mathrm{-var}}<\infty$ and $|| M_{.}^{-1} [\nabla\sigma(\phi_{.}^{0})|\nabla\hat \sigma(\phi_{.}^{0})]||_{q\mathrm{-var}}<\infty$. Therefore
	\begin{eqnarray*}\label{4-8}
	\left\|\int^{\cdot}_{0} d\left[M_{u}^{-1} [\sigma(\phi_{u}^{0})|\hat \sigma(\phi_{u}^{0})]\right] (f_u,\hat f_u)^\mathrm{T}\right\|_{q\mathrm{-var}} &\leq& c_{1}\left\|M_{.}^{-1} [\sigma(\phi_{.}^{0})|\hat \sigma(\phi_{\cdot}^{0})]\right\|_{q\mathrm{-var}}\|(f,\hat f)\|_{p\mathrm{-var}}\cr
	 &\leq& c_{2}\|(f,\hat f)\|_{p\mathrm{-var}},
	\end{eqnarray*}
for some constant $c_2$.
	Moreover,
	\begin{eqnarray*}\label{4-9}
	&&\bigg\|M . \int_{0} M_{s}^{-1} [\nabla \sigma(\phi_{s}^{0})|\nabla\hat \sigma(\phi_{s}^{0})]\bigg\langle M_{s} \int_{0}^{.} d\left[M_{u}^{-1} [\sigma(\phi_{u}^{0})|\hat \sigma(\phi_{u}^{0})]\right] (f_u,\hat f_u)^\mathrm{T},d(k_s,\hat k_s)^\mathrm{T}\bigg\rangle\bigg\|_{p\mathrm{-var}} \cr
	& \leq& c_{3}\|(f,\hat f)\|_{p\mathrm{-var}}\|(k,\hat k)\|_{p\mathrm{-var}} ,
	\end{eqnarray*}
for some constant $c_3$. 
Obviously, $((f,\hat f)^\mathrm{T},(k,\hat k)^\mathrm{T})\in C^{p\mathrm{-var}}_0(\mathbb{R}^{d_1+d_2} )\times C^{p\mathrm{-var}}_0(\mathbb{R}^{d_1+d_2} )\longmapsto R_2\in C^{p\mathrm{-var}}_0(\mathbb{R}^{n} )$ is a bounded bilinear map. In particular, $\vartheta \circ R_2$ is a bounded bilinear form on $C^{p\mathrm{-var}}_0(\mathbb{R}^{n} )$. The Goodman’s theorem (Theorem 4.6, [\cite{1997Janson}]) shows that , its restriction on the Cameron-Martin space ${\mathcal{H}}$ is of trace class.
	
	Next, we intend  to prove the second argument.
	Since   ${\mathcal{H}^{H,d_1}} \hookrightarrow W_{0}^{\delta, 2} \cong L_{0, \text { real }}^{\delta, 2}$ [\cite{2013Inahama}, Proposition 3.4], it has  ${\mathcal{H}} \hookrightarrow  L_{0, \text { real }}^{\delta, 2}(\mathbb{R}^{d_1})\oplus{\mathcal{H}^{\frac{1}{2},d_2}}$.
	Then choose the ONB $(f,\hat f)^\mathrm{T}_{0,i}=1\cdot \mathbf{e}_i$,$(f,\hat f)^\mathrm{T}_{m,i}= \sqrt{2}\frac{\cos \left( m\pi t \right)}{\left( 1+m^2 \right) ^{\text{1/}2q}}\mathbf{e}_i$ for $i=\text{1,}\cdots ,d_1$,  $(f,\hat f)^\mathrm{T}_{0,i}= t\cdot \mathbf{e}_i$,  $(f,\hat f)^\mathrm{T}_{2m-1,i}= \frac{\cos \left( 2m\pi t \right)-1}{2m\pi}\cdot \mathbf{e}_i$ and $(f,\hat f)^\mathrm{T}_{2m,i}= \frac{\sin \left( 2m\pi t \right)}{2m\pi}\cdot \mathbf{e}_i$ for $i=d_1+1,\cdots ,d_1+d_2$.
	
	Fix $i=1,\cdots,d_1$, $(f,\hat f)^\mathrm{T}_{m,i}= \sqrt{2}\frac{\cos \left( m\pi t \right)}{\left( 1+m^2 \right) ^{\text{1/}2q}}$. For the jog-free function, it is easy to check that $\left\|\cos (n \pi t)-1\right\|_{p  {\mathrm{-var} }}=2 n^{1 / p}$\cite{2013Inahama}. Therefore,
	\begin{eqnarray}\label{4-10}
	\big\|(f,\hat f)^\mathrm{T}_{m, i}\big\|_{p-\operatorname{var}} \leq\left(1+m^{2}\right)^{-\delta / 2} \sqrt{2}\left(1+2 m^{1 / p}\right) \leq c\left(\frac{1}{1+m}\right)^{1 / q-1 / p},
	\end{eqnarray}
	due to  $2(1 / q-1 / p)>1$, it leads to 
	\begin{eqnarray*}\label{4-11}
	\sum_{i, j=1}^{d_1} \sum_{m, m'=0}^{\infty}\big\|R_{2}\big\langle (f,\hat f)_{m, i}, (f,\hat f)_{m', j}\big\rangle\big\|_{p\mathrm{-var}}^{2}
	& \leq &c \sum_{i, j=1}^{d} \sum_{m, m'=0}^{\infty}\big\|(f,\hat f)_{m, i}\big\|_{p\mathrm{-var}}^{2}\big\|(f,\hat f)_{m', j}\big\|_{p\mathrm{-var}}^{2}\cr
	& \leq & c \sum_{m=0}^{\infty}\left(\frac{1}{1+m}\right)^{2(1 / q-1 / p)} \sum_{m'=0}^{\infty}\left(\frac{1}{1+m'}\right)^{2(1 / q-1 / p)}\cr
	&<&\infty.
	\end{eqnarray*}
	Fix $i=d_1+1,\cdots,d_1+d_2$.  Likewise, it has $\left\|\sin (n \pi t)\right\|_{p  {\mathrm{-var} }}\le c (1/2^p +n)^{1 / p}$. Then
	\begin{eqnarray}\label{4-12}
	\big\|(f,\hat f)^\mathrm{T}_{2m, i}\big\|_{p-\operatorname{var}} \leq c(2m\pi)^{- 1} (1/2^p +m)^{1 / p} \leq c\left(\frac{1}{1+m}\right)^{1-1 / p},
	\end{eqnarray}
	Fix $i=d_1+1,\cdots,d_1+d_2,j=1,\cdots,d_1$, since $2(1 -1 / p)>1$ and $2(1 / q-1 / p)>1$, we have
	\begin{eqnarray*}\label{4-13}
	&&\sum_{i=d_1+1}^{d_1+d_2}	\sum_{j=1}^{d_1} \sum_{m, m'=0}^{\infty}\big\|R_{2}\big\langle (f,\hat f)_{2m, i}, (f,\hat f)_{m', j}\big\rangle\big\|_{p\mathrm{-var}}^{2}\cr
	&&+	\sum_{i=d_1+1}^{d_1+d_2}	\sum_{j=1}^{d_1} \sum_{m, m'=0}^{\infty}\big\|R_{2}\big\langle (f,\hat f)_{2m-1, i}, (f,\hat f)_{m', j}\big\rangle\big\|_{p\mathrm{-var}}^{2}\cr
	& \leq &c \sum_{i=d_1+1}^{d_1+d_2}	\sum_{j=1}^{d_1}\sum_{m, m'=0}^{\infty}\big\|(f,\hat f)_{2m, i}\big\|_{p\mathrm{-var}}^{2}\big\|(f,\hat f)_{m', j}\big\|_{p\mathrm{-var}}^{2}\cr
	&&+c\sum_{i=d_1+1}^{d_1+d_2}	\sum_{j=1}^{d_1} \sum_{m, m'=0}^{\infty}\big\|(f,\hat f)_{2m-1, i}\big\|_{p\mathrm{-var}}^{2}\big\|(f,\hat f)_{m', j}\big\|_{p\mathrm{-var}}^{2}\cr
	& \leq & c \sum_{m=0}^{\infty}\left(\frac{1}{1+m}\right)^{2(1 -1 / p)} \sum_{m'=0}^{\infty}\left(\frac{1}{1+m'}\right)^{2(1 / q-1 / p)}\cr
	&<&\infty.
	\end{eqnarray*}
	Fix $i=1,\cdots,d_1,j=d_1,\cdots,d_1+d_2$, in the same manner,  it deduces that $$\sum_{i=1}^{d_1}\sum_{j=d_1+1}^{d_1+d_2}	 \sum_{m, m'=0}^{\infty}\big\|R_{2}\big\langle (f,\hat f)_{m, i}, (f,\hat f)_{2m'-1, j}\big\rangle\big\|_{p\mathrm{-var}}^{2}<\infty,$$
	and 
	$$\sum_{i=1}^{d_1}\sum_{j=d_1+1}^{d_1+d_2} \sum_{m, m'=0}^{\infty}\big\|R_{2}\big\langle (f,\hat f)_{m, i}, (f,\hat f)_{2m', j}\big\rangle\big\|_{p\mathrm{-var}}^{2}<\infty.$$
	
	Fix $i=d_1+1,\cdots,d_1+d_2,j=d_1+1,\cdots,d_1+d_2$,
	since $2(1 -1 / p)>1$, we have
		\begin{eqnarray*}\label{4-14}
		&&\sum_{i=d_1+1}^{d_1+d_2}	\sum_{j=d_1+1}^{d_1+d_2} \sum_{m, m'=0}^{\infty}\big\|R_{2}\big\langle (f,\hat f)_{2m-1, i}, (f,\hat f)_{2m'-1, j}\big\rangle\big\|_{p\mathrm{-var}}^{2}\cr
		&&+	\sum_{i=d_1+1}^{d_1+d_2}	\sum_{j=d_1+1}^{d_1+d_2} \sum_{m, m'=0}^{\infty}\big\|R_{2}\big\langle (f,\hat f)_{2m-1, i}, (f,\hat f)_{2m', j}\big\rangle\big\|_{p\mathrm{-var}}^{2}\cr
		&&+		\sum_{i=d_1+1}^{d_1+d_2}	\sum_{j=d_1+1}^{d_1+d_2} \sum_{m, m'=0}^{\infty}\big\|R_{2}\big\langle (f,\hat f)_{2m, i}, (f,\hat f)_{2m'-1, j}\big\rangle\big\|_{p\mathrm{-var}}^{2}\cr
		&&+	\sum_{i=d_1+1}^{d_1+d_2}	\sum_{j=d_1+1}^{d_1+d_2} \sum_{m, m'=0}^{\infty}\big\|R_{2}\big\langle (f,\hat f)_{2m, i}, (f,\hat f)_{2m', j}\big\rangle\big\|_{p\mathrm{-var}}^{2}\cr
		& \leq &c \sum_{i=d_1+1}^{d_1+d_2}	\sum_{j=d_1+1}^{d_1+d_2}\sum_{m, m'=0}^{\infty}\big\|(f,\hat f)_{2m-1, i}\big\|_{p\mathrm{-var}}^{2}\big\|(f,\hat f)_{2m'-1, j}\big\|_{p\mathrm{-var}}^{2}\cr
		&&+c\sum_{i=d_1+1}^{d_1+d_2}	\sum_{j=d_1+1}^{d_1+d_2} \sum_{m, m'=0}^{\infty}\big\|(f,\hat f)_{2m-1, i}\big\|_{p\mathrm{-var}}^{2}\big\|(f,\hat f)_{2m', j}\big\|_{p\mathrm{-var}}^{2}\cr
		&&+c \sum_{i=d_1+1}^{d_1+d_2}	\sum_{j=d_1+1}^{d_1+d_2}\sum_{m, m'=0}^{\infty}\big\|(f,\hat f)_{2m, i}\big\|_{p\mathrm{-var}}^{2}\big\|(f,\hat f)_{2m'-1, j}\big\|_{p\mathrm{-var}}^{2}\cr
		&&+c\sum_{i=d_1+1}^{d_1+d_2}	\sum_{j=d_1+1}^{d_1+d_2} \sum_{m, m'=0}^{\infty}\big\|(f,\hat f)_{2m, i}\big\|_{p\mathrm{-var}}^{2}\big\|(f,\hat f)_{2m', j}\big\|_{p\mathrm{-var}}^{2}\cr
		&\leq & c \sum_{m=0}^{\infty}\left(\frac{1}{1+m}\right)^{2(1 -1 / p)} \sum_{m'=0}^{\infty}\left(\frac{1}{1+m'}\right)^{2(1 -1 / p)}
		<\infty.
		\end{eqnarray*}

	From the Banach-Steinhaus theorem, $\left\|\vartheta_{l}-\vartheta\right\|_{C^{p\mathrm{-var}, *}} \leq c$, then we have
	\begin{eqnarray}\label{4-15}
	\big|\left(\vartheta_{l}-\vartheta\right) \circ R_{2}\big\langle (f,\hat f)_{m, i}, (f,\hat f)_{m', j}\big\rangle\big|^{2} &\leq& c^{2}\big\|R_{2}\big\langle (f,\hat f)_{m, i}, (f,\hat f)_{m', j}\big\rangle\big\|_{p{\mathrm{-var}}}^{2},
	\end{eqnarray}
	similar results could be obtained for $	i=1,\cdots,d_1, j=d_1+1,\cdots,d_1+d_2$,  $	i=d_1+1,\cdots,d_1+d_2, j=1,\cdots,d_1$ and $	i=d_1+1,\cdots,d_1+d_2, j=d_1+1,\cdots,d_1+d_2$, and  omited  for simiplicity.
	
	By the dominated convergence theorem, it is clear to see that
	\begin{eqnarray}\label{4-16}
	\left\|\vartheta_{l} \circ R_{2}-\vartheta \circ R_{2}\right\|_{{\mathrm{HS}}\mathrm{-}\mathcal{H}^{H,d_1}\oplus\mathcal{H}^{\frac{1}{2},d_1}} \leq\|\iota\|_{\mathrm{op}}\left\|\iota^{*}\right\|_{\mathrm{op}}\left\|\vartheta_{l} \circ R_{2}-\vartheta \circ R_{2}\right\|_{\mathrm{HS}\mathrm{-}L_{\mathrm{real}}^{\delta, 2}(\mathbb{R}^{d_1})\oplus\mathcal{H}^{\frac{1}{2},d_2}} \rightarrow 0,
	\end{eqnarray}
	with $l \to 0$.
Hence, if $\vartheta_l$ is $weak \star$ convergent to $\vartheta$
as $l \to \infty$ in $C_0^{p\mathrm{-var}}(\mathbb{R}^n) ^\star$
, then $\vartheta_l \circ R_2$ converges to $\vartheta \circ R_2$ as  $l \to \infty$ in the
Hilbert-Schmidt norm.	
	The proof is completed. \qed	
	
	\begin{lem}\label{thm4-1-4} 
		As a bilinear form on Cameron-Martin space $
		{\mathcal{H}}$, $\nabla F\left(\phi^{0}\right)  \circ R_1$ is Hilbert-Schmidt. Moreover, if $\vartheta_l$ is $weak \star$ convergent to $\vartheta$
		as $l \to \infty$ in $C_0^{p\mathrm{-var}}(\mathbb{R}^n) ^\star$, then $\vartheta_l \circ R_1$ converges to $\vartheta \circ R_1$ as  $l \to \infty$ in the
		Hilbert-Schmidt norm.
	\end{lem}	
	\para{Proof}. Similar to \lemref{thm4-1-3}, it is sufficient to make the following observation,
	\begin{eqnarray}\label{4-17}
\sum_{i=1}^{d_1}	\sum_{j=1}^{d_1} \sum_{m, m'=0}^{\infty}\big\|R_{1}\big\langle (f,\hat f)_{m, i}, (f,\hat f)_{m', j}\big\rangle\big\|_{p\mathrm{-var}}^{2}<\infty,
	\end{eqnarray}
	and 	similar  for $	i=1,\cdots,d_1, j=d_1,\cdots,d_1+d_2$,  $	i=d_1,\cdots,d_1+d_2, j=1,\cdots,d_1$ and $	i=d_1,\cdots,d_1+d_2, j=d_1,\cdots,d_1+d_2$.
	
	We divide the rest of the proof into four steps.
	
	\textbf{Step 1}. Fix $i=\text{1,}\cdots ,d_1, j=\text{1,}\cdots ,d_1$.
	
		Now,  $(f,\hat f)^\mathrm{T}_{m,i}= \sqrt{2}\frac{\cos \left( m\pi t \right)}{\left( 1+m^2 \right) ^{\text{1/}2q}},  (f,\hat f)^\mathrm{T}_{m',j}= \sqrt{2}\frac{\cos \left( m'\pi t \right)}{\left( 1+m'^2 \right) ^{\text{1/}2q}}$,  $m\ne m'$,
	\begin{eqnarray*}\label{4-18}
	\sqrt{2}\cos \left( m\pi t \right) d\big[ \sqrt{2}\cos \left( m'\pi t \right) \big] 
	&=& m'd\left[ \frac{\cos ( ( m'+m ) \pi t )}{m'+m}+\frac{\cos ( ( m'-m ) \pi t )}{m'-m} \right].
	\end{eqnarray*}
	Then,
	\begin{eqnarray*}\label{4-19}
	\left. R_1\big< (f,\hat f)_{m,i}, \right. (f,\hat f)_{m',j} \big> _t
	&=&M_t\int_0^t{M}_{s}^{-1}[\nabla \sigma(\phi_{s}^{0})|\nabla\hat \sigma(\phi_{s}^{0})] \langle [\sigma(\phi_{s}^{0})|\hat \sigma(\phi_{s}^{0})] \mathbf{e}_i,\mathbf{e}_{j} \rangle \sqrt{2}\frac{\cos \left( m\pi t \right)}{\left( 1+m^2 \right) ^{\text{1/}2q}}\cr
    &&\quad\times d\left[ \frac{\sqrt{2}\cos \left( m'\pi s \right)}{\left( 1+m^{'2} \right) ^{\text{1/}2q}} \right]\cr
	&=&\frac{m'}{\left( 1+m^2 \right) ^{\text{1/}2q}\left( 1+m^{'2} \right) ^{\text{1/}2q}\left( m'+m \right)}\cr
	&&\quad\times M_t\int_0^t{M}_{s}^{-1}[\nabla \sigma(\phi_{s}^{0})|\nabla\hat \sigma(\phi_{s}^{0})] \langle [\sigma(\phi_{s}^{0})|\hat \sigma(\phi_{s}^{0})] \mathbf{e}_i,\mathbf{e}_{j} \rangle d\left[ \cos \left( \left( m'+m \right) \pi s \right) \right]\cr
	&&+\frac{m'}{\left( 1+m^2 \right) ^{\text{1/}2q}\left( 1+m^{'2} \right) ^{\text{1/}2q}\left( m'-m \right)}\cr
	&&\quad \times M_t\int_0^t{M}_{s}^{-1}[\nabla \sigma(\phi_{s}^{0})|\nabla\hat \sigma(\phi_{s}^{0})] \langle [\sigma(\phi_{s}^{0})|\hat \sigma(\phi_{s}^{0})] \mathbf{e}_i,\mathbf{e}_{j} \rangle d\left[ \cos \left( \left( m'-m \right) \pi s \right) \right]\cr
	&=:&I_1 +I_2.
	\end{eqnarray*}
	For the jog-free function, it is easy to check that $\left\|\cos (n \pi t)-1\right\|_{p  {\mathrm{-var} }}=2 n^{1 / p}$\cite{2013Inahama}. Next by the Young integral,  $2(1 / q-1 / p)>1$, $2 / q>1$, $4(1 / q-1 / 2)>1$, and $2(1-1 / p)>1$, we can obtain 
	\begin{eqnarray*}\label{4-20}
	\sum_{0 \leq m, m'<\infty, m \neq m'}\left\|I_1\right\|_{p  {\mathrm{-var} }}^{2}
	&\le& c \sum_{m \in \mathbf{Z}} \frac{1}{(1+|m|)^{2 / q}}\left(\sum_{m' \in \mathbf{Z}} \frac{1}{\left(1+\left|m'+m\right|\right)^{2(1 / q-1 / p)}}\right) \cr
	&+&c \sum_{m \in \mathbf{Z}} \frac{1}{(1+|m|)^{4(1 / q-1 / 2)}}\left(\sum_{m' \in \mathbf{Z}} \frac{1}{\left(1+\left|m'+m\right|\right)^{2(1-1 / p)}}\right)<\infty.
	\end{eqnarray*}
 Likewise, it has $\sum_{0 \leq m, m'<\infty, m \neq m'}\left\|I_2\right\|_{p  {\mathrm{-var} }}^{2}
 <\infty.$
	
	Next, we turn to the case that  $m=m'$,
	$$	\sqrt{2}\cos \left( m\pi t \right) d[ \sqrt{2}\cos \left( m\pi t \right) ] =d\left[ \cos \left( 2m\pi t \right) \right] /2,$$
	therefore
		\begin{eqnarray*}\label{4-21}
	R_1\big< (f,\hat f)_{m,i},(f,\hat f)_{m,j} \big> _t
	&=&M_t\int_0^t{M}_{s}^{-1}[\nabla \sigma(\phi_{s}^{0})|\nabla\hat \sigma(\phi_{s}^{0})] \langle [\sigma(\phi_{s}^{0})|\hat \sigma(\phi_{s}^{0})] \mathbf{e}_i,\mathbf{e}_{j} \rangle \sqrt{2}\frac{\cos \left( m\pi t \right)}{\left( 1+m^2 \right) ^{\text{1/}2q}}\cr
    &&\quad\times d\left[ \frac{\sqrt{2}\cos \left( m\pi s \right)}{\left( 1+m^2 \right) ^{\text{1/}2q}} \right]\cr 
	&=&\frac{\text{1/}2}{\left( 1+m^2 \right) ^{\text{1/}q}}\cr
    &&\quad\times M_t\int_0^t{M}_{s}^{-1}[\nabla \sigma(\phi_{s}^{0})|\nabla\hat \sigma(\phi_{s}^{0})] \langle [\sigma(\phi_{s}^{0})|\hat \sigma(\phi_{s}^{0})] \mathbf{e}_i,\mathbf{e}_{j} \rangle d\left[ \cos \left( 2m\pi s \right) \right] .
	\end{eqnarray*}
	Since $4 / q-2 / p>1$, 
	\begin{eqnarray*}\label{4-22}
	\sum_{i,j=1}^{d_1}\sum_{m=0}^{\infty}\big\|R_{1}\big\langle (f,\hat f)_{m, i}, (f,\hat f)_{m, j}\big\rangle\big\|_{p  {\mathrm{-var} }}^{2}  &\leq&\frac{c}{\left(1+m^{2}\right)^{2 / q}}\|\cos (2 m \pi \cdot)-1\|_{p\mathrm{-var}}^{2} \cr
	&\leq &\sum_{m=0}^{\infty}\frac{c}{(1+m)^{4 / q-2 / p}}<\infty.
	\end{eqnarray*}

	\textbf{Step 2}.
	Fix $i=1,\cdots ,d_1, j=d_1+1,\cdots ,d_1+d_2$. 
	
	Consider that $(f,\hat f)^\mathrm{T}_{m,i}= \sqrt{2}\frac{\cos \left( m\pi t \right)}{\left( 1+m^2 \right) ^{\text{1/}2q}},  (f,\hat f)^\mathrm{T}_{2m'-1,j}= \frac{\cos \left( 2m'\pi t \right)-1}{2m'\pi}$, and $m \ne 2m'$, then
	\begin{eqnarray*}\label{4-23}
	\sqrt{2}\cos \left( m\pi t \right) d\left[ \frac{\cos \left( 2m'\pi t \right)-1}{2m'\pi} \right] 
	&=& \frac{1}{\sqrt{2}}d\left[ \frac{\cos ( ( 2m'+m ) \pi t )}{(2m'+m)\pi}-\frac{\cos ( ( 2m'-m ) \pi t )}{(2m'-m)\pi} \right],
	\end{eqnarray*}
	then
	\begin{eqnarray*}\label{4-24}
	R_1\big< (f,\hat f)_{m,i},  (f,\hat f)_{2m'-1,j} \big> _t
	&=&M_t\int_0^t{M}_{s}^{-1}[\nabla \sigma(\phi_{s}^{0})|\nabla\hat \sigma(\phi_{s}^{0})] \langle [\sigma(\phi_{s}^{0})|\hat \sigma(\phi_{s}^{0})] \mathbf{e}_i,\mathbf{e}_{j} \rangle \sqrt{2}\frac{\cos \left( m\pi t \right)}{\left( 1+m^2 \right) ^{\text{1/}2q}}\cr
   &&\quad\times d\left[ \frac{\cos \left( 2m'\pi s \right)-1}{2m'\pi} \right]\cr
	&=&\frac{1}{\sqrt{2}\left( 1+m^2 \right) ^{\text{1/}2q}\left( 2m'+m \right)\pi}\cr
	&&\quad\times M_t\int_0^t{M}_{s}^{-1}[\nabla \sigma(\phi_{s}^{0})|\nabla\hat \sigma(\phi_{s}^{0})] \langle [\sigma(\phi_{s}^{0})|\hat \sigma(\phi_{s}^{0})] \mathbf{e}_i,\mathbf{e}_{j} \rangle d\left[ \cos \left( \left( 2m'+m \right) \pi s \right) \right]\cr
	&&+\frac{1}{\sqrt{2}\left( 1+m^2 \right) ^{\text{1/}2q}\left( 2m'-m \right)\pi}\cr
	&&\quad \times M_t\int_0^t{M}_{s}^{-1}[\nabla \sigma(\phi_{s}^{0})|\nabla\hat \sigma(\phi_{s}^{0})] \langle [\sigma(\phi_{s}^{0})|\hat \sigma(\phi_{s}^{0})] \mathbf{e}_i,\mathbf{e}_{j} \rangle d\left[ \cos \left( \left( 2m'-m \right) \pi s \right) \right]\cr
	&=:&I_3+I_4.
	\end{eqnarray*}
    Next by the Young integral, and  $2 / q>1$,  we can obtain 
		\begin{eqnarray*}\label{4-25}
		\sum_{0 \leq m, m'<\infty, m \neq 2m'}\left\|I_3\right\|_{p  {\mathrm{-var} }}^{2}
		&\le& c \sum_{m \in \mathbf{Z}} \frac{1}{(1+|m|)^{2 / q}}\left(\sum_{m' \in \mathbf{Z}} \frac{1}{\left(1+\left|2m'+m\right|\right)^{2 / q}}\right)<\infty.
		\end{eqnarray*}
		Likewise, it has $\sum_{0 \leq m, m'<\infty, m \neq 2m'}\left\|I_4\right\|_{p  {\mathrm{-var} }}^{2}
		<\infty.$
	
	When $(f,\hat f)^\mathrm{T}_{m,i}= \sqrt{2}\frac{\cos \left( m\pi t \right)}{\left( 1+m^2 \right) ^{\text{1/}2q}},  (f,\hat f)^\mathrm{T}_{2m'-1,j}= \frac{\cos \left( 2m'\pi t \right)-1}{2m'\pi}$, and $m =2 m'$, then
	\begin{eqnarray*}\label{4-26}
	\sqrt{2}\cos \left( m\pi t \right) d\left[ \frac{\cos \left( m\pi t \right)-1}{m\pi} \right] 
	&=& \frac{1}{\sqrt{2}}d\left[ \frac{\cos ( 2m \pi t )}{2m\pi} \right],
	\end{eqnarray*}
	then
	\begin{eqnarray*}\label{4-27}
	\left. R_1\big< (f,\hat f)_{m,i}, \right. (f,\hat f)_{2m'-1,j} \big> _t
	&=&M_t\int_0^t{M}_{s}^{-1}[\nabla \sigma(\phi_{s}^{0})|\nabla\hat \sigma(\phi_{s}^{0})] \langle [\sigma(\phi_{s}^{0})|\hat \sigma(\phi_{s}^{0})] \mathbf{e}_i,\mathbf{e}_{j} \rangle \sqrt{2}\frac{\cos \left( m\pi t \right)}{\left( 1+m^2 \right) ^{\text{1/}2q}}\cr
	&&\quad\times d\left[ \frac{\cos \left( 2m'\pi s \right)-1}{2m'\pi} \right]\cr
	&=&\frac{1}{2\sqrt{2}\left( 1+m^2 \right) ^{\text{1/}2q} m\pi}\cr
	&&\quad\times M_t\int_0^t{M}_{s}^{-1}[\nabla \sigma(\phi_{s}^{0})|\nabla\hat \sigma(\phi_{s}^{0})] \langle [\sigma(\phi_{s}^{0})|\hat \sigma(\phi_{s}^{0})] \mathbf{e}_i,\mathbf{e}_{j} \rangle d\left[ \cos \left( 2m\pi s \right) \right]\cr
		&=:&I_5.
	\end{eqnarray*}
	
    In the same manner,  we can deduce that 
		 $\sum_{0 \leq m<\infty}\left\|I_5\right\|_{p  {\mathrm{-var} }}^{2}
		<\infty.$

	Consider that $(f,\hat f)^\mathrm{T}_{m,i}= \sqrt{2}\frac{\cos \left( m\pi t \right)}{\left( 1+m^2 \right) ^{\text{1/}2q}},  (f,\hat f)^\mathrm{T}_{2m',j}=  \frac{\sin \left( 2m'\pi t \right)}{2m'\pi}$, and $m \ne 2m'$, then
	\begin{eqnarray*}\label{4-28}
	\sqrt{2}\cos \left( m\pi t \right) d\left[ \frac{\sin \left( 2m'\pi t \right)}{2m'\pi} \right] 
	&=& \frac{1}{\sqrt{2}}d\left[ \frac{\sin ( ( 2m'+m ) \pi t )}{(2m'+m)\pi}\right]+\frac{1}{\sqrt{2}}d\left[\frac{\sin ( ( 2m'-m ) \pi t )}{(2m'-m)\pi} \right],
	\end{eqnarray*}
it follows that
	\begin{eqnarray*}\label{4-29}
	\left. R_1\big< (f,\hat f)_{m,i}, \right. (f,\hat f)_{2m',j} \big> _t
	&=&M_t\int_0^t{M}_{s}^{-1}[\nabla \sigma(\phi_{s}^{0})|\nabla\hat \sigma(\phi_{s}^{0})] \langle [\sigma(\phi_{s}^{0})|\hat \sigma(\phi_{s}^{0})] \mathbf{e}_i,\mathbf{e}_{j} \rangle \sqrt{2}\frac{\cos \left( m\pi t \right)}{\left( 1+m^2 \right) ^{\text{1/}2q}}d\left[ \frac{\sin \left( 2m'\pi s \right)}{2m'\pi} \right]\cr
	&=&\frac{1}{\sqrt{2}\left( 1+m^2 \right) ^{\text{1/}2q}\left( 2m'+m \right)\pi}\cr
	&&\quad\times M_t\int_0^t{M}_{s}^{-1}[\nabla \sigma(\phi_{s}^{0})|\nabla\hat \sigma(\phi_{s}^{0})] \langle [\sigma(\phi_{s}^{0})|\hat \sigma(\phi_{s}^{0})] \mathbf{e}_i,\mathbf{e}_{j} \rangle d\left[ \sin \left( \left( 2m'+m \right) \pi s \right) \right]\cr
	&&+\frac{1}{\sqrt{2}\left( 1+m^2 \right) ^{\text{1/}2q}\left( 2m'-m \right)\pi}\cr
	&&\quad \times M_t\int_0^t{M}_{s}^{-1}[\nabla \sigma(\phi_{s}^{0})|\nabla\hat \sigma(\phi_{s}^{0})] \langle [\sigma(\phi_{s}^{0})|\hat \sigma(\phi_{s}^{0})] \mathbf{e}_i,\mathbf{e}_{j} \rangle d\left[ \sin \left( \left( 2m'-m \right) \pi s \right) \right]\cr
	&=:&I_6 +I_7.
	\end{eqnarray*}
	
By the Young integral and $2 / q>1$, we can obtain 
		\begin{eqnarray}\label{4-30}
		\sum_{0 \leq m, m'<\infty, m \neq 2m'}\left\|I_6\right\|_{p  {\mathrm{-var} }}^{2}
		&\le& c \sum_{m \in \mathbf{Z}} \frac{1}{(1+|m|)^{2 / q}}\left(\sum_{m' \in \mathbf{Z}} \frac{1}{\left(1+\left|2m'+m\right|\right)^{2 / q}}\right) <\infty.
		\end{eqnarray}
		Likewise, it has $\sum_{0 \leq m, m'<\infty, m \neq 2m'}\left\|I_6\right\|_{p  {\mathrm{-var} }}^{2}
		<\infty.$

	Consider that $(f,\hat f)^\mathrm{T}_{m,i}= \sqrt{2}\frac{\cos \left( m\pi t \right)}{\left( 1+m^2 \right) ^{\text{1/}2q}},  (f,\hat f)^\mathrm{T}_{2m',j}=  \frac{\sin \left( 2m'\pi t \right)}{2m'\pi}$, and $m = 2m'$, then
	\begin{eqnarray*}\label{4-31}
	\sqrt{2}\cos \left( m\pi t \right) d\left[ \frac{\sin \left( m\pi t \right)}{m\pi} \right] 
	&=& \frac{1}{\sqrt{2}}d\left[ \frac{\sin ( 2m \pi t )}{2m\pi} \right]+\frac{1}{\sqrt{2}}dt.
	\end{eqnarray*}
	Therefore,
	\begin{eqnarray*}\label{4-32}
	R_1\big< (f,\hat f)_{m,i},  (f,\hat f)_{m',j} \big> _t
	&=&M_t\int_0^t{M}_{s}^{-1}[\nabla \sigma(\phi_{s}^{0})|\nabla\hat \sigma(\phi_{s}^{0})] \langle [\sigma(\phi_{s}^{0})|\hat \sigma(\phi_{s}^{0})] \mathbf{e}_i,\mathbf{e}_{j} \rangle \sqrt{2}\frac{\cos \left( m\pi t \right)}{\left( 1+m^2 \right) ^{\text{1/}2q}}d\left[\frac{\sin \left( m\pi s \right)}{m\pi}  \right]\cr
	&=&\frac{1}{2\sqrt{2}\left( 1+m^2 \right) ^{\text{1/}2q} m\pi}\cr
	&&\quad\times M_t\int_0^t{M}_{s}^{-1}[\nabla \sigma(\phi_{s}^{0})|\nabla\hat \sigma(\phi_{s}^{0})] \langle [\sigma(\phi_{s}^{0})|\hat \sigma(\phi_{s}^{0})] \mathbf{e}_i,\mathbf{e}_{j} \rangle d\left[ \cos \left( 2m\pi s \right) \right]\cr
	&&+\frac{1}{\sqrt{2}\left( 1+m^2 \right) ^{\text{1/}2q} }\times M_t\int_0^t{M}_{s}^{-1}[\nabla \sigma(\phi_{s}^{0})|\nabla\hat \sigma(\phi_{s}^{0})] \langle [\sigma(\phi_{s}^{0})|\hat \sigma(\phi_{s}^{0})] \mathbf{e}_i,\mathbf{e}_{j} \rangle dt\cr
		&=:&I_8 +I_9.
	\end{eqnarray*}
	
		By that   $2(1 / q-1 / p)>1$,  we can obtain 
		\begin{eqnarray*}\label{4-33}
		\sum_{0 \leq m<\infty}\left\|I_8\right\|_{p  {\mathrm{-var} }}^{2}
		&\le& c \sum_{m \in \mathbf{Z}} \frac{1}{m^2(1+|m|)^{2(1/ q-1/p)}} <\infty.
		\end{eqnarray*}
		Likewise, it has $\sum_{0 \leq m<\infty}\left\|I_9\right\|_{p  {\mathrm{-var} }}^{2}
		<\infty.$

	\textbf{Step 3}.
	Fix $i=d_1+1,\cdots ,d_1+d_2, j=1,\cdots ,d_1$.
	
	Consider that $(f,\hat f)^\mathrm{T}_{2m-1,i}= \frac{\cos \left( 2m\pi t \right)-1}{2m\pi},  (f,\hat f)^\mathrm{T}_{m',j}= \sqrt{2}\frac{\cos \left( m'\pi t \right)}{\left( 1+m'^2 \right) ^{\text{1/}2q}} $, and $2m \ne m'$, then
	\begin{eqnarray*}\label{4-34}
	\frac{\cos \left( 2m\pi t \right)-1}{2m\pi} d\left[ \sqrt{2}\cos \left( m'\pi t \right) \right] 
	&=& \frac{m'}{2\sqrt{2}m}d\left[ \frac{\cos ( ( 2m +m') \pi t )}{(2m +m')\pi}-\frac{\cos ( ( 2m -m') \pi t )}{(2m -m')\pi} \right]\cr
	&&-\frac{1}{2m\pi}d\left[ \sqrt{2}\cos \left( m'\pi t \right) \right] ,
	\end{eqnarray*}
	hence
	\begin{eqnarray*}\label{4-35}
	R_1\big< (f,\hat f)_{2m,i},  (f,\hat f)_{m',j} \big> _t
	&=&M_t\int_0^t{M}_{s}^{-1}[\nabla \sigma(\phi_{s}^{0})|\nabla\hat \sigma(\phi_{s}^{0})] \langle [\sigma(\phi_{s}^{0})|\hat \sigma(\phi_{s}^{0})] \mathbf{e}_i,\mathbf{e}_{j} \rangle \sqrt{2}\frac{\cos \left( m\pi t \right)}{\left( 1+m^2 \right) ^{\text{1/}2q}}\cr
	&&\quad \times d\left[ \frac{\cos \left( 2m'\pi s \right)-1}{2m'\pi} \right]\cr
	&=&\frac{m'}{2\sqrt{2}\left( 1+m'^2 \right) ^{\text{1/}2q}(2m +m')\pi}\cr
	&&\quad\times M_t\int_0^t{M}_{s}^{-1}[\nabla \sigma(\phi_{s}^{0})|\nabla\hat \sigma(\phi_{s}^{0})] \langle [\sigma(\phi_{s}^{0})|\hat \sigma(\phi_{s}^{0})] \mathbf{e}_i,\mathbf{e}_{j} \rangle d\left[ \cos \left( (2m +m') \pi s \right) \right]\cr
	&&-\frac{m'}{\sqrt{2}\left( 1+m'^2 \right) ^{\text{1/}2q}(2m -m')\pi}\cr
	&&\quad \times M_t\int_0^t{M}_{s}^{-1}[\nabla \sigma(\phi_{s}^{0})|\nabla\hat \sigma(\phi_{s}^{0})] \langle [\sigma(\phi_{s}^{0})|\hat \sigma(\phi_{s}^{0})] \mathbf{e}_i,\mathbf{e}_{j} \rangle d\left[ \cos \left((2m -m') \pi s \right) \right]\cr
	&&-\frac{1}{\sqrt{2}\left( 1+m'^2 \right) ^{\text{1/}2q}m\pi}\cr
	&&\quad \times M_t\int_0^t{M}_{s}^{-1}[\nabla \sigma(\phi_{s}^{0})|\nabla\hat \sigma(\phi_{s}^{0})] \langle [\sigma(\phi_{s}^{0})|\hat \sigma(\phi_{s}^{0})] \mathbf{e}_i,\mathbf{e}_{j} \rangle d\left[ \cos \left( m' \pi s \right) \right]\cr
	&=:&I_{10}+I_{11}+I_{12}.
	\end{eqnarray*}
	
Next by the Young integral,   $2 / q>1$ and $2(1-1 / p)>1$, we have
		\begin{eqnarray*}\label{4-36}
		\sum_{0 \leq m, m'<\infty, 2m \neq m'}\left\|I_{10}\right\|_{p  {\mathrm{-var} }}^{2}
		&\le& c \sum_{0 \leq m, m'<\infty, 2m \neq m'}\frac{m'(1+|m'+2m|)^{2/p}}{(1+|m'|)^{2/q}\left|2m'+m\right|^2} \cr
		&\le& c \sum_{m \in \mathbf{Z}} \frac{1}{(1+|m|)^2}\left(\sum_{m' \in \mathbf{Z}} \frac{1}{\left(1+\left|2m'+m\right|\right)^{2(1 / q-1 / p)}}\right) \cr
		&&+c \sum_{m \in \mathbf{Z}} \frac{1}{(1+|m|)^{2 / q}}\left(\sum_{m' \in \mathbf{Z}} \frac{1}{\left(1+\left|2m'+m\right|\right)^{2(1-1 / p)}}\right)\cr
		&<&\infty.
		\end{eqnarray*}
		Likewise, it has $\sum_{0 \leq m, m'<\infty, 2m \neq m'}\left\|I_{11}\right\|_{p  {\mathrm{-var} }}^{2}
		<\infty$, and $\sum_{0 \leq m, m'<\infty, 2m \neq m'}\left\|I_{12}\right\|_{p  {\mathrm{-var} }}^{2}
		<\infty$.

	Consider that $(f,\hat f)^\mathrm{T}_{2m-1,i}= \frac{\cos \left( 2m\pi t \right)-1}{2m\pi},  (f,\hat f)^\mathrm{T}_{m',j}= \sqrt{2}\frac{\cos \left( m'\pi t \right)}{\left( 1+m'^2 \right) ^{\text{1/}2q}} $, and $2m = m'$, then
	\begin{eqnarray*}\label{4-37}
	\frac{\cos \left( m\pi t \right)-1}{2m\pi} d\left[ \sqrt{2}\cos \left( m\pi t \right) \right] 
	&=& \frac{1}{2\sqrt{2}}d\left[ \frac{\cos ( 2m \pi t )}{2m\pi}\right]-\frac{1}{2m\pi}d\left[ \sqrt{2}\cos \left( 2m\pi t \right) \right] ,
	\end{eqnarray*}
	then
	\begin{eqnarray*}\label{4-38}
    R_1\big< (f,\hat f)_{2m-1,i}, (f,\hat f)_{m',j} \big> _t
	&=&M_t\int_0^t{M}_{s}^{-1}[\nabla \sigma(\phi_{s}^{0})|\nabla\hat \sigma(\phi_{s}^{0})] \langle [\sigma(\phi_{s}^{0})|\hat \sigma(\phi_{s}^{0})] \mathbf{e}_i,\mathbf{e}_{j} \rangle \frac{\cos \left( m\pi t \right)-1}{m\pi}d\left[ \frac{ \sqrt{2}\cos \left( m\pi s \right)}{\left( 1+m^2 \right) ^{\text{1/}2q}} \right]\cr
	&=&\frac{1}{4\sqrt{2}\left( 1+m^2 \right) ^{\text{1/}2q} m\pi}\cr
	&&\quad\times M_t\int_0^t{M}_{s}^{-1}[\nabla \sigma(\phi_{s}^{0})|\nabla\hat \sigma(\phi_{s}^{0})] \langle [\sigma(\phi_{s}^{0})|\hat \sigma(\phi_{s}^{0})] \mathbf{e}_i,\mathbf{e}_{j} \rangle d\left[ \cos \left( 2m\pi s \right) \right]\cr
	&&-\frac{\sqrt{2}}{2\left( 1+m^2 \right) ^{\text{1/}2q} m\pi}\cr
	&&\quad\times M_t\int_0^t{M}_{s}^{-1}[\nabla \sigma(\phi_{s}^{0})|\nabla\hat \sigma(\phi_{s}^{0})] \langle [\sigma(\phi_{s}^{0})|\hat \sigma(\phi_{s}^{0})] \mathbf{e}_i,\mathbf{e}_{j} \rangle d\left[ \cos \left( 2m\pi s \right) \right]\cr
	&=:&I_{13}+I_{14}.
	\end{eqnarray*}
	
	Next by the Young integral and   $2(1 / q-1 / p)>1$,  we can obtain 
		\begin{eqnarray*}\label{4-39}
		\sum_{0 \leq m<\infty}\left\|I_{13}\right\|_{p  {\mathrm{-var} }}^{2}
		&\le& c \sum_{m \in \mathbf{Z}} \frac{1}{m^2{(1+m)^{2(1 / q-1 / p)}} } <\infty.
		\end{eqnarray*}
		Likewise, it has $\sum_{0 \leq m<\infty}\left\|I_{14}\right\|_{p  {\mathrm{-var} }}^{2}
		<\infty.$

	Consider that $(f,\hat f)^\mathrm{T}_{2m,i}= \frac{\sin \left( 2m\pi t \right)}{2m\pi},  (f,\hat f)^\mathrm{T}_{m',j}= \sqrt{2}\frac{\cos \left( m'\pi t \right)}{\left( 1+m'^2 \right) ^{\text{1/}2q}} $, and $2m \ne m'$, then
	\begin{eqnarray*}\label{4-40}
\frac{\sin \left( 2m\pi t \right)}{2m\pi} d\left[ \sqrt{2} \cos \left( m'\pi t \right)\right] 
	&=& \frac{\sqrt{2}m'}{4m}d\left[ \frac{\sin ( ( 2m+m' ) \pi t )}{(2m+m')\pi}-\frac{\sin ( ( 2m-m' ) \pi t )}{(2m-m')\pi} \right],
	\end{eqnarray*}
	then it has
	\begin{eqnarray*}\label{4-41}
	 R_1\big< (f,\hat f)_{2m,i}, (f,\hat f)_{m',j} \big> _t
	&=&M_t\int_0^t{M}_{s}^{-1}[\nabla \sigma(\phi_{s}^{0})|\nabla\hat \sigma(\phi_{s}^{0})] \langle [\sigma(\phi_{s}^{0})|\hat \sigma(\phi_{s}^{0})] \mathbf{e}_i,\mathbf{e}_{j} \rangle \frac{\sin \left( 2m\pi t \right)}{2m\pi} d\left[ \sqrt{2} \frac{\cos \left( m'\pi s \right)}{\left( 1+m'^2 \right) ^{\text{1/}2q}}\right] \cr
	&=&\frac{\sqrt{2}m'}{4\left( 1+m'^2 \right) ^{\text{1/}2q}\left( 2m+m' \right)\pi}\cr
	&&\quad\times M_t\int_0^t{M}_{s}^{-1}[\nabla \sigma(\phi_{s}^{0})|\nabla\hat \sigma(\phi_{s}^{0})] \langle [\sigma(\phi_{s}^{0})|\hat \sigma(\phi_{s}^{0})] \mathbf{e}_i,\mathbf{e}_{j} \rangle d\left[ \sin \left( \left( 2m+m' \right) \pi s \right) \right]\cr
	&&-\frac{\sqrt{2}m'}{4\left( 1+m'^2 \right) ^{\text{1/}2q}\left( 2m-m' \right)\pi}\cr
	&&\quad \times M_t\int_0^t{M}_{s}^{-1}[\nabla \sigma(\phi_{s}^{0})|\nabla\hat \sigma(\phi_{s}^{0})] \langle [\sigma(\phi_{s}^{0})|\hat \sigma(\phi_{s}^{0})] \mathbf{e}_i,\mathbf{e}_{j} \rangle d\left[ \sin \left( \left( 2m-m' \right) \pi s \right) \right]\cr
	&=:&I_{15}+I_{16}.
	\end{eqnarray*}
	
	Similar to (\ref{4-30}),   we have
	$\sum_{0 \leq m, m'<\infty, m \neq m'}\left\|I_{15}\right\|_{p  {\mathrm{-var} }}^{2}
	<\infty$, and $\sum_{0 \leq m, m'<\infty, m \neq m'}\left\|I_{16}\right\|_{p  {\mathrm{-var} }}^{2}
	<\infty.$

	Consider that $(f,\hat f)^\mathrm{T}_{2m,i}= \frac{\sin \left( 2m\pi t \right)}{2m\pi},  (f,\hat f)^\mathrm{T}_{m',j}= \sqrt{2}\frac{\cos \left( m'\pi t \right)}{\left( 1+m'^2 \right) ^{\text{1/}2q}} $, and $2m = m'$, then
	\begin{eqnarray*}\label{4-42}
	\frac{\sin \left( 2m\pi t \right)}{2m\pi} d\left[ \sqrt{2} \cos \left( 2m\pi t \right)\right] 
	&=& \frac{\sqrt{2}}{8m\pi}d\left[ \sin (  4m  \pi t ) \right]-\frac{1}{\sqrt{2}}dt,
	\end{eqnarray*}
	Therefore,
	\begin{eqnarray}\label{4-43}
	 R_1\big< (f,\hat f)_{2m,i},  (f,\hat f)_{m',j} \big> _t
	&=&M_t\int_0^t{M}_{s}^{-1}[\nabla \sigma(\phi_{s}^{0})|\nabla\hat \sigma(\phi_{s}^{0})] \langle [\sigma(\phi_{s}^{0})|\hat \sigma(\phi_{s}^{0})] \mathbf{e}_i,\mathbf{e}_{j} \rangle \frac{\sin \left( 2m\pi t \right)}{2m\pi} \cr
	&&\quad \times d\left[ \sqrt{2} \frac{\cos \left( 2m\pi s \right)}{\left( 1+(2m)^2 \right) ^{\text{1/}2q}}\right] \cr
	&=&\frac{1}{4\sqrt{2}\left( 1+(2m)^2 \right) ^{\text{1/}2q} m\pi} \cr
	&&\quad \times M_t\int_0^t{M}_{s}^{-1}[\nabla \sigma(\phi_{s}^{0})|\nabla\hat \sigma(\phi_{s}^{0})] \langle [\sigma(\phi_{s}^{0})|\hat \sigma(\phi_{s}^{0})] \mathbf{e}_i,\mathbf{e}_{j} \rangle d\left[ \sin \left( 4m\pi s \right) \right]\cr
	&&-\frac{1}{\sqrt{2}\left( 1+(2m)^2 \right) ^{\text{1/}2q} } \cr
	&&\quad \times M_t\int_0^t{M}_{s}^{-1}[\nabla \sigma(\phi_{s}^{0})|\nabla\hat \sigma(\phi_{s}^{0})] \langle [\sigma(\phi_{s}^{0})|\hat \sigma(\phi_{s}^{0})] \mathbf{e}_i,\mathbf{e}_{j} \rangle dt\cr
		&=:&I_{17}+I_{18}.
	\end{eqnarray}
	
	 Analysis similar to the proof in (\ref{4-30}),  we claim that
	 $\sum_{0 \leq m<\infty}\left\|I_{17}\right\|_{p  {\mathrm{-var} }}^{2}
		<\infty$, and $\sum_{0 \leq m<\infty}\left\|I_{18}\right\|_{p  {\mathrm{-var} }}^{2}
		<\infty$.

	\textbf{Step 4}.	
	Fix $i=d_1+1,\cdots ,d_1+d_2, j=d_1+1,\cdots ,d_1+d_2$.

    When $(f,\hat f)^\mathrm{T}_{2m-1,i}= \frac{\cos \left( 2m\pi t \right)-1}{2m\pi},  (f,\hat f)^\mathrm{T}_{2m'-1,j}= \frac{\cos \left( 2m'\pi t \right)-1}{2m'\pi}$ and $m \ne m'$,
	\begin{eqnarray*}\label{4-47}
	\left[ \frac{\cos \left( 2m\pi t \right)-1}{2m\pi} \right] d\left[ \frac{\cos \left( 2m'\pi t \right)-1}{2m'\pi} \right] 
	&=& \frac{1}{4m\pi}d\left[ \frac{\cos ( (2m+2m')\pi t )}{(2m+2m')\pi}\right]-\frac{1}{4m\pi}d\left[ \frac{\cos ( (2m-2m')\pi t )}{(2m-2m')\pi}\right]\cr
    &&-\frac{1}{4mm'\pi^2}d\left[ \cos ( 2m'\pi t )\right].
	\end{eqnarray*}
	then
	\begin{eqnarray*}\label{4-48}
	 R_1\big< (f,\hat f)_{2m-1,i}, (f,\hat f)_{2m'-1,j} \big> _t
	&=&M_t\int_0^t{M}_{s}^{-1}[\nabla \sigma(\phi_{s}^{0})|\nabla\hat \sigma(\phi_{s}^{0})] \langle [\sigma(\phi_{s}^{0})|\hat \sigma(\phi_{s}^{0})] \mathbf{e}_i,\mathbf{e}_{j} \rangle 	\left[ \frac{\cos \left( 2m\pi t \right)-1}{2m\pi} \right] \cr
	&&\quad \times d\left[ \frac{\cos \left( 2m'\pi s \right)-1}{2m'\pi} \right] \cr
	&=&\frac{1}{4m(2m+2m')\pi^2}\cr
	&&\quad\times M_t\int_0^t{M}_{s}^{-1}[\nabla \sigma(\phi_{s}^{0})|\nabla\hat \sigma(\phi_{s}^{0})] \langle [\sigma(\phi_{s}^{0})|\hat \sigma(\phi_{s}^{0})] \mathbf{e}_i,\mathbf{e}_{j} \rangle d\left[\cos ( (2m+2m')\pi s )\right]\cr
	&&-\frac{1}{4m(2m-2m')\pi^2}\cr
	&&\quad\times M_t\int_0^t{M}_{s}^{-1}[\nabla \sigma(\phi_{s}^{0})|\nabla\hat \sigma(\phi_{s}^{0})] \langle [\sigma(\phi_{s}^{0})|\hat \sigma(\phi_{s}^{0})] \mathbf{e}_i,\mathbf{e}_{j} \rangle d\left[ \cos ( (2m-2m')\pi s )\right]\cr
	&&-\frac{1}{4m'm\pi^2} M_t\int_0^t{M}_{s}^{-1}[\nabla \sigma(\phi_{s}^{0})|\nabla\hat \sigma(\phi_{s}^{0})] \langle [\sigma(\phi_{s}^{0})|\hat \sigma(\phi_{s}^{0})] \mathbf{e}_i,\mathbf{e}_{j} \rangle d\left[ \cos ( 2m'\pi s )\right]\cr
			&=:&I_{19}+I_{20}+I_{21}.
	\end{eqnarray*}
	
	The Young integral and $2 / q>1$ yield that 
		\begin{eqnarray}\label{4-49}
		\sum_{0 \leq m, m'<\infty, m \neq m'}\left\|I_{19}\right\|_{p  {\mathrm{-var} }}^{2}
		&\le& c \sum_{m \in \mathbf{Z}} \frac{1}{(1+|m|)^2}\left(\sum_{m' \in \mathbf{Z}} \frac{1}{\left(1+\left|m'+m\right|\right)^{2 / q}}\right) <\infty.
		\end{eqnarray}
		Likewise, it has $\sum_{0 \leq m, m'<\infty, m \neq m'}\left\|I_{20}\right\|_{p  {\mathrm{-var} }}^{2}
		<\infty$, $\sum_{0 \leq m, m'<\infty, m \neq m'}\left\|I_{21}\right\|_{p  {\mathrm{-var} }}^{2}
		<\infty.$
	
	When $(f,\hat f)^\mathrm{T}_{2m-1,i}= \frac{\cos \left( 2m\pi t \right)-1}{2m\pi},  (f,\hat f)^\mathrm{T}_{2m'-1,j}= \frac{\cos \left( 2m'\pi t \right)-1}{2m'\pi}$ and $m=m'$,
	\begin{eqnarray*}\label{4-44}
		\left[ \frac{\cos \left( 2m\pi t \right)-1}{2m\pi} \right] d\left[ \frac{\cos \left( 2m\pi t \right)-1}{2m\pi} \right] 
		&=& d\left[ \frac{\cos ( 4m\pi t )}{16m^2\pi^2}\right]-d\left[\frac{\cos ( 2m\pi t )-1}{4m^2\pi^2} \right],
	\end{eqnarray*}
	then
		\begin{eqnarray*}\label{4-45}
		R_1\big< (f,\hat f)_{2m-1,i},  (f,\hat f)_{2m-1,j} \big> _t
		&=&M_t\int_0^t{M}_{s}^{-1}[\nabla \sigma(\phi_{s}^{0})|\nabla\hat \sigma(\phi_{s}^{0})] \langle [\sigma(\phi_{s}^{0})|\hat \sigma(\phi_{s}^{0})] \mathbf{e}_i,\mathbf{e}_{j} \rangle \left[ \frac{\cos \left( 2m\pi t \right)-1}{2m\pi} \right] \cr
		&&\quad \times d\left[ \frac{\cos \left( 2m\pi s \right)-1}{2m\pi} \right]\cr
		&=&\frac{1}{16m^2\pi^2} M_t\int_0^t{M}_{s}^{-1}[\nabla \sigma(\phi_{s}^{0})|\nabla\hat \sigma(\phi_{s}^{0})] \langle [\sigma(\phi_{s}^{0})|\hat \sigma(\phi_{s}^{0})] \mathbf{e}_i,\mathbf{e}_{j} \rangle d\left[ \cos \left( 4m \pi s \right) \right]\cr
		&&-\frac{1}{4m^2\pi^2} M_t\int_0^t{M}_{s}^{-1}[\nabla \sigma(\phi_{s}^{0})|\nabla\hat \sigma(\phi_{s}^{0})] \langle [\sigma(\phi_{s}^{0})|\hat \sigma(\phi_{s}^{0})] \mathbf{e}_i,\mathbf{e}_{j} \rangle d\left[ \cos \left( 2m \pi s \right) \right]\cr
		&=:&I_{22}+I_{23}.
	\end{eqnarray*}
	According to that $2 / q>1$,  we can prove 
	\begin{eqnarray*}\label{4-46}
		\sum_{0 \leq m<\infty}\left\|I_{22}\right\|_{p  {\mathrm{-var} }}^{2}
		&\le& c \sum_{m \in \mathbf{Z}} \frac{1}{m^{2 }(1+m)^{2/q}}<\infty.
	\end{eqnarray*}
	Likewise, it has $\sum_{0 \leq m<\infty}\left\|I_{23}\right\|_{p  {\mathrm{-var} }}^{2}
	<\infty.$

    When $(f,\hat f)^\mathrm{T}_{2m-1,i}= \frac{\cos \left( 2m\pi t \right)-1}{2m\pi},  (f,\hat f)^\mathrm{T}_{2m',j}= \frac{\sin \left( 2m'\pi t \right)}{2m'\pi}$ and $m\ne m'$,
	\begin{eqnarray*}\label{4-50}
	\left[ \frac{\cos \left( 2m\pi t \right)-1}{2m\pi} \right]d\left[ \frac{\sin \left( 2m'\pi t \right)}{2m'\pi} \right] 
	&=& \frac{1}{4m\pi}d\left[ \frac{\sin ( (2m+2m')\pi t )}{(2m+2m')\pi}\right]+\frac{1}{4m\pi}d\left[ \frac{\sin ( (2m-2m')\pi t )}{(2m-2m')\pi}\right]\cr
	&&-\frac{1}{4mm'\pi^2}d\left[ \sin ( 2m'\pi t )\right],
	\end{eqnarray*}
	then it has
	\begin{eqnarray*}\label{4-51}
	 R_1\big< (f,\hat f)_{2m-1,i},  (f,\hat f)_{2m',j} \big> _t
	&=&M_t\int_0^t{M}_{s}^{-1}[\nabla \sigma(\phi_{s}^{0})|\nabla\hat \sigma(\phi_{s}^{0})] \langle [\sigma(\phi_{s}^{0})|\hat \sigma(\phi_{s}^{0})] \mathbf{e}_i,\mathbf{e}_{j} \rangle \sqrt{2}\frac{\cos \left( m\pi s \right)}{\left( 1+m^2 \right) ^{\text{1/}2q}}d\left[ \frac{\sin \left( 2m'\pi s \right)}{2m'\pi} \right]\cr
	&=&\frac{1}{4m(2m+2m')\pi^2}\cr
	&&\quad\times M_t\int_0^t{M}_{s}^{-1}[\nabla \sigma(\phi_{s}^{0})|\nabla\hat \sigma(\phi_{s}^{0})] \langle [\sigma(\phi_{s}^{0})|\hat \sigma(\phi_{s}^{0})] \mathbf{e}_i,\mathbf{e}_{j} \rangle d\left[\sin ( (2m+2m')\pi s )\right]\cr
	&&+\frac{1}{4m(2m-2m')\pi^2}\cr
	&&\quad\times M_t\int_0^t{M}_{s}^{-1}[\nabla \sigma(\phi_{s}^{0})|\nabla\hat \sigma(\phi_{s}^{0})] \langle [\sigma(\phi_{s}^{0})|\hat \sigma(\phi_{s}^{0})] \mathbf{e}_i,\mathbf{e}_{j} \rangle d\left[ \sin ( (2m-2m')\pi s )\right]\cr
	&&-\frac{1}{4m'm\pi^2} M_t\int_0^t{M}_{s}^{-1}[\nabla \sigma(\phi_{s}^{0})|\nabla\hat \sigma(\phi_{s}^{0})] \langle [\sigma(\phi_{s}^{0})|\hat \sigma(\phi_{s}^{0})] \mathbf{e}_i,\mathbf{e}_{j} \rangle d\left[ \sin ( 2m'\pi s )\right]\cr
				&=:&I_{24}+I_{25}+I_{26}.
	\end{eqnarray*}
	
	Similar arguments apply to the (\ref{4-49}), we can prove 
		 $$\sum_{0 \leq m, m'<\infty, m \neq m'}\left\|I_{24}\right\|_{p  {\mathrm{-var} }}^{2}
		 <\infty, \quad \sum_{0 \leq m, m'<\infty, m \neq m'}\left\|I_{25}\right\|_{p  {\mathrm{-var} }}^{2}
		<\infty, \quad \sum_{0 \leq m, m'<\infty, m \neq m'}\left\|I_{26}\right\|_{p  {\mathrm{-var} }}^{2}
		<\infty.$$

    When $(f,\hat f)^\mathrm{T}_{2m-1,i}= \frac{\cos \left( 2m\pi t \right)-1}{2m\pi},  (f,\hat f)^\mathrm{T}_{2m',j}= \frac{\sin \left( 2m'\pi t \right)}{2m'\pi}$ and $m= m'$,
    \begin{eqnarray*}\label{4-52}
    \left[ \frac{\cos \left( 2m\pi t \right)-1}{2m\pi} \right]d\left[ \frac{\sin \left( 2m\pi t \right)}{2m\pi} \right] 
    &=& \frac{1}{16m^2\pi^2}d\left[\sin (4m\pi t )\right]+\frac{1}{4m\pi}dt-\frac{1}{4m^2\pi^2}d\left[ \sin ( 2m\pi t )\right],
    \end{eqnarray*}
     then it has
    \begin{eqnarray*}\label{4-53}
   R_1\big< (f,\hat f)_{2m-1,i}, (f,\hat f)_{2m',j} \big> _t
   &=&M_t\int_0^t{M}_{s}^{-1}[\nabla \sigma(\phi_{s}^{0})|\nabla\hat \sigma(\phi_{s}^{0})] \langle [\sigma(\phi_{s}^{0})|\hat \sigma(\phi_{s}^{0})] \mathbf{e}_i,\mathbf{e}_{j} \rangle     \left[ \frac{\cos \left( 2m\pi t \right)-1}{2m\pi} \right] \cr
   &&\quad \times d\left[ \frac{\sin \left( 2m\pi s \right)}{2m\pi} \right] \cr
    &=&\frac{1}{16m^2\pi^2}\times M_t\int_0^t{M}_{s}^{-1}[\nabla \sigma(\phi_{s}^{0})|\nabla\hat \sigma(\phi_{s}^{0})] \langle [\sigma(\phi_{s}^{0})|\hat \sigma(\phi_{s}^{0})] \mathbf{e}_i,\mathbf{e}_{j} \rangle d\left[\sin ( 4m\pi s )\right]\cr
   &&+\frac{1}{2m\pi}\times M_t\int_0^t{M}_{s}^{-1}[\nabla \sigma(\phi_{s}^{0})|\nabla\hat \sigma(\phi_{s}^{0})] \langle [\sigma(\phi_{s}^{0})|\hat \sigma(\phi_{s}^{0})] \mathbf{e}_i,\mathbf{e}_{j} \rangle ds\cr
   &&-\frac{1}{4m^2\pi^2} M_t\int_0^t{M}_{s}^{-1}[\nabla \sigma(\phi_{s}^{0})|\nabla\hat \sigma(\phi_{s}^{0})] \langle [\sigma(\phi_{s}^{0})|\hat \sigma(\phi_{s}^{0})] \mathbf{e}_i,\mathbf{e}_{j} \rangle d\left[ \sin ( 2m\pi s )\right]\cr
   &=:&I_{27}+I_{28}+I_{29}.
   \end{eqnarray*}
   Likewise, we can prove that 
   $\sum_{0 \leq m<\infty}\left\|I_{27}\right\|_{p  {\mathrm{-var} }}^{2}
   <\infty$, $\sum_{0 \leq m<\infty}\left\|I_{28}\right\|_{p  {\mathrm{-var} }}^{2}
   <\infty$ and $\sum_{0 \leq m<\infty}\left\|I_{29}\right\|_{p  {\mathrm{-var} }}^{2}
   <\infty$.
   
   When $(f,\hat f)^\mathrm{T}_{2m,i}= \frac{\sin \left( 2m\pi t \right)}{2m\pi},  (f,\hat f)^\mathrm{T}_{2m'-1,j}= \frac{\cos \left( 2m'\pi t \right)-1}{2m'\pi}$ and $m \ne m'$,
   \begin{eqnarray*}\label{4-54}
   \frac{\sin \left( 2m\pi t \right)}{2m\pi}d\left[ \frac{\cos \left( 2m'\pi t \right)-1}{2m'\pi} \right] 
   &=& \frac{1}{4m\pi}d\left[ \frac{\sin ( (2m+2m')\pi t )}{(2m+2m')\pi}\right]-\frac{1}{4m\pi}d\left[ \frac{\sin ( (2m-2m')\pi t )}{(2m-2m')\pi}\right].
   \end{eqnarray*}
   then
   \begin{eqnarray*}\label{4-55}
   R_1\big< (f,\hat f)_{2m,i},  (f,\hat f)_{2m'-1,j} \big> _t
   &=&M_t\int_0^t{M}_{s}^{-1}[\nabla \sigma(\phi_{s}^{0})|\nabla\hat \sigma(\phi_{s}^{0})] \langle [\sigma(\phi_{s}^{0})|\hat \sigma(\phi_{s}^{0})] \mathbf{e}_i,\mathbf{e}_{j} \rangle 	\left[ \frac{\cos \left( 2m\pi t \right)-1}{2m\pi} \right]  \cr
   &&\quad \times d\left[ \frac{\cos \left( 2m'\pi s \right)-1}{2m'\pi} \right] \cr
   &=&\frac{1}{4m(2m+2m')\pi^2}\cr
   &&\quad\times M_t\int_0^t{M}_{s}^{-1}[\nabla \sigma(\phi_{s}^{0})|\nabla\hat \sigma(\phi_{s}^{0})] \langle [\sigma(\phi_{s}^{0})|\hat \sigma(\phi_{s}^{0})] \mathbf{e}_i,\mathbf{e}_{j} \rangle d\left[\sin ( (2m+2m')\pi s )\right]\cr
   &&-\frac{1}{4m(2m-2m')\pi^2}\cr
   &&\quad\times M_t\int_0^t{M}_{s}^{-1}[\nabla \sigma(\phi_{s}^{0})|\nabla\hat \sigma(\phi_{s}^{0})] \langle [\sigma(\phi_{s}^{0})|\hat \sigma(\phi_{s}^{0})] \mathbf{e}_i,\mathbf{e}_{j} \rangle d\left[ \sin ( (2m-2m')\pi s )\right]\cr
   &=:&I_{30}+I_{31}.
   \end{eqnarray*}
   Similar arguments apply to the (\ref{4-30}), we can prove 
   $$\sum_{0 \leq m, m'<\infty, m \neq m'}\left\|I_{30}\right\|_{p  {\mathrm{-var} }}^{2}
   <\infty, \quad \sum_{0 \leq m, m'<\infty, m \neq m'}\left\|I_{31}\right\|_{p  {\mathrm{-var} }}^{2}
   <\infty.$$
   
   When $(f,\hat f)^\mathrm{T}_{2m,i}= \frac{\sin \left( 2m\pi t \right)}{2m\pi},  (f,\hat f)^\mathrm{T}_{2m'-1,j}= \frac{\cos \left( 2m'\pi t \right)-1}{2m'\pi}$ and $m = m'$,

	\begin{eqnarray*}\label{4-56}
	\left[ \frac{\sin \left( 2m\pi t \right)}{2m\pi} \right] d\left[ \frac{\cos \left( 2m\pi t \right)-1}{2m\pi} \right] 
	&=& \frac{1}{16m^2\pi^2}d\left[\sin ( 4m\pi t )\right]-\frac{1}{4m\pi}dt,
	\end{eqnarray*}
	then
	\begin{eqnarray*}\label{4-57}
	R_1\big< (f,\hat f)_{2m,i}, (f,\hat f)_{2m-1,j} \big> _t
	&=&M_t\int_0^t{M}_{s}^{-1}[\nabla \sigma(\phi_{s}^{0})|\nabla\hat \sigma(\phi_{s}^{0})] \langle [\sigma(\phi_{s}^{0})|\hat \sigma(\phi_{s}^{0})] \mathbf{e}_i,\mathbf{e}_{j} \rangle \frac{\sin \left( 2m\pi t \right)}{2m\pi} d\left[ \frac{\cos \left( 2m\pi s \right)-1}{2m\pi} \right]\cr
	&=&\frac{1}{16m^2\pi^2} M_t\int_0^t{M}_{s}^{-1}[\nabla \sigma(\phi_{s}^{0})|\nabla\hat \sigma(\phi_{s}^{0})] \langle [\sigma(\phi_{s}^{0})|\hat \sigma(\phi_{s}^{0})] \mathbf{e}_i,\mathbf{e}_{j} \rangle d\left[ \sin \left( 4m \pi s \right) \right]\cr
	&&-\frac{1}{4m\pi} M_t\int_0^t{M}_{s}^{-1}[\nabla \sigma(\phi_{s}^{0})|\nabla\hat \sigma(\phi_{s}^{0})] \langle [\sigma(\phi_{s}^{0})|\hat \sigma(\phi_{s}^{0})] \mathbf{e}_i,\mathbf{e}_{j} \rangle ds\cr
	&=:&      		I_{32}+I_{33}.
	\end{eqnarray*}
	 Similar to the (\ref{4-30}), we can show 
	$\sum_{0 \leq m<\infty}\left\|I_{32}\right\|_{p  {\mathrm{-var} }}^{2}
	<\infty$ and $\sum_{0 \leq m<\infty}\left\|I_{33}\right\|_{p  {\mathrm{-var} }}^{2}
	<\infty$.
	
	When $(f,\hat f)^\mathrm{T}_{2m,i}= \frac{\sin \left( 2m\pi t \right)}{2m\pi},  (f,\hat f)^\mathrm{T}_{2m',j}= \frac{\sin \left( 2m'\pi t \right)}{2m'\pi}$ and $m\ne m'$,
	\begin{eqnarray*}\label{4-58}
	 \frac{\sin \left( 2m\pi t \right)}{2m\pi}d\left[ \frac{\sin \left( 2m'\pi t \right)}{2m'\pi} \right] 
	&=&-\frac{1}{4m\pi}d\left[ \frac{\cos ( (2m+2m')\pi t )}{(2m+2m')\pi}\right]-\frac{1}{4m\pi}d\left[ \frac{\cos ( (2m-2m')\pi t )}{(2m-2m')\pi}\right],
	\end{eqnarray*}
	then it has
	\begin{eqnarray*}\label{4-59}
	 R_1\big< (f,\hat f)_{2m,i},  (f,\hat f)_{2m',j} \big> _t
	&=&M_t\int_0^t{M}_{s}^{-1}[\nabla \sigma(\phi_{s}^{0})|\nabla\hat \sigma(\phi_{s}^{0})] \langle [\sigma(\phi_{s}^{0})|\hat \sigma(\phi_{s}^{0})] \mathbf{e}_i,\mathbf{e}_{j} \rangle \sqrt{2}\frac{\cos \left( m\pi s \right)}{\left( 1+m^2 \right) ^{\text{1/}2q}}\cr
	&&\quad \times d\left[ \frac{\sin \left( 2m'\pi s \right)}{2m'\pi} \right]\cr
	&=&-\frac{1}{4m(2m+2m')\pi^2}\cr
	&&\quad\times M_t\int_0^t{M}_{s}^{-1}[\nabla \sigma(\phi_{s}^{0})|\nabla\hat \sigma(\phi_{s}^{0})] \langle [\sigma(\phi_{s}^{0})|\hat \sigma(\phi_{s}^{0})] \mathbf{e}_i,\mathbf{e}_{j} \rangle d\left[\cos ( (2m+2m')\pi s )\right]\cr
	&&-\frac{1}{4m(2m-2m')\pi^2}\cr
	&&\quad\times M_t\int_0^t{M}_{s}^{-1}[\nabla \sigma(\phi_{s}^{0})|\nabla\hat \sigma(\phi_{s}^{0})] \langle [\sigma(\phi_{s}^{0})|\hat \sigma(\phi_{s}^{0})] \mathbf{e}_i,\mathbf{e}_{j} \rangle d\left[ \cos ( (2m-2m')\pi s )\right]\cr
	&=:& I_{34}+I_{35}.
	\end{eqnarray*}
	
	 Similarly, we can obtain that 
	$\sum_{0 \leq m, m'<\infty, m \neq m'}\left\|I_{34}\right\|_{p  {\mathrm{-var} }}^{2}
	<\infty$ and $\sum_{0 \leq m, m'<\infty, m \neq m'}\left\|I_{35}\right\|_{p  {\mathrm{-var} }}^{2}
	<\infty$.

	When $(f,\hat f)^\mathrm{T}_{2m,i}= \frac{\sin \left( 2m\pi t \right)}{2m\pi},  (f,\hat f)^\mathrm{T}_{2m',j}= \frac{\sin \left( 2m'\pi t \right)}{2m'\pi}$ and $m= m'$,
	\begin{eqnarray*}\label{4-60}
	\frac{\sin \left( 2m\pi t \right)}{2m\pi}d\left[ \frac{\sin \left( 2m\pi t \right)}{2m\pi} \right] 
	&=& -\frac{1}{16m^2\pi^2}d\left[ \cos ( 4m\pi t )\right],
	\end{eqnarray*}
	then it has
	\begin{eqnarray*}\label{4-61}
	R_1\big< (f,\hat f)_{2m,i},  (f,\hat f)_{2m',j} \big> _t
	&=&M_t\int_0^t{M}_{s}^{-1}[\nabla \sigma(\phi_{s}^{0})|\nabla\hat \sigma(\phi_{s}^{0})] \langle [\sigma(\phi_{s}^{0})|\hat \sigma(\phi_{s}^{0})] \mathbf{e}_i,\mathbf{e}_{j} \rangle     	\frac{\sin \left( 2m\pi t \right)}{2m\pi}d\left[ \frac{\sin \left( 2m\pi s \right)}{2m\pi} \right] \cr
	&=&-\frac{1}{16m^2\pi^2}\times M_t\int_0^t{M}_{s}^{-1}[\nabla \sigma(\phi_{s}^{0})|\nabla\hat \sigma(\phi_{s}^{0})] \langle [\sigma(\phi_{s}^{0})|\hat \sigma(\phi_{s}^{0})] \mathbf{e}_i,\mathbf{e}_{j} \rangle d\left[\cos ( 4m\pi s )\right]\cr
	&=:&I_{36}.
	\end{eqnarray*}
	Similar to the (\ref{4-30}), it has
	$\sum_{0 \leq m<\infty}\left\|I_{36}\right\|_{p  {\mathrm{-var} }}^{2}
	<\infty$ with $4-2/p>1$.
	
	According to    estimates in Step 1 to Step 4, it deduces that 
	\begin{eqnarray*}\label{4-62}
	\sum_{i, j=1}^{d_1} \sum_{m, m'=0}^{\infty}\big\|R_{1}\big\langle (f,\hat f)_{m, i}, (f,\hat f)_{m', j}\big\rangle\big\|_{p\mathrm{-var}}^{2}<\infty,
	\end{eqnarray*}
	and 	similar estimates for $	i=1,\cdots,d_1, j=d_1,\cdots,d_1+d_2$,  $	i=d_1,\cdots,d_1+d_2, j=1,\cdots,d_1$ and $	i=d_1,\cdots,d_1+d_2, j=d_1,\cdots,d_1+d_2$.
		
    The proof is completed. \qed

	\para{Proof of \thmref{thm4-1-1}}. Combine  \lemref{thm4-1-2}, \lemref{thm4-1-3}, and \lemref{thm4-1-4}, 
$\nabla^{2}(F \circ \Psi)(\gamma,\eta)\langle(\cdot, \cdot),(\cdot, \cdot)\rangle$   is a symmetric Hilbert-Schmidt bilinear form on $
{\mathcal{H}}$ for all $(\gamma,\eta)
\in{\mathcal{H}}$.

	\qed

	\subsection{A probabilistic representation of  Hessian}\label{sec-42}
	Denote $A_1$ be a self-adjoint Hilbert-Schmidt operator ${\mathcal{H}}$ which corresponds to 
    \begin{eqnarray*}
	\nabla F\left(\phi^{0}\right)\left\langle V_{1}((\cdot,\cdot), (\cdot,\cdot))\right\rangle.
    \end{eqnarray*}
    Let $A-A_1$ be a self-adjoint Hilbert-Schmidt operator ${\mathcal{H}}$ which corresponds to 
    \begin{eqnarray*}
	\nabla F\left(\phi^{0}\right)\left\langle V_{2}((\cdot,\cdot), (\cdot,\cdot))\right\rangle+\nabla^{2} F\left(\phi^{0}\right)\langle\chi((\cdot,\cdot)), \chi((\cdot,\cdot))\rangle.
    \end{eqnarray*}
    Then combining with these all, it deduces that
    \begin{eqnarray}\label{4-2-1}
    (k,\hat k) \in {\mathcal{H}} \mapsto\langle A (k,\hat k), (k,\hat k)\rangle_{
	{\mathcal{H}}}&=&\nabla F\left(\phi^{0}\right)\langle\psi((k,\hat k), (k,\hat k))\rangle\cr
    &&+\nabla^{2} F\left(\phi^{0}\right)\langle\chi((k,\hat k)), \chi((k,\hat k))\rangle
    \end{eqnarray}
extends to a continuous map on $G \Omega_{p}\left(\mathbb{R}^{d_1 +d_2}\right)$.

     Let  $\mathcal{X}$ be the closure of $\mathcal{H}$ with respect to the $p$-variation, then  $\big(\mathcal{X}, \mathcal{H}, \mu^{H}\big)$ is the abstract Wiener space.	Any
     $\langle (k,\hat k)^\mathrm{T}, (\cdot,\cdot )^\mathrm{T}\rangle \in\mathcal{H}^{*}=({\mathcal{H}^{H,d_1}})^{*}\oplus({\mathcal{H}^{\frac{1}{2},d_2}})^{*}$
     extends to a measurable linear functional on $\mathcal{X}$, which is denoted by $\langle (k,\hat k)^\mathrm{T}, (b^H,w )^\mathrm{T}\rangle$.   Denote that $\mathcal{C}_{n}=\mathcal{C}_{n}\left(\mu^{H}\right)(n=0,1,2, \ldots)$ be the $n$th Wiener chaos of $(b^H,w)^\mathrm{T}$. 
     
  For a cylinder function $F\left(b^{H},w\right)=f\big(\big\langle (k,\hat k)_{1}, (b^{H},w)^\mathrm{T}\big\rangle, \ldots,\big\langle (k,\hat k)_{d_1+d_2}, (b^{H},w)^\mathrm{T}\big\rangle\big)$, where $f:\mathbb{R}^{d_1+d_2} \to \mathbb{R}$ is a bounded smooth function with bounded derivatives,  define that 
  \begin{eqnarray*}\label{4-2-2}
  D_{(k,\hat k)} F\left(b^{H},w\right):=\sum_{j=1}^{d_1+d_2} \partial_{j} f\big(\big\langle (k,\hat k)_{1}, (b^{H},w)^\mathrm{T}\big\rangle, \ldots,\big\langle (k,\hat k)_{d_1+d_2}, (b^{H},w)^\mathrm{T}\big\rangle\big)\big((\hat k,\hat k)_{j}, (k,\hat k)\big)_{\mathcal{H}}, 
  \end{eqnarray*}
  for $(k,\hat k)^\mathrm{T} \in {\mathcal{H}} $.

  Firstly, we consider the stochastic integration of the kernel associated with $A_1$. 
    \begin{lem}\label{thm4-2-1} 
	For each fixed $t$, $V_{1}\left((b^H (m),w(m)), (b^H (m),w(m))\right)_{t}^i$, $ i\in \{1,\cdots,n\}$ converges almost surely and in $L^2(\mu^H)$ with $m\to \infty$. Morever, 
	$$\lim _{m \rightarrow \infty} V_{1}\left((b^H (m),w(m)), (b^H (m),w(m))\right)_{t}^{i}=\Theta_{t}^{i}+\Lambda_{t}^{i},$$
	where $\Lambda_{t}^{i}=\lim _{m \rightarrow \infty} \mathbb{E}\big[V_{1}\left((b^H (m),w(m)), (b^H (m),w(m))\right)_{t}^{i}\big]$ is of finite $p-$variation, and  $\Theta_{t}^{i} \in \mathcal{C}_{2}$ which corresponds to the symmetric Hilbert-Schmidt bilinear form $V_{1}\left((\cdot,\cdot), (\cdot,\cdot)\right)_{t}^{i}$.
   \end{lem}
	\para{Proof}. Consider the rough path $V'(B^H ,W)$ of $V_{1}\left((b^H ,w), (b^H ,w)\right)$, and $V'(B^H ,W)^1=V_{1}\left((b^H ,w), (b^H ,w)\right)$ being the first level path. Moreover, $V'(B^H ,W)$ and  $V_{1}\left((b^H ,w), (b^H ,w)\right)$ is of second order, it follows that
	\begin{eqnarray}\label{41-1}
	\left\|V^{\prime}(B^H ,W)^{1}\right\|_{p-\operatorname{var}}  \leq c\left(1+{\vertiii{\varepsilon(B^H,W)}^2_{p\mathrm{-var}}}\right),
	\end{eqnarray} 
	and
     \begin{eqnarray}\label{41-2}
	\big\|V'(B^H ,W)^{1}-V^{\prime}(\hat B^H ,\hat W)^{1}\big\|_{p-\operatorname{var}}  \leq c\left(1+{\vertiii{\varepsilon(B^H,W)}^c_{p\mathrm{-var}}}\right) \sum^{[p]}\big\|(B^H ,W)^{j}-(\hat B^H ,\hat W)^{j}\big\|_{p / j\mathrm{-var}},
	\end{eqnarray}
	where $c$ is some constant, and $\vertiii{(B^H,W)}_{p\mathrm{-var}}$ is defined in Section 3. With \propref{prop3-2-1}, it is easy to see that $V'(b^H(m),W(m))^{1}=V_{1}\left((b^H(m) ,w(m)), (b^H(m) ,w(m))\right)$ converges to $V_{1}\left((b^H ,w), (b^H ,w)\right)$ almost surely. Consequently, we also have
	\begin{eqnarray}\label{41-3}
    \mathbb{E}\left[\left\|V'(B^{H}(m),W(m))^{1}-V'(B^H ,W)^{1}\right\|_{p\mathrm{-var}}^{2}\right] \rightarrow 0 \quad \text { as } m \rightarrow \infty.
	\end{eqnarray}
	Likewise, with (\ref{41-1}) and \propref{prop3-2-1}, it shows that 
	\begin{eqnarray}\label{5-4}
     \mathbb{E}\left[\left\|V'(B^H ,W)^{1}\right\|_{p\mathrm{-var}}\right]<\infty.
	\end{eqnarray}
	Hence, it is clear to see that $\Lambda_{t}^{i}=\lim _{m \rightarrow \infty} \mathbb{E}\big[V_{1}\left((b^H (m),w(m)), (b^H (m),w(m))\right)_{t}^{i}\big]$ is of finite $p-$variation.
	
	Our next objective is to prove that $\Theta_{t}^{i} \in \mathcal{C}_{2}$ which corresponds to the symmetric Hilbert-Schmidt bilinear form $V_{1}\left((\cdot,\cdot), (\cdot,\cdot)\right)_{t}^{i}$.
	Due to the definition of  derivative operator, 
	\begin{eqnarray*}\label{41-5}
		D_{(k,\hat k)} \chi(b^{H}(m),w(m) )_{t} &=&M_{t} \int_{0}^{t} M_{s}^{-1} [\sigma(\phi_{s}^{0})|\hat \sigma(\phi_{s}^{0})] D_{(k,\hat k)} d (b^{H}(m),w(m) )^\mathrm{T}_{s} \cr
		&=&M_{t} \int_{0}^{t} M_{s}^{-1} [\sigma(\phi_{s}^{0})|\hat \sigma(\phi_{s}^{0})] d(k(m),\hat k(m) )^\mathrm{T}_{s}\cr
		&=&\chi(k(m),\hat k(m) )_{t} .
		\end{eqnarray*}
		As $m\to \infty$, it has $\|(k(m),\hat k(m))-(k, \hat k)\|_{q\mathrm{-var}}$ with $2>q>1$. Combined with the result that $(k,\hat k)^\mathrm{T}\in C_0^{q\mathrm{-var}}(\mathbb{R}^{d_1+d_2}) \mapsto \chi(k,\hat k ) \in C_0^{q\mathrm{-var}}(\mathbb{R}^{n})$ is a bounded linear map, it follows that $\|\chi(k(m),\hat k(m))-\chi(k, \hat k)\|_{q\mathrm{-var}} \to 0$ with $m\to \infty$.

		Similarly, we have 
\begin{eqnarray*}\label{41-6}	
	&&\frac{1}{2} D_{(k,\hat k)} V_{1}\left((b^{H}(m),w(m) ), (b^{H}(m),w(m) )\right)_{t} \cr
	&=& M_{t} \int_{0}^{t} M_{s}^{-1}\left\{[\nabla \sigma(\phi_{s}^{0})|\nabla\hat \sigma(\phi_{s}^{0})]\big\langle D_{(k,\hat k)} \chi(b^{H}(m),w(m) )_{s}, d (b^{H}(m),w(m) )^\mathrm{T}_{s}\big\rangle\right. \cr
	&&\left.\quad+[\nabla \sigma(\phi_{s}^{0})|\nabla\hat \sigma(\phi_{s}^{0})]\big\langle\chi(b^{H}(m),w(m) )_{s}, D_{(k,\hat k)} d (b^{H}(m),w(m) )^\mathrm{T}_{s}\big\rangle\right\} \cr
	&=&M_{t} \int_{0}^{t} M_{s}^{-1}\left\{[\nabla \sigma(\phi_{s}^{0})|\nabla\hat \sigma(\phi_{s}^{0})]\big\langle\chi(k(m),\hat k(m) )_{s}, d (b^{H}(m),w(m) )^\mathrm{T}_{s}\big\rangle\right. \cr
	&&\left.\quad+[\nabla \sigma(\phi_{s}^{0})|\nabla\hat \sigma(\phi_{s}^{0})]\big\langle\chi(b^{H}(m),w(m) )_{s}, d (k(m),\hat k(m) )^\mathrm{T}_{s}\big\rangle\right\} \cr
	&=&V_{1}\big((k(m),\hat k(m) ), (b^{H}(m),w(m) )\big)_{t} ,
		\end{eqnarray*}
when $m \to \infty$, we have $||(b^{H}(m),w(m) )-(b^{H},w )||_{p\mathrm{-var}} \to 0$ almost surely and in $L^r$ for any $r>0$, it is clear to see that $\frac{1}{2} D_{(k,\hat k)} V_{1}\left((b^{H}(m),w(m) ), (b^{H}(m),w(m) )\right)_{t} =V_{1}\big((k(m),\hat k(m) ), (b^{H}(m),w(m) )\big)_{t} \to V_{1}\big((k,\hat k ), (b^{H},w)\big)_{t}$ almost surely as $m \to \infty$.

 In the same manner, it deduces that
$$\frac{1}{2} D_{(\tilde k,\tilde{\hat k})}D_{(k,\hat k)} V_{1}\left((b^{H}(m),w(m) ), (b^{H}(m),w(m) )\right)_{t} =V_{1}\big((k(m),\hat k(m) ), (\tilde k(m),\tilde{\hat k}(m) )\big)_{t} \to V_{1}\big((k,\hat k ), (\tilde k,\tilde{\hat k})\big)_{t},$$ almost surely as $m \to \infty$. 

Meanwhile,	by the definition of  derivative operator, it has 
\begin{eqnarray*}\label{41-7}
\frac{1}{2} D_{(\tilde k,\tilde{\hat k})}D_{(k,\hat k)} V'\left(B^{H},W\right)_{t}^{1, i}=V_{1}\big((k,\hat k ), (\tilde k,\tilde{\hat k})\big)_{t}^{i}\quad i\in\{1,\cdots,n\},
\end{eqnarray*}
and
\begin{eqnarray*}\label{41-8}
\frac{1}{2} D_{(k,\hat k)} V'\left(B^{H},W\right)_{t}^{1, i}=V_{1}\big((k,\hat k), (b^H,w)\big)_{t}^{i},\quad i\in\{1,\cdots,n\},
\end{eqnarray*}
which claims that all  elements in the first Wiener chaos $\mathcal{C}_{1}$ equal to $0$. Furthermore, combined with \thmref{thm4-1-1}, and the fact  that the second Wiener chaos $\mathcal{C}_{2}$ is unitarily
isometric with the space of symmetric Hilbert-Schmidt operators in $\mathcal{H}\otimes \mathcal{H}$, it suffices to show that $\Theta_{t}^{i} \in \mathcal{C}_{2}$ which corresponds to the symmetric Hilbert-Schmidt bilinear form $V_{1}\left((\cdot,\cdot), (\cdot,\cdot)\right)_{t}^{i}$.
	
    The proof is completed. \qed	
    
    \begin{lem}\label{thm4-2-3} 
    	Let $p'>p$ and $F: C_{0}^{p'-\operatorname{var}}\left(\mathbb{R}^{n}\right)$ be a Fr\'echet differentiable function, then $\nabla F\left(\phi^{0}\right)\langle\Theta\rangle \in \mathcal{C}_{2}\left(\mu^{H}\right)$ which is related to the symmetric Hilbert-Schmidt bilinear form $\nabla F\left(\phi^{0}\right) \circ V_{1}=\nabla F\left(\phi^{0}\right)\left\langle V_{1}((\cdot,\cdot), (\cdot,\cdot))\right\rangle$ in Cameron-Martin space $
    	{\mathcal{H}}$.
    \end{lem}
    \para{Proof}. Denote $g_K$ the elements of $\mathcal{C}_2(\mu^H)$, where $K$ is sysmetric Hilbert-Schmidt bilinear form, define that
    \begin{eqnarray}\label{4-2-3}
    M:=\left\{\vartheta \in C_{0}^{p-\operatorname{val}}\left(\mathbb{R}^{n}\right)^{*} \mid \vartheta\langle\Theta(w)\rangle=g_{\vartheta \circ V_{1}}(w) \text { a.a.} w\left(\mu^{H}\right)\right\},
    \end{eqnarray}
    $M$ is a linear subspace. From \lemref{thm4-1-3} and \lemref{thm4-1-4}, if $\vartheta_1\in M$ converges to $\vartheta\in M$ in the $weak^\star$-limit as $l\to \infty$, it deduces that ${\vartheta_l \circ V_{1}}$ converges to ${\vartheta \circ V_{1}}$ almost surely as $l\to \infty$.  Therefore, $M$ is closed under $weak^\star$-limit.
    
    Since $y_{k / 2^{m}}^{i}\left(1 \leq k \leq 2^{m}, 1 \leq i \leq n\right)$ is the dyadic approximation of $y \in C_{0}^{p\mathrm{-var}}\left(\mathbb{R}^{n}\right)$,  $\nabla F\left(\phi^{0}\right)\langle y(m) \rangle$ is the linear combination of $y_{k / 2^{m}}^{i}\left(1 \leq k \leq 2^{m}, 1 \leq i \leq n\right)$. Furthermore, $\nabla F\left(\phi^{0}\right) \circ \pi(m)\in M$.
    Due to that $\nabla F\left(\phi^{0}\right) \circ \pi(m) \rightarrow \nabla F\left(\phi^{0}\right) $ in the $weak^\star$-limit, it follows that $\nabla F\left(\phi^{0}\right)\in M$. 
    The proof is completed. \qed
    
      Next, we turn to the stochastic integration of the kernel associated with $A-A_1$.
	\begin{lem}\label{thm4-2-2} 
		$A-A_1$ being a self-adjoint Hilbert-Schmidt operator ${\mathcal{H}}$ which corresponds to 
		\begin{eqnarray*}\label{4-2-4}
			\nabla F\left(\phi^{0}\right)\left\langle V_{2}((\cdot,\cdot), (\cdot,\cdot))\right\rangle+\nabla^{2} F\left(\phi^{0}\right)\langle\chi((\cdot,\cdot)), \chi((\cdot,\cdot))\rangle.
		\end{eqnarray*}
		is of trace class, moreover, 
		\begin{eqnarray}\label{4-2-5}
		\left\langle\left(A-A_{1}\right) (B^{H},W), (B^{H},W)\right\rangle&:=&			\nabla F\left(\phi^{0}\right)\left\langle V_{2}((B^{H},W), (B^{H},W))\right\rangle\cr
		&&+\nabla^{2} F\left(\phi^{0}\right)\langle\chi(B^{H},W), \chi(B^{H},W)\rangle,
		\end{eqnarray}
		is the sum of $Tr(A-A_1)$ and the second Wiener chaos corresponding to $A-A_1$ which is denoted by $\hat{K}_{A-{A_1}}\left((B^H, W)^{1}\right)$.
	\end{lem}
	\para{Proof}. With the aid of \lemref{thm4-1-2}, $A-A_1$, which is a self-adjoint Hilbert-Schmidt operator ${\mathcal{H}}$ 
	is of trace class  for a Fr\'echet differentiable function $F$. It yields that $\left\langle\left(A-A_{1}\right) (B^{H},W), (B^{H},W)\right\rangle$ 		can be rewritten as  the sum of $Tr(A-A_1)$ and the second Wiener chaos corresponding to $A-A_1$,
	\begin{eqnarray}\label{4-2-6}
    \left\langle\left(A-A_{1}\right) (B^{H},W), (B^{H},W)\right\rangle=\hat{K}_{A-{A_1}}\left((B^H, W)^{1}\right)+\operatorname{Tr}(A-{A_1}) .
	\end{eqnarray}
	
	The proof is completed. \qed

  \begin{thm}\label{thm4-2-4} 
	Let $\alpha >1$ be such that $\mathrm{Id}_{
		{\mathcal{H}}}+\alpha A$ is strictly positive in the form sense, then 
	\begin{eqnarray}\label{4-2-7}
	&&\int_{G \Omega_{p}\left(\mathbb{R}^{d_1+d_2}\right)} \exp \big(-\frac{\alpha}{2}\langle A (B^H, W), (B^H, W)\rangle\big) \mathbb{P}^{H}(d (B^H, W))\cr
	&=&\exp[-\frac{\alpha}{2}(\operatorname{Tr}\left(A-A_{1}\right)+\nabla F\left(\phi^{0}\right)\langle\Lambda\rangle)]\cdot\operatorname{det}\big(\mathrm{Id}_{
		{\mathcal{H}}}+\alpha A\big)^{-1 / 2},
	\end{eqnarray}
	where
	 $\operatorname{det}_{2}$ denotes  the Carleman-Fredholm determinant.  
  \end{thm}
	\para{Proof}. 
   Firstly, we have
	\begin{eqnarray*}\label{4-2-8}
	&\nabla F(\phi) \psi(B^H, W)+\nabla^{2} F(\phi)\left\langle \chi(B^H, W), \chi(B^H, W)\right\rangle \cr
	&\quad=\nabla F(\phi)\langle V_1(B^H, W)\rangle+\langle(A-{A_1}) (B^H, W), (B^H, W)\rangle,
	\end{eqnarray*}
	where
	\begin{eqnarray*}\label{4-2-9}
	\left\langle\left(A-A_{1}\right) (B^{H},W), (B^{H},W)\right\rangle&:=&			\nabla F\left(\phi^{0}\right)\left\langle V_{2}((B^{H},W), (B^{H},W))\right\rangle\cr
	&&+\nabla^{2} F\left(\phi^{0}\right)\langle\chi(B^{H},W), \chi(B^{H},W)\rangle,
	\end{eqnarray*}
	then by \lemref{thm4-2-2} and \lemref{thm4-2-3}, it follows that 	$\left\langle\left(A-A_{1}\right) (B^{H},W), (B^{H},W)\right\rangle$	is the sum of $Tr(A-A_1)$ and the second Wiener chaos corresponding to $A-A_1$ which is denoted by $\hat{K}_{A-{A_1}}\left((B^H, W)^{1}\right)$, and $\nabla F\left(\phi^{0}\right)\langle\Theta\rangle \in \mathcal{C}_{2}\left(\mu^{H}\right)$. Hence, define that $\left\{\lambda_{j}\right\} $  and  $\{\varkappa_{j}\}$ are eigenvalues and corresponding (orthonormal) eigenvectors of $A-{A_1}$, and $\{\hat \lambda_{j}\} $  and  $\{\hat \varkappa_{j}\}$ are eigenvalues and corresponding (orthonormal) eigenvectors of $A$, we have
	\begin{eqnarray*}\label{4-2-10}
	&&\nabla F(\phi)\langle V_1(B^H, W)\rangle+\langle(A-{A_1}) (B^H, W), (B^H, W)\rangle\cr 
	&&=\nabla F(\phi)\langle V_1(B^H, W)\rangle+\sum_{j} \lambda_{j}\left\langle e_{j}, (B^H, W)_{1}\right\rangle^2 \cr
	&&=\nabla F(\phi)\langle V_1(B^H, W)\rangle+\hat{K}_{A-{A_1}}\left((B^H, W)^{1}\right)+\operatorname{Tr}(A-{A_1}) \cr
	&&=\sum \hat{\lambda}_{j}\big(\left\langle\hat{\varkappa}_{j}, (B^H, W)^{1}\right\rangle^2-1\big)+\operatorname{Tr}(A-{A_1})+\nabla F\left(\phi^{0}\right)\langle\Lambda\rangle,
	\end{eqnarray*}
	where $\left\{\hat{\varkappa}_{j}\right\}_{j=1,2, \ldots}$ are ONB of Cameron-Martin space ${\mathcal{H}}$, and  $\Lambda$ defined in \lemref{thm4-2-1} is of finite $p$-variation. Take the above proof into consideration, (\ref{4-2-7}) is obtained.
	
	The proof is completed. \qed

	\section{Main proof}\label{sec-5}
   
	Now, we give the proof of our main result.
	\para{Proof of \thmref{thm2-1}}. By Assumption {\rm (A2)}, $F_{\Lambda}:=F\circ\Phi+||(\cdot,\cdot)||^2_{\mathcal{H}} /2$ attains its minimum at a unique point $(\gamma,\eta)^\mathrm{T}
	\in{\mathcal{H}}$. 
	In the nerghborhood $O \subset G \Omega_{p}\left(\mathbb{R}^{d_1+d_2}\right)$ of $(\gamma,\eta)$, there exists a constant $C$ such that
    \begin{eqnarray}\label{5-0}
    \int_{O^{c}} \exp \big(-F\big(\hat{\Phi}_{\varepsilon}((B^H,W), \lambda)^{1}\big) / \varepsilon^{2}\big) \mathbb{P}_{\varepsilon}^{H}(d (B^H,W)) \leq C e^{-(a+\delta) / \varepsilon^{2}},
    \end{eqnarray}
    for $\varepsilon \in(0,1]$. The above term (\ref{5-0}) converges to zero as $\varepsilon \to 0$.
	
	According to the \lemref{thm3-1-1}, set $O=\gamma+U_{\rho}$ where $U_{\rho}=\left\{(B^H,W) \in G \Omega_{p}\left(\mathbb{R}^{d_1 +d_2}\right) \mid \vertiii{(B^H,W)}_{p\mathrm{-var}}<\rho\right\}$ for $\rho>0$, 
	\begin{eqnarray}\label{5-1}
		&&\int_{\gamma+U_{\rho}}  \exp \big(-F\big(\hat{\Phi}_{\varepsilon}((B^H,W), \lambda)^{1}\big) / \varepsilon^{2}\big) \mathbb{P}_{\varepsilon}^{H}(d (B^H,W))\cr 
		&=& \int_{U_{\rho}} \exp \big(-F\big(\hat{\Phi}_{\varepsilon}((B^H,W)+(\gamma,\eta)^\mathrm{T}, \lambda)^{1}\big) / \varepsilon^{2}\big) \cr
		&& \quad\times \exp \bigg(-\frac{1}{\varepsilon^{2}}\left\langle(\gamma,\eta)^\mathrm{T}, (B^H,W)^{1}\right\rangle-\frac{1}{2 \varepsilon^{2}}\|(\gamma,\eta)\|_{\mathcal{H}}^{2}\bigg) \mathbb{P}_{\varepsilon}^{H}(d (B^H,W))\cr
		&=&\int_{\{\vertiii{(B^H,W)}_{p\mathrm{-var}}<\rho\}} \exp \bigg(-\frac{F\left(\phi^{(\varepsilon)}\right)}{\varepsilon^{2}}-\frac{1}{\varepsilon}\left\langle(\gamma,\eta)^\mathrm{T}, (B^H,W)^{1}\right\rangle\bigg.\cr
		&&\quad\bigg.-\frac{1}{2 \varepsilon^{2}}\|(\gamma,\eta)\|_{\mathcal{H}}^{2}\bigg) \mathbb{P}^{H}(d(B^H,W)).
	\end{eqnarray} 
	Take the stochastic Taylor expansion for $\phi^{(\varepsilon)}$ in the neighbourhood of $\phi^{(0)}$, one can get
	\begin{eqnarray} \label{5-2}
		F(\phi^{(\varepsilon)})=& F\left(\phi^{0}\right)+\nabla F\left(\phi^{0}\right)\left\langle\phi^{(\varepsilon)}-\phi^{0}\right\rangle+\frac{1}{2} \nabla^{2} F\left(\phi^{0}\right)\left\langle\phi^{(\varepsilon)}-\phi^{0}, \phi^{(\varepsilon)}-\phi^{0}\right\rangle \cr &+\frac{1}{6} \int_{0}^{1} d \theta \nabla^{3} F\left(\theta \phi^{(\varepsilon)}+(1-\theta) \phi^{0}\right)\left\langle\phi^{(\varepsilon)}-\phi^{0}, \phi^{(\varepsilon)}-\phi^{0}, \phi^{(\varepsilon)}-\phi^{0}\right\rangle \cr
		=& F\left(\phi^{0}\right)+\nabla F\left(\phi^{0}\right)\left\langle\varepsilon \phi^{1}+\varepsilon^{2} \phi^{2}\right\rangle+\frac{1}{2} \nabla^{2} F\left(\phi^{0}\right)\left\langle\varepsilon \phi^{1}, \varepsilon \phi^{1}\right\rangle+Q_{\varepsilon}^{3} .
	\end{eqnarray}
	Due to the result in \lemref{lem3-2-2},  there exists some constanct $C>0$ such that $\left|Q_{\varepsilon}^{3}\right| \leq C(\varepsilon+\vertiii{\varepsilon(B^H,W)}_{p\mathrm{-var}})^{3} $ on the set $\big\{\vertiii{\varepsilon(B^H,W)}_{p\mathrm{-var}}<\rho_{0}\big\}$.
	
	Consider the term with order $-2$ in (\ref{5-1}), due to the Assumption (A2)
	$$-\frac{1}{\varepsilon^{2}}\big(F\left(\phi^{0}\right)+\frac{1}{2}\|(\gamma,\eta)\|_{\mathcal{H}}^{2}\big)=-\frac{a}{\varepsilon^{2}}.$$
	For the term of order $-1$ in (\ref{5-1}), we have $\phi^1:G\Omega_p(\mathbb{R}^{d_1+d_2})\mapsto C_0^{p\mathrm{-var}}(\mathbb{R}^n)$, and it satisfies the following differential equation,
	\begin{eqnarray}\label{5-81}
		&&d \phi_{t}^{1} -[\nabla \sigma(\phi_{s}^{0})|\nabla\hat \sigma(\phi_{s}^{0})]\left\langle\phi_{t}^{1}, d (\gamma_{t},\eta_{t})^\mathrm{T}\right\rangle-\nabla_{y} \beta\left(0, \phi_{t}^{0}\right)\left\langle\phi_{t}^{1},1\right\rangle d t \cr
		&=&[\sigma(\phi_{t}^{0})|\hat \sigma(\phi_{t}^{0})] d (k_{t},\hat k_t)^\mathrm{T}+\nabla_{\varepsilon} \beta\left(0, \phi_{t}^{0}\right) d t, \quad \phi_{0}^{1}=0
		\end{eqnarray}
	Then, set  $\theta^1(k,\hat k)_{t}=\phi^1(k,\hat k)_t-\chi(k,\hat k)_{t}$, where 	\begin{eqnarray}\label{5-20}
	\chi(k,\hat k)_{t}=M_{t} \int_{0}^{t} M_{s}^{-1} [\sigma(\phi_{s}^{0})|\hat \sigma(\phi_{s}^{0})] d (k_{s},\hat k_s)^\mathrm{T}.
	\end{eqnarray} 
	By the  Assumption (A2), we conclude that
	\begin{eqnarray}\label{5-8}
	\langle (k,\hat k), (\gamma,\eta)\rangle_{\mathcal{H}}+\nabla F\left(\phi^{0}\right)\langle\chi(k,\hat k)\rangle=0,
	\end{eqnarray}
	
	Then, $\theta^1(k,\hat k)_{t}$ satisfies the following differential equation,
	\begin{eqnarray}
	&&d \theta_{t}^{1}-[\nabla \sigma(\phi_{s}^{0})|\nabla\hat \sigma(\phi_{s}^{0})]\left\langle\theta_{t}^{1}, d (\gamma_{t},\eta_{t})^\mathrm{T}\right\rangle-\nabla_{y} \beta\left(0, \phi_{t}^{0}\right)\left\langle\theta_{t}^{1},1\right\rangle d t \cr
	&=&\nabla_{\varepsilon} \beta\left(0, \phi_{t}^{0}\right) d t, \quad \theta_{0}^{1}=0,
	\end{eqnarray}
	its solution can be rewritten as 
$
	\theta_{t}^{1}=M_{t} \int_{0}^{t} M_{s}^{-1} \nabla_{\varepsilon} \beta\left(0, \phi_{s}^{0}\right) d s,
$
	then, $\theta^1(B^H,W)_{t}$ is independent of $(B^H,W)$ and  of finite $q$-variation.
	
Hence, with (\ref{5-8}) and straightforward computation, it deduces that 
	$$-\frac{1}{\varepsilon}[\langle (k,\hat k), (\gamma,\eta)\rangle_{\mathcal{H}}+\nabla F\left(\phi^{0}\right)\langle\phi^1(k,\hat k)\rangle]=-\frac{\nabla F\left(\phi^{0}\right)\langle\theta^1(k,\hat k)\rangle}{\varepsilon}.$$
	Next, we focus on  the term of order $0$ in (\ref{5-1}), it has
	\begin{eqnarray}\label{5-31}
		&&d \phi_{t}^{2}- [\nabla \sigma(\phi_{s}^{0})|\nabla\hat \sigma(\phi_{s}^{0})]\left\langle\phi_{t}^{2}, d (\gamma_{t},\eta_{t})^\mathrm{T}\right\rangle-\nabla_{y} \beta\left(0, \phi_{t}^{0}\right)\left\langle\phi_{t}^{2},1\right\rangle d t \cr
		&=& [\nabla \sigma(\phi_{s}^{0})|\nabla\hat \sigma(\phi_{s}^{0})]\big\langle\phi_{t}^{1}, d (k_{t},\hat k_t)^\mathrm{T}\big\rangle+\frac{1}{2} [\nabla^2 \sigma(\phi_{t}^{0})|\nabla^2\hat \sigma(\phi_{t}^{0})]\left\langle\phi_{t}^{1}, \phi_{t}^{1}, d (\gamma_{t},\eta_{t})^\mathrm{T}\right\rangle \cr
		&&+\frac{1}{2} \nabla_{y}^{2} \beta\left(0,\phi_{t}^{0}\right)\left\langle\phi_{t}^{1}, \phi_{t}^{1}\right\rangle d t +\nabla_{y} \nabla_{\varepsilon} \beta\left(0,\phi_{t}^{0}\right)\left\langle\phi_{t}^{1},1\right\rangle d t\cr
		&&+\frac{1}{2} \nabla_{\varepsilon}^{2} \beta\left(0, \phi_{t}^{0}\right) d t, \quad \phi_{0}^{2}=0 .
		\end{eqnarray}
	Then,  $\chi$ and $\psi$ extend to continuous maps from $G\Omega_p(\mathbb{R}^{d_1+d_2})$, and rewrite that $\chi(B^H,W)$ and $\psi((B^H,W), (B^H,W))$.
	
	Consequently, we set $\theta^2(k,\hat k)_{t}=\phi^{2}(k,\hat k)_t-\psi^{2}((k,\hat k), (k,\hat k))_t / 2$
, where	\begin{eqnarray} \label{5-3}
	&&d\psi _t- [\nabla \sigma(\phi_{t}^{0})|\nabla\hat \sigma(\phi_{t}^{0})] \langle \psi _t,d(\gamma_t,\eta_t)^\mathrm{T} \rangle -\nabla_y \beta\left(0, \phi _{t}^{0} \right) \langle \psi _t,1 \rangle dt\cr
	&=&2[\nabla \sigma(\phi_{t}^{0})|\nabla\hat \sigma(\phi_{t}^{0})] \big< \chi(k,\hat k)_t,d(k_t,\hat k_t)^\mathrm{T} \big> \cr
	&&+[\nabla^2 \sigma(\phi_{t}^{0})|\nabla^2\hat \sigma(\phi_{t}^{0})] \big< \chi(k,\hat k)_t,\chi(k,\hat k)_t,d(\gamma_t,\eta_t)^\mathrm{T} \big>\cr
	&&+\nabla ^2\beta_0\left( \phi _{t}^{0} \right) \big< \chi(k,\hat k) _t,\chi(k,\hat k) _t \big> dt,\quad
	\psi _0=0.
	\end{eqnarray}
And $\theta^2(k,\hat k)_{t}$	statifies the differential equation as follows,	
\begin{eqnarray} \label{5-32}
		&&d \theta_{t}^{2}- [\nabla \sigma(\phi_{t}^{0})|\nabla\hat \sigma(\phi_{t}^{0})]\left\langle\theta_{t}^{2}, d (\gamma_t,\eta_t)^\mathrm{T}\right\rangle-\nabla_{y} \beta\left(0, \phi_{t}^{0}\right)\left\langle\theta_{t}^{2},1\right\rangle d t \cr
		&=& [\nabla \sigma(\phi_{t}^{0})|\nabla\hat \sigma(\phi_{t}^{0})]\left\langle\theta_{t}^{1}, d (k_t,\hat k_t)^\mathrm{T} \right\rangle+\frac{1}{2} [\nabla^2 \sigma(\phi_{t}^{0})|\nabla^2\hat \sigma(\phi_{t}^{0})]\left\langle\theta_{t}^{1}, \theta_{t}^{1}, d (\gamma_t,\eta_t)^\mathrm{T}\right\rangle\cr
		&&+[\nabla^2 \sigma(\phi_{t}^{0})|\nabla^2\hat \sigma(\phi_{t}^{0})]\left\langle\theta_{t}^{1}, \chi_{t}, d (\gamma_t,\eta_t)^\mathrm{T}\right\rangle \cr
		&&+\frac{1}{2} \nabla_{y}^{2} \beta\left(\phi_{t}^{0}\right)\left\langle\theta_{t}^{1}, \theta_{t}^{1}\right\rangle d t+\nabla_{y}^{2} \beta\left(\phi_{t}^{0}\right)\left\langle\theta_{t}^{1}, \chi_{t}\right\rangle d t \cr
		&&+\nabla_{y} \nabla_{\varepsilon} \beta\left(\phi_{t}^{0}\right)\left\langle\theta_{t}^{1}+\chi_{t}\right\rangle d t+\frac{1}{2} \nabla_{\varepsilon}^{2} \beta\left(0, \phi_{t}^{0}\right) d t, \quad \theta_{0}^{2}=0 .
		\end{eqnarray}
		Equivalently, its solution can be rewritten as 
	\begin{eqnarray}
		\theta_{t}^{2}&=&M_{t} \int_{0}^{t} M_{s}^{-1}\big[ [\nabla \sigma(\phi_{s}^{0})|\nabla\hat \sigma(\phi_{s}^{0})]\big\langle\theta_{s}^{1}, d (k_{s},\hat k_s)^\mathrm{T}\big\rangle+\frac{1}{2} [\nabla^2 \sigma(\phi_{s}^{0})|\nabla^2\hat \sigma(\phi_{s}^{0})]\left\langle\theta_{s}^{1}, \theta_{s}^{1}, d (\gamma_{s},\eta_{s})^\mathrm{T}\right\rangle \big.\cr
		&&\quad+[\nabla^2 \sigma(\phi_{s}^{0})|\nabla^2\hat \sigma(\phi_{s}^{0})]\left\langle\theta_{s}^{1}, \chi_{s}, d (\gamma_{s},\eta_{s})^\mathrm{T}\right\rangle \cr
		&&\quad+\frac{1}{2} \nabla_{y}^{2} \beta\left(\phi_{s}^{0}\right)\left\langle\theta_{s}^{1}, \theta_{s}^{1}\right\rangle d s+\nabla_{y}^{2} \beta\left(\phi_{s}^{0}\right)\left\langle\theta_{s}^{1}, \chi_{s}\right\rangle d s \cr
		&&\quad\big.+\nabla_{y} \nabla_{\varepsilon} \beta\left(\phi_{s}^{0}\right)\left\langle\theta_{s}^{1}+\chi_{s}\right\rangle d s+\frac{1}{2} \nabla_{\varepsilon}^{2} \beta\left(0, \phi_{s}^{0}\right) d s\big] .
		\end{eqnarray}
	So, $\theta^2((k,\hat k), (k,\hat k))$ extends a map from $G\Omega_p(\mathbb{R}^{d_1+d_2})$ to $ C_0^{p\mathrm{-var}}(\mathbb{R}^n)$. Moreover, from \lemref{lem3-2-2},  for some constant $C>0$, it has  $\left\|\theta^{2}(B^H,W)\right\|_{p-\operatorname{var}} \leq C(1+\vertiii{(B^H,W)}_{p\mathrm{-var}})$ with $(B^H,W)\in G\Omega_p(\mathbb{R}^{d_1+d_2})$. By the  \lemref{prop3-2-3}, it follows that $\theta^2(B^H,W)$ is  exponential integrable.
	
 By \lemref{lem3-2-2}, \lemref{prop3-2-3} and \thmref{thm4-1-1}, it deduces that
	$$\exp \big(-\nabla F\left(\phi^{0}\right)\left\langle\phi^{2}\right\rangle-\frac{1}{2} \nabla^{2} F\left(\phi^{0}\right)\left\langle\phi^{1}, \phi^{1}\right\rangle\big) \in L^{r}\left(G \Omega_{p}\left(\mathbb{R}^{d_1 +d_2}\right), \mathbb{P}^{H}\right),\quad r>1.$$
	When $\varepsilon \le \rho$,   it has
	\begin{eqnarray}\label{5-5}
		&\mathbf{1}_{\{\vertiii{\varepsilon(B^H,W)}_{p\mathrm{-var}}<\rho\}} \exp \big(-\nabla F\left(\phi^{0}\right)\left\langle\phi^{2}\right\rangle-\frac{1}{2} \nabla^{2} F\left(\phi^{0}\right)\left\langle\phi^{1}, \phi^{1}\right\rangle\big) \exp \left(-\varepsilon^{-2} Q_{\varepsilon}^{3}\right) \cr
		&\leq \exp \big(-\nabla F\left(\phi^{0}\right)\left\langle\phi^{2}\right\rangle-\frac{1}{2} \nabla^{2} F\left(\phi^{0}\right)\left\langle\phi^{1}, \phi^{1}\right\rangle\big) \exp \left[2 C \rho(1+\vertiii{(B^H,W)}_{p\mathrm{-var}})^{2}\right].
	\end{eqnarray}
	The dominated convergence theorem yields that
	\begin{eqnarray}\label{5-6}
	\lim _{\varepsilon \to 0} \int_{\{\vertiii{\varepsilon(B^H,W)}_{p\mathrm{-var}}<\rho\}} \exp \big(-\nabla F\left(\phi^{0}\right)\left\langle\phi^{2}\right\rangle-\frac{1}{2} \nabla^{2} F\left(\phi^{0}\right)\left\langle\phi^{1}, \phi^{1}\right\rangle-\frac{1}{\varepsilon^{2}} Q_{\varepsilon}^{3}\big) \mathbb{P}^{H}(d (B^H,W)) \cr
	\quad=\int_{G \Omega_{p}\left(\mathbb{R}^{d_1 +d_2}\right)} \exp \big(-\nabla F\left(\phi^{0}\right)\left\langle\phi^{2}\right\rangle-\frac{1}{2} \nabla^{2} F\left(\phi^{0}\right)\left\langle\phi^{1}, \phi^{1}\right\rangle\big) \mathbb{P}^{H}(d (B^H,W)).
	\end{eqnarray}
	By \lemref{thm4-2-4}, the right-hand side of (\ref{5-6}) exists, then the coefficient $\alpha_0$ can be defined in the following sense,
	\begin{eqnarray}\label{5-7}
	\alpha_0=\exp[-\frac{1}{2}(\operatorname{Tr}\left(A-A_{1}\right)+\nabla F\left(\phi^{0}\right)\langle\Lambda\rangle)]\cdot\operatorname{det}\big(\mathrm{Id}_{
		{\mathcal{H}}}+ A\big)^{-1 / 2},
	\end{eqnarray} 
	where $\Lambda$ is of finite $p-$variation.
	
	The proof is completed. \qed

	\section*{Acknowledgments}
     This work was partly supported by  the NSF of China (Grant 12072264),    the Shaanxi Provincial Key R\&D Program (Grants 2020KW-013, 2019TD-010).


\section*{References}
	
	\begin{thebibliography}{}
		 \bibitem{2003Adams}Adams A., Fournier J.,  Sobolev spaces. Elsevier, 2003.
		 
		 
		 \bibitem{1982Azencott}Azencott R.,  Formule de Taylor stochastique et développement asymptotique d’intégrales de Feynmann. S\'eminaire de Probabilités XVI, 1980/81 Supplément: Géométrie Diff\'erentielle Stochastique. Springer, Berlin, Heidelberg, 1982: 237-285.
		 
		 \bibitem{2006Aida}Aida S.,  Notes on proofs of continuity theorem in rough path analysis. Preprint of Osaka University, 2006.
		 

		
		
		\bibitem{1988Arous}Ben Arous G., Methods de Laplace et de la phase stationnaire sur l'espace de Wiener. Stochastics, 1988, 25(3): 125-153.
		
		\bibitem{2008Biagini}Biagini F., Hu Y.,  Oksendal B., et al. Stochastic calculus for fractional Brownian motion and applications. Springer Science \& Business Media, 2008.
		
		\bibitem{2015Bailleul}Bailleul I.,  Flows driven by rough paths. Revista matematica iberoamericana, 2015, 31(3): 901-934.
		
	    \bibitem{2002Coutin}Coutin L., 	Qian Z., Stochastic analysis, rough path analysis and fractional Brownian motions. Probability Theory and Related Fields, 2002, 122(1): 108-140.
	    
	    \bibitem{1998Dembo}Dembo A., Zeitouni O., Large Deviations Techniques and Applications, second Ed., Springer-Verlag, 1998.
	    
	   \bibitem{2021Friz}Friz P., Gassiat P., Pigato P. Precise asymptotics: robust stochastic volatility models. The Annals of Applied Probability, 2021, 31(2): 896-940.

        \bibitem{2020Friz}Friz P.,  Hairer M.,  A course on rough paths. Springer International Publishing, 2020.	

        \bibitem{2007Friz}Friz P.,  Victoir N., Large deviation principle for enhanced Gaussian processes. Annales de l'Institut Henri Poincare (B) Probability and Statistics,  2007, 43(6): 775-785.


         \bibitem{2010Friz}Friz P.,  Victoir N., Multidimensional stochastic processes as rough paths: theory and applications. Cambridge University Press, 2010.	


        \bibitem{2022Friz}Friz P., Klose T.,  Precise Laplace asymptotics for singular stochastic PDEs: The case of 2D gPAM. Journal of Functional Analysis, 2022, 109446.  
        
        

        
        \bibitem{2006InahamaKawabi}Inahama Y., Kawabi H.,  Large deviations for heat kernel measures on loop spaces via rough paths. Journal of the London Mathematical Society, 2006, 73(3): 797-816.
        
        
        
        \bibitem{2006Inahama}Inahama Y.,  Laplace's method for the laws of heat processes on loop spaces. Journal of Functional Analysis, 2006, 232(1): 148-194.	
        
        
        \bibitem{2007Inahama}Inahama Y., Kawabi H., Asymptotic expansions for the Laplace approximations for It\^o functionals of Brownian rough paths. Journal of Functional Analysis, 2007, 243(1): 270-322.
        
        \bibitem{2013Inahama}Inahama Y.,  Laplace approximation for rough differential equation driven by fractional Brownian motion. The Annals of Probability, 2013, 41(1): 170-205.
        
        
        
        \bibitem{2010Inahama}Inahama Y.,   A stochastic Taylor-like expansion in the rough path theory. Journal of Theoretical Probability, 2010, 23(3): 671-714.  
				 
				 
		 \bibitem{1997Janson}Janson S.,  Gaussian Hilbert Spaces. Cambridge University Press, Cambridge, 1997.
		
		\bibitem{1991Kusuoka}Kusuoka S.,  Stroock D., Precise asymptotics of certain Wiener functionals. Journal of Functional Analysis, 1991, 99(1): 1-74.
		
		\bibitem{2008Osajima}Kusuoka S., Osajima Y.,  A remark on the asymptotic expansion of density function of Wiener functionals. Journal of Functional Analysis, 2008, 255(9): 2545-2562.	
			 
		\bibitem{2002Ledoux}Ledoux M., Qian Z.,  Zhang T.,  Large deviations and support theorem for diffusion processes via rough paths. Stochastic Processes and Their Applications, 2002, 102(2): 265-283.	
		
		 \bibitem{1998Lyons}Lyons T.,  Differential equations driven by rough signals. Revista Matemática Iberoamericana, 1998, 14(2): 215-310.
		
		\bibitem{2002Lyons}Lyons T., Qian Z.,  System control and rough paths. Oxford University Press, 2002.
		
		    		 
	    \bibitem{2006Millet}Millet A., Sanz-Solé M.,  Large deviations for rough paths of the fractional Brownian motion. Annales de l'IHP Probabilités et statistiques, 2006, 42(2): 245-271.
	    
	   \bibitem{2008Mishura}Mishura Y.,  Stochastic calculus for fractional Brownian motion and related processes. Springer Science \& Business Media, 2008.
	    
	     \bibitem{2021Pei}Pei B., Inahama Y., Xu Y.,  Averaging principle for fast-slow system driven by mixed fractional Brownian rough path. Journal of Differential Equations, 2021, 301: 202-235.
	     
	    \bibitem{2000Rovira}Rovira C., Tindel S., Sharp Laplace asymptotics for a parabolic SPDE. Stochastics: An International Journal of Probability and Stochastic Processes, 2000, 69(1-2): 11-30.   
		       

		\bibitem{1993Takanobu}Takanobu S.,  Watanabe S.,  Asymptotic expansion formulas of the Schilder-type for a class of conditional Wiener functional integrations. Asymptotic problems in probability theory: Wiener functionals and asymptotics, Proceedings of the Taniguchi International Symposium, Sanda and Kyoto, 1990. Longman Sci. Tech., 1993: 194-241. 
		 
		 

		 

		 
		 
		
		
		
		
				  	 

		 

		 		 
		 		 

		 

		 
		 
		 \bibitem{1987Watanabe}Watanabe S.,  Analysis of Wiener functionals (Malliavin calculus) and its applications to heat kernels. The Annals of Probability, 1987, 15(1): 1-39.
		 

		 
		 
		        



		        

		        
		 
		 
		 

		 
		
		 
		 		 

		 

		 

		

		

		

		
		
		
		
		

		
		
	
			

			

			



		

			

	

		

		

		

		







		

        

        



        

        


















	\end{thebibliography}
\end{document}